\documentclass[12pt,a4paper]{article}
\usepackage{amsmath,amssymb}
\usepackage{latexsym}
\usepackage{graphicx}
\usepackage{color}
\usepackage{multirow}
\usepackage{indentfirst}
\usepackage{subfigure}
\usepackage{mathrsfs}
\usepackage[pdftex]{hyperref}
 \usepackage[numbers,sort&compress]{natbib}

\newtheorem{exam}{\hspace{6mm}Example}[section]

\begin{document}
\baselineskip=2pc

\begin{center}
{\large \bf A quasi-conservative DG-ALE method for multi-component flows using the non-oscillatory kinetic flux}

\end{center}

%\thanks{The reasearch is supported partly by Science Challenge Project (China), No. TZ 2016002 and National Natural Science Foundation-Joint Fund (China) grant U1630247.}
%\footnote{The reasearch is supported partly by Science Challenge Project (China), No. TZ 2016002 and National Natural Science Foundation-Joint Fund (China) grant U1630247.}
\centerline{
Dongmi Luo%
\footnote{Institute of Applied Physics and Computational Mathematics, Beijing 100088, China. E-mail: dongmiluo@stu.xmu.edu.cn.},
Shiyi Li%
\footnote{Institute of Applied Physics and Computational Mathematics, Beijing 100088, China. E-mail: lishiyi14@tsinghua.org.cn.},
Weizhang Huang%
\footnote{Department of Mathematics, University of Kansas, Lawrence, Kansas 66045, U.S.A.
E-mail: whuang@ku.edu.},
Jianxian Qiu%
\footnote{School of Mathematical Sciences and Fujian Provincial
Key Laboratory of Mathematical Modeling and High-Performance
Scientific Computing, Xiamen University, Xiamen, Fujian 361005, China. E-mail: jxqiu@xmu.edu.cn.},
Yibing Chen%
\footnote{Institute of Applied Physics and Computational Mathematics, Beijing 100088, China, 
E-mail: chen\_yibing@iapcm.ac.cn.}
}

\vspace{20pt}

\begin{abstract}
A high-order quasi-conservative discontinuous Galerkin (DG) method is proposed for the numerical simulation of compressible multi-component flows. A distinct feature of the method is a predictor-corrector strategy to define the grid velocity. A Lagrangian mesh is first computed based on the flow velocity and then used as an initial mesh in a moving mesh method (the moving mesh partial differential equation or MMPDE method ) to improve its quality. The fluid dynamic equations are discretized in the direct arbitrary Lagrangian-Eulerian framework using DG elements and the non-oscillatory kinetic flux while the species equation is discretized using a quasi-conservative DG scheme to avoid numerical oscillations near material interfaces. A selection of one- and two-dimensional examples are presented to verify the convergence order and the constant-pressure-velocity preservation property of the method. They also demonstrate that the incorporation of the Lagrangian meshing with the MMPDE moving mesh method works well to concentrate mesh points in regions of shocks and material interfaces.  
\end{abstract}


\pagenumbering{arabic}

\newpage

\section{Introduction}
\label{sec1}
\setcounter{equation}{0}
\setcounter{figure}{0}
\setcounter{table}{0}

The numerical simulation of compressible multi-component flows plays a significant role in the computational fluid dynamics
due to their application in a wide range of fields, including inertial confinement fusion, combustion, and
flame propagation in detonation. Those flows contain rich physical features and structures such as
shock waves, contact waves, and material interfaces. How to resolve these structures accurately and efficiently presents a challenge
in the numerical simulation of compressible multi-component flows. 

There exist two types of frameworks to describe fluid motion, i.e.,
the Eulerian framework and the Lagrangian framework. In the former framework, the mesh is fixed,
which makes methods based on the framework suitable for problems with large deformations.
But contact discontinuities are smeared typically in the Eulerian framework.
On the other hand, methods based on the Lagrangian framework or called Lagrangian methods
(e.g., see \cite{cheng2007,cheng2014,maire2007}), where the mesh moves with the fluid,
are more suitable for flows with contact discontinuities. Unfortunately, the corresponding computation can easily fail due to mesh distortion caused by large flow deformations. To overcome this problem, an arbitrary Lagrangian-Eulerian (ALE) method was first proposed by Hirt et al. \cite{hirt1997} where the grid velocity can be chosen independently from the local flow velocity and therefore the mesh points can be moved to maintain the mesh quality. Roughly speaking, ALE methods can be split into two major groups, indirect ALE methods (e.g., see \cite{maire2008,kucharik2012}) and direct ALE methods (e.g., see \cite{ren2016,boscheri2016, klingenberg2016}). An indirect ALE method typically includes three steps, i.e., the Lagrangian, rezoning, and remapping steps. In the Lagrangian step, the physical variables and the mesh are updated at the same time and the mesh update is based on the flow velocity. The mesh quality is improved using a remeshing strategy in the rezoning step. Finally, the physical solutions are projected from the old mesh to the new one in the remapping step. In a direct ALE method, the mesh movement is taken into account directly in the discretization of the governing equations and the remapping step is avoided.
More high-order methods have been developed with the direct ALE approach since it avoids the remapping step
and the need for high-order conservative remapping schemes.
For example, a high-order moving mesh kinetic scheme has been proposed \cite{xu2013} for compressible flows on unstructured meshes based
on the direct ALE approach. Moreover, a new class of high-order ALE one-step finite volume schemes has been developed in \cite{boscheri2016}
for two- and three-dimensional hyperbolic partial differential equations (PDEs).

Discontinuous Galerkin (DG) methods have been incorporated into the ALE framework; e.g, see \cite{lomtev1999,nguyen2010,persson2009,pandare2016,ren2016,klingenberg2016,zhao2020}. 
DG methods, first applied to neutron transport equations \cite{reed1973} and later extended to general nonlinear systems of hyperbolic conservation laws in a series of papers by Cockburn and Shu and coworkers \cite{cockburn1989, cockburn1990,cockburn1989jcp,cockburn1998}, have been successfully used for multi-material flows \cite{qiu2007, qiu2008, zhu2011,wang2010,henrydefrahan2015, saleem2018,cheng2020,luo2020}.
They are known to be high-order and compact and can readily be combined with various mesh adaptation strategies.
Examples of DG-ALE methods for compressible multiphase flows include that by Pandare et al. \cite{pandare2016} where
the HLLC numerical flux is employed but a constant or a linearly varying grid velocity is assumed in the computation.
Recently, an ALE Runge-Kutta DG method was proposed by Zhao et al. \cite{zhao2020} for multi-material flows where a conservative equation related to
a $\gamma$-model is coupled with the system of fluid equations and the grid velocity is obtained by a variational approach.

In the last four decades a variety of $r$-adaptive or (adaptive) moving mesh methods have been developed 
in the context of the numerical solution of general PDEs; e.g., see the books/review articles \cite{Bai94a,Baines-2011,BHR09,huang2011,Tan05}
and references therein. Moving mesh methods share many common features with ALE methods and in particular,
their rezoning version and quasi-Lagrangian version can be compared directly with indirect and direct ALE methods, respectively;
see discussion in \cite{huang2011}.
The main difference between moving mesh methods and ALE methods lies in how adaptive meshes are generated.
As mentioned earlier, an adaptive mesh is obtained in ALE methods through the flow velocity (Lagrangian meshing) and a remeshing strategy.
This approach has the advantage that mesh points (viewed as fluid particles) are naturally moved to and
concentrated in regions of shocks and contact/material interfaces. Its main disadvantage is that a Lagrangian mesh can often be
distorted and even tangling and thus cannot be used directly in the numerical solution of PDEs. This is also the reason why a remeshing strategy must be used
to effectively improve the quality of the mesh while maintaining a level of mesh concentration in the Lagrangian mesh.
On the other hand, an adaptive mesh is generated in moving mesh methods by either integrating a set of moving mesh equations or minimizing
a meshing functional. This approach has the advantage that it generally generates a mesh of good quality.
The mesh adaptation in this approach is typically associated with an error estimate instead of the flow velocity,
which can be viewed as a drawback of the approach when applied to the computational fluid dynamics.

The objective of the current work is to study a high-order DG-ALE method for multi-component flows
by incorporating the Lagrangian meshing with a moving mesh method. The incorporation is realized with a predictor-corrector strategy, i.e.,
a Lagrangian mesh is first computed and then used as an initial mesh for a moving mesh method to improve mesh quality.
We use the moving mesh PDE (MMPDE) method \cite{huang1994jcp,huang1994siam,huang2011,luo2019} for this purpose.
The MMPDE method moves a mesh through a set of ordinary differential equations (cf. (\ref{MM})) while concentrating points
according to a user-supplied metric tensor that provides the magnitude and directional information to control
the size, shape, and orientation of mesh elements throughout the physical domain.
With a proper choice of the metric tensor,  the MMPDE method is expected to maintain a similar level of point concentration
as in the Lagrangian mesh while improving the mesh quality.
Moreover, it has been shown analytically and numerically in \cite{HK2018} that a formulation of the MMPDE method, at semi-discrete or
fully discrete level, leads to a non-singular (tangling-free) mesh for any domain in any dimension
when the metric tensor is bounded and the initial mesh is non-singular.
This provides a justification on the use of the MMPDE method to improve the quality of the Lagrangian mesh.
With the incorporation of the Lagrangian meshing and the MMPDE method, we hope that the proposed DG-ALE method can have the advantages
of both methods and is able to track rapid changes in the flow field including shocks and material interfaces.
It is noted that this computation of the grid velocity is different from those in \cite{pandare2016,zhao2020}.
(The mesh is assumed to move at a constant or a linearly varying mesh velocity in \cite{pandare2016}
while the grid velocity is obtained by a variational approach in \cite{zhao2020}.)
We adopt the non-oscillatory kinetic (NOK) flux \cite{chen2011,liu2017,luo2020} for the DG discretization of the Euler equations
since  it provides more physical information of the flow (than other numerical fluxes such as the commonly used Lax-Friedrich flux)
and avoids the need to construct any Riemann solver. Moreover, we use a quasi-conservative DG discretization \cite{abgrall1996}
for the species equation to avoid numerical oscillations near material interfaces.

The organization of the paper is as follows. The DG method for conservative systems in the ALE framework and the discretization of the species equation are described in Section~\ref{secmethod}. In Section~\ref{comvel}, the computation of the grid velocity and the MMPDE moving mesh method are discussed.
One- and two-dimensional numerical examples are presented to demonstrate the accuracy and the mesh adaptation capability of the proposed DG-ALE method in Section~\ref{secnum}. Conclusions are drawn in Section~\ref{seccon}.

%% section 2
\section{A DG-ALE method for multi-component flows}
\label{secmethod}
\setcounter{equation}{0}
\setcounter{figure}{0}
\setcounter{table}{0}

\subsection{The govorning equations}
In this work we consider the 4-equation model for a compressible multi-component flow in the dimensionless form
\begin{equation}
\label{eq1}
\begin{cases}
W_t+\nabla\cdot F(W)=0,\\
\frac{\partial Y}{\partial t}+U\frac{\partial Y}{\partial x}+V\frac{\partial Y}{\partial y}=0,
\end{cases}
\end{equation}
where
\[
W=\begin{pmatrix} \rho\\ \rho U \\ \rho V\\ E\end{pmatrix},\quad
F(W)=\begin{pmatrix} \rho U, & \rho V \\ \rho U^2+P, & \rho UV \\ \rho UV, & \rho V^2+P \\ U(E+P), & V(E+P) \end{pmatrix} ,
\]
$\rho$ is the density, $P$ is the pressure, $\vec{U}=(U,V)$ is the flow velocity,
$E=\rho e+\frac{1}{2}\rho(U^2+V^2)$ is the total energy, $\rho e$ is the internal energy, and $Y$ is the volume fraction.
We consider the stiffened gas equation of state (EOS)
\begin{equation}
\rho e=\frac{P+\gamma B}{\gamma-1},
\label{eos-1}
\end{equation}
where $\gamma$ is the ratio of specific heats and $B$ is a prescribed pressure-like constant. 
In this work we consider a two-component flow and assume that the material parameters for fluid 1 and fluid 2
are $\gamma_1,\, B_1$ and $\gamma_2,\, B_2$, respectively. The material parameters for mixing cells are computed as
(e.g., see \cite{shyue1998})
\begin{align}
\label{kappa-1}
\frac{1}{\gamma-1}&=\frac{Y}{\gamma_1-1}+\frac{1-Y}{\gamma_2-1},\\
\label{chi-1}
\frac{\gamma B}{\gamma-1}&=\frac{Y\gamma_1 B_1}{\gamma_1-1}+\frac{(1-Y)\gamma_2B_2}{\gamma_2-1}.
\end{align}
Finally, the mixing sound speed is computed as
\[
c=\sqrt{\frac{\gamma(P+B)}{\rho}}.
\]

\subsection{DG discretization of the fluid dynamic equations on moving meshes}
\label{SEC:NOK}

We now consider the DG discretization of the first four equations in (\ref{eq1}) on a moving triangular mesh
with an arbitrary but given grid velocity. (See Section~\ref{comvel} for the determination of the grid velocity.)
The resulting scheme reduces to the Eulerian and Lagrangian forms when the grid velocity is taken to be zero and the flow velocity, respectively.

For the moment we assume that a moving triangular mesh $\mathscr{T}_h$ for the domain $\Omega$ is given at time instants
\[
0=t_0<t_1<\cdots<t_n<t_{n+1}<\cdots\leqslant T.
\]
The appearances are denoted by $\mathscr{T}_h^n,\; n=0,1,\cdots$. Since they belong to the same mesh, they have the same number of elements ($N$) and vertices ($N_v$) and the same connectivity, and differ only in the location of the vertices. Denote the coordinates of the vertex of $\mathscr{T}_h^n$ by $\textbf{x}_j^n,\; j=1,2,\cdots,N_v$. For any $t \in [t_n, t_{n+1}]$, the coordinates and velocities of the vertices of the mesh $\mathscr{T}_h(t)$ are defined as
\begin{align}
\textbf{x}_j(t)&=\frac{t-t_n}{\Delta t_n}\textbf{x}_j^{n+1}+\frac{t_{n+1}-t}{\Delta t_n}\textbf{x}_j^{n},\qquad j=1,2,\cdots,N_v\\
\label{gridve}
\dot{\textbf{x}}_j(t)&=\frac{\textbf{x}_j^{n+1}-\textbf{x}_j^{n}}{\Delta t_n},\qquad j=1,2,\cdots,N_v
\end{align}
where $\Delta t_n=t_{n+1}-t_{n}$. Then the DG finite element space is defined as
\[
V_h^k=\left \{p(\textbf{x},t):p|_K\in P^k(K), \forall K\in \mathscr{T}_h(t)\right \},
\]
where $P^k(K)$ is the set of polynomials defined on $K$ of degree no more than $k$.
%{\color{blue}Notice that $P^k(K)$ can be expressed as
%\[
%P^k(K)=\text{span}\{\phi_1,\cdots,\phi_L\},
%\]
%where $\phi_l = \phi_l(\textbf{x},t)$, $l = 1, ..., L$ are orthogonal basis functions and $L=k+1$ for the one-dimensional case and $L=\frac{(k+1)(k+2)}{2}$
%for the two-dimensional case. Moreover, the time dependence of the basis functions comes from the time dependence of the location
%of the vertices. Then the DG approximation of $W$ can be expressed as
%\begin{align}
%\label{ns}
%W_h(\textbf{x},t)=\sum\limits_{l=1}^LW_l^K(t)\phi_l(\textbf{x},t),\quad \forall \textbf{x}\in K, \quad \forall K\in \mathscr{T}_h(t) .
%\end{align}
%where $W_l^K(t)=(\rho_l^K(t),(\rho U)_l^K,(\rho V)_l^K,E_l^K)^T.$
%}
The semi-discrete DG approximation of (\ref{eq1}) is to find $W_h(\cdot,t)\in V_h^k, t\in (0,T]$ such that
\begin{equation}
\label{semieq}
\int_K\frac{\partial W_h}{\partial t}\psi dxdy +\int_{\partial K} F\cdot \vec{n}\psi ds-\int_K F\cdot\nabla\psi dxdy=0,\quad \forall \psi\in P^k(K),
\quad \forall K\in \mathscr{T}_h(t)
\end{equation}
where $\vec{n}=(n_x,n_y)$ is the outward unit normal vector of the triangular boundary $\partial K$. 
It is not difficult to show \cite{jimack1991} that 
\[
\frac{\partial \psi(\textbf{x},t)}{\partial t}=-\nabla\psi(\textbf{x},t)\cdot \dot{X}(\textbf{x},t),\quad \forall \psi\in P^k(K)
\]
where $\dot{X}(\textbf{x},t)$ is the piecewise linear interpolant of the nodal velocities.
From this and the Reynolds transport theorem, we have
\begin{align*}
\frac{d}{dt}\int_KW_h\psi dxdy&=\int_K\frac{\partial (W_h\psi )}{\partial t}dxdy+\int_{\partial K}W_h\dot{X}\cdot\vec{n}\psi ds\\
& =\int_K\frac{\partial W_h}{\partial t}\psi dxdy+\int_K\frac{\partial \psi}{\partial t}W_hdxdy+\int_{\partial K}W_h\dot{X}\cdot\vec{n}\psi ds\\
&=\int_K\frac{\partial W_h}{\partial t}\psi dxdy-\int_KW_h\dot{X}\cdot \nabla \psi dxdy+\int_{\partial K}W_h\dot{X}\cdot\vec{n}\psi ds ,
\end{align*}
which yields
\[
\int_K\frac{\partial W_h}{\partial t}\psi dxdy=\frac{d}{dt}\int_KW_h\psi dxdy+\int_KW_h\dot{X}\cdot \nabla \psi dxdy
-\int_{\partial K}W_h\dot{X}\cdot\vec{n}\psi ds.
\]
Inserting this to \eqref{semieq} yields
\begin{equation}
\label{intesemi}
\frac{d}{dt}\int_KW_h\psi dxdy+\int_{\partial K}H(W_h)\psi ds-\int_KH_1(W_h)dxdy=0,
\end{equation}
 where
 \begin{align*}
H(W_h) & = (F(W_h)-W_h\dot{X})\cdot\vec{n}=\begin{pmatrix}
\rho (\vec{U}-\dot{X})\cdot \vec{n}\\
\rho U(\vec{U}-\dot{X})\cdot\vec{n}+Pn_x\\
\rho V(\vec{U}-\dot{X})\cdot\vec{n}+Pn_y\\
E(\vec{U}-\dot{X})\cdot\vec{n}+P\vec{U}\cdot\vec{n}
\end{pmatrix},
\\
H_1(W_h)&=(F(W_h)-W_h\dot{X})\cdot \nabla\psi=\begin{pmatrix}
\rho (\vec{U}-\dot{X}) \cdot \nabla\psi\\
\rho U(\vec{U}-\dot{X}) \cdot \nabla\psi+P\psi_x\\
\rho V(\vec{U}-\dot{X}) \cdot \nabla\psi+P\psi_y\\
E(\vec{U}-\dot{X}) \cdot \nabla\psi+P\vec{U}\cdot \nabla\psi
\end{pmatrix}.
\end{align*}
Recall that these expressions reduce to the Eulerian form when $\dot{X}=0$ and the Lagrangian form when $\dot{X}=\vec{U}$.
 
 Applying a Gaussian quadrature rule to the second and third terms, we get
\begin{align*}
  \int_{\partial K} H(W_h)\psi ds &\approx \sum\limits_{e}
 \sum\limits_{G_e} H(W_h(\textbf{x}_{G_e}))\psi(\textbf{x}_{G_e})w_{G_e}|e|,  \\
 \int_{K} H_1(W_h)dxdy &\approx \sum\limits_{G}H_1(W_h(\textbf{x}_{G}))w_G|K|,
 \end {align*}
where $e$ denotes an edge of the element $K$,  $|K|$ is the volume of $K$, $|e|$ is the length of $e$,
$\textbf{x}_{G}$ and $\textbf{x}_{G_e}$ denote the Gauss points on $K$ and $e$, respectively,
and the summations $\sum_{e}$, $\sum_{G}$, and $\sum_{G_e}$ are taken over the edges of $\partial K$,
Gauss points on $K$, and Gauss points on $e$, respectively.
Replacing the flux $H$ by a numerical flux $\Hat H$, we obtain
\begin{align}
%\begin{cases}
& \frac{d}{dt}\int_{K} W_h \psi dxdy+\sum\limits_e\sum\limits_{G_e} \Hat H\psi(\textbf{x}_{G_e}^{int})w_{G_e}|e| 
-\sum\limits_{G}H_1(W_h(\textbf{x}_{G}))w_G|K|=0,
\notag \\
& \qquad \qquad \qquad. \qquad \qquad \qquad \qquad \qquad \qquad \qquad \qquad \qquad
 \forall K\in \mathscr{T}_h,\;\psi\in V_h^k .
\label{vv}
%\int_K (W_h(\textbf{x},0)-W_0(\textbf{x}))\psi dxdy=0,\;\forall K\in \mathscr{T}_h,\;\psi\in V_h^k
%\end{cases}
\end{align}
Generally speaking, a numerical flux has the form $\hat H=\Hat H(W_h(\textbf{x}_G^{int}),W_h(\textbf{x}_G^{ext}))$,
where $W_h(\textbf{x}_G^{int})$ and $W_h(\textbf{x}_G^{ext})$ are defined as the values of $W_h$ at $\textbf{x} = \textbf{x}_G$
from the interior and exterior of $K$, respectively.

In this work, we use the non-oscillatory kinetic (NOK) flux \cite{chen2009,chen2011,luo2020}
that is a generalization of the gas kinetic scheme (GKS) \cite{xu1997} that avoids the need to construct
Riemann solvers for the numerical solution of hyperbolic PDEs. Different from traditional Godunov-type methods,
GKS is based on the Boltzmann equation and describes the flux function of the governing equations
by particle collision of the transport process. Once the particle distribution function
on a cell interface is obtained, the numerical flux can be calculated directly. The GKS has been successfully applied to multi-component flows \cite{xu1997}.
Combined with a full conservative scheme, the GKS can work well when the difference between two species is small although
oscillations may occur near material interfaces. Chen and Jiang \cite{chen2009} proposed a quasi-conservative scheme
(called the non-oscillatory kinetic scheme or NOK) to overcome the problem and further extended it to more general material problems
in \cite{chen2011}.

The NOK flux is defined using the local coordinates. Denote
 \[
 \tilde{U} = (\vec{U}-\dot{X})\cdot\vec{n},\quad \tilde{V} = (\vec{U}-\dot{X})\cdot\vec{s},
 \quad \tilde{E} = \rho e + \frac{1}{2} \rho ( \tilde{U}^2 +  \tilde{V}^2),
 \]
 where $\vec{s}=(-n_y,n_x)$. Notice that $\tilde{U}$ and $\tilde{V}$ are the components of the vector
 $\vec{U}-\dot{X}$ along the normal and tangential directions of an edge $e$, respectively.
 Denoting $\dot{X}=(U_g,V_g)$, we have
 \begin{align*}
& U=U_g+\tilde{U}n_x-\tilde{V}n_y,\quad V=V_g+\tilde{U}n_y+\tilde{V}n_x ,
 \\
 & U^2+V^2=(U_g^2+V_g^2)+2(U_gn_x+V_gn_y)\tilde{U}+2(V_gn_x-U_gn_y)\tilde{V}+\tilde{U}^2+\tilde{V}^2,
 \\
& E =  \tilde{E}+\frac{1}{2}\rho(U_g^2+V_g^2)+\rho(U_gn_x+V_gn_y)\tilde{U}+\rho(V_gn_x-U_gn_y)\tilde{V}.
\end{align*}
Introducing the new variable
 \begin{equation}
 \xi=
 \begin{pmatrix}
 \xi_1\\
 \xi_2\\
 \xi_3\\
 \xi_4
 \end{pmatrix}
 =
  \begin{pmatrix}
 \rho\tilde{U}\\
 \rho\tilde{U}^2+P\\
 \rho\tilde{U}\tilde{V}\\
 \tilde{U}(\tilde{E}+P)
 \end{pmatrix}
,
 \end{equation}
 we can rewrite the flux along $e$ as
  \begin{equation*}
 \hat{H}(W_h)=
 \begin{pmatrix}
 \xi_1\\
 U_g\xi_1+n_x\xi_2-n_y\xi_3\\
 V_g\xi_1+n_y\xi_2+n_x\xi_3\\
 \frac{1}{2}(U_g^2+V_g^2)\xi_1+(U_gn_x+V_gn_y)\xi_2+(V_gn_x-U_gn_y)\xi_3+\xi_4
 \end{pmatrix} .
 \end{equation*}
This implies that $\hat{H}(W_h)$ can be computed on any edge $e$ if we can compute the variable $\xi$ on the edge.
Moreover, recalling that $\tilde{U}$ and $\tilde{V}$ are the components of the vector
 $\vec{U}-\dot{X}$ along the normal and tangential directions, we know that $\tilde{U}$ has a jump whereas $\tilde{V}$ is continuous
 across the edge. Thus, the computation of $\xi$ on $e$ relies essentially on the computation of $\tilde{U}$ on $e$, which is
 a one-dimensional problem.
 
% \[
% \xi = =  \tilde{U} \begin{pmatrix}
% \rho\\
% \rho\tilde{U}\\
% \rho\tilde{V}\\
% \tilde{E}
% \end{pmatrix}
%  +
%   \begin{pmatrix}
%0\\
%P\\
% 0\\
% \tilde{U} P
% \end{pmatrix}
% ,
%\]

In the NOK flux, $\xi$ is computed as
\[
\xi=\xi_L^++\xi_R^-,
\]
where the subscripts $L$ and $R$ denote the values of the corresponding variable from the interior and exterior of $K$, respectively, and
\begin{align}
\label{NOK-1}
& \xi_L^+=\langle u^1\rangle_+\begin{pmatrix}
\rho\\
\rho\tilde{U}\\
\rho\tilde{V}\\
\tilde{E}
\end{pmatrix}_L+
\begin{pmatrix}
0\\
P_L\langle u^0\rangle_+\\
0\\
\frac{1}{2}P_L\langle u^1\rangle_++\frac{1}{2}P_L\tilde{U}_L\langle u^0\rangle_+
\end{pmatrix},
\\
\label{NOK-2}
& \xi_R^-=\langle u^1\rangle_-\begin{pmatrix}
\rho\\
\rho\tilde{U}\\
\rho\tilde{V}\\
\tilde{E}
\end{pmatrix}_R+
\begin{pmatrix}
0\\
P_R\langle u^0\rangle_-\\
0\\
\frac{1}{2}P_R\langle u^1\rangle_-+\frac{1}{2}P_R\tilde{U}_R\langle u^0\rangle_-
\end{pmatrix},
\\
\label{NOK-3}
&
\langle u^1\rangle_+=\tilde{U}_L\langle u^0\rangle_++\frac{1}{2}\frac{e^{-\lambda\tilde{U}_L^2}}{\sqrt{\pi \lambda}},
\qquad 
\langle u^1\rangle_-=\tilde{U}_R\langle u^0\rangle_--\frac{1}{2}\frac{e^{-\lambda\tilde{U}_R^2}}{\sqrt{\pi \lambda}},
\\
\label{NOK-4}
&
\langle u^0\rangle_+=\frac{1}{2}\text{erfc}(-\sqrt{\lambda}\tilde{U}_L),
\qquad
\langle u^0\rangle_-=\frac{1}{2}\text{erfc}(\sqrt{\lambda}\tilde{U}_R),
\end{align}
$\lambda=\min\{\frac{1}{c_L^2},\frac{1}{c_R^2}\}$, $c$ is sound speed, and $\text{erfc}(\cdot)$ is
the complementary error function. We refer the reader to \cite{chen2011} for the derivation of these formulas.

Finally we remark that the scheme (\ref{vv}) involves the grid velocity $\dot{X}=(U_g,V_g)$ which is defined as
the piecewise linear interpolant of the nodal grid velocity.  The latter is computed as the change rate of the location
of vertices between the old and new meshes $\mathscr{T}_h^{n}$ and $\mathscr{T}_h^{n+1}$; cf. (\ref{gridve}).
The generation of the new mesh will be discussed in Section~\ref{comvel}.

%From the equation \eqref{intesemi} , one can find that the grid velocities are unknown
%in the numerical flux. In the following section, we will describe the computation of the grid velocities.

\subsection{DG-ALE for species equation}
\label{seczf}

In this subsection we consider the DG discretization of the species equation in (\ref{eq1}). Numerical experiments show that
a direct DG discretization of the equation can lead to oscillations in the pressure and velocity near material interfaces.
The analysis by Abgrall \cite{abgrall1996} shows that those oscillations can be prevented if the scheme preserves
constant velocity and constant pressure. We follow the procedure of Abgrall \cite{abgrall1996} and develop such a DG discretization
of the species equation on a moving mesh in the following.

To start with, we assume that the velocity and the pressure are constant, i.e., $\vec{U}=\vec{U}_0=(U_0,V_0)$ and $P=P_0$.
Denote
$
\tilde{\vec{U}}_0=\vec{U}_0-\dot{X} .
$
We can rewrite $H(W_h)$ and $H_1(W_h)$ as
\begin{align}
H(W_h)=\begin{pmatrix}
\rho \tilde{\vec{U}}_0\cdot \vec{n}\\
\rho U_0\tilde{\vec{U}}_0\cdot\vec{n}+P_0n_x\\
\rho V_0\tilde{\vec{U}}_0\cdot\vec{n}+P_0n_y\\
E\tilde{\vec{U}}_0\cdot\vec{n}+P_0\vec{U}_0\cdot\vec{n}
\end{pmatrix},
\quad
H_1(W_h)=\begin{pmatrix}
\rho \tilde{\vec{U}}_0 \cdot \nabla\psi\\
\rho U_0\tilde{\vec{U}}_0 \cdot \nabla\psi+P_0\psi_x\\
\rho V_0\tilde{\vec{U}}_0 \cdot \nabla\psi+P_0\psi_y\\
E\tilde{\vec{U}}_0 \cdot \nabla\psi+P_0\vec{U}_0 \cdot \nabla\psi
\end{pmatrix}.
\end{align}
Then, the continuity equation in \eqref{intesemi} becomes
\begin{align}
\label{con}
\frac{d}{dt}\int_K \rho \psi dx d y =-\int_{\partial K} \rho\tilde{\vec{U}}_0\cdot\vec{n}\psi ds
+\int_K\rho\tilde{\vec{U}}_0 \cdot \nabla\psi dxdy .
\end{align}
Similarly, the momentum equations read as
\begin{align*}
\frac{d}{dt}\int_K \rho U \psi dx d y&=-\int_{\partial K}(\rho U_0\tilde{\vec{U}}_0\cdot\vec{n}+P_0n_x)\psi ds
+\int_K(\rho U_0\tilde{\vec{U}}_0\cdot \nabla\psi +P_0\psi_x)dxdy,\\%\frac{\partial\psi}{\partial x}dxdy,\\
\frac{d}{dt}\int_K \rho V \psi dx d y&=-\int_{\partial K}(\rho V_0\tilde{\vec{U}}_0\cdot\vec{n}+P_0n_y)\psi ds
+\int_K(\rho V_0\tilde{\vec{U}}_0\cdot \nabla\psi +P_0\psi_y)dxdy.%\frac{\partial\psi}{\partial y}dxdy .
\end{align*}
From the divergence theorem, we have $\int_K P_0\nabla\psi dxdy=\int_{\partial K}P_0\psi \vec{n}ds$.
Inserting this into the above equations, we get
\begin{align*}
\frac{d}{dt}\int_K \rho U \psi dx d y&=-\int_{\partial K}\rho U_0\tilde{\vec{U}}_0\cdot\vec{n} \psi ds
+\int_K\rho U_0\tilde{\vec{U}}_0\cdot \nabla\psi dxdy,\\
\frac{d}{dt}\int_K \rho V \psi dx d y&=-\int_{\partial K}\rho V_0\tilde{\vec{U}}_0\cdot\vec{n} \psi ds
+\int_K\rho V_0\tilde{\vec{U}}_0\cdot \nabla\psi dxdy .
\end{align*}
Comparing these with (\ref{con}), we have
\begin{equation}
\frac{d}{dt}\int_K \rho U \psi dx d y = U_0 \frac{d}{dt}\int_K \rho \psi dx d y, \quad
\frac{d}{dt}\int_K \rho V \psi dx d y = V_0 \frac{d}{dt}\int_K \rho \psi dx d y.
\label{moment-1}
\end{equation}
From the Reynolds transport theorem, we have
\begin{align*}
\frac{d}{dt}\int_K \rho U \psi dx d y & = \int_K U \frac{\partial (\rho \psi)}{\partial t} dx d y
+ \int_K \frac{\partial U}{\partial t} \rho \psi  dx d y + \int_{\partial K} U \rho \psi \vec{n} \cdot \dot{X} ds
\\
& = \int_K U_0 \frac{\partial (\rho \psi)}{\partial t} dx d y
+ \int_K \frac{\partial U}{\partial t} \rho \psi  dx d y + \int_{\partial K} U_0 \rho \psi \vec{n} \cdot \dot{X} ds
\\
& = U_0 \frac{d}{dt}\int_K \rho \psi dx d y + \int_K \frac{\partial U}{\partial t} \rho \psi  dx d y .
\end{align*}
Similarly, we have
\[
\frac{d}{dt}\int_K \rho V \psi dx d y = V_0 \frac{d}{dt}\int_K \rho \psi dx d y + \int_K \frac{\partial V}{\partial t} \rho \psi  dx d y .
\]
Combining the above results with (\ref{moment-1}), we get 
\[
\int_K \rho \frac{\partial U}{\partial t} \psi dx d y = 0, \quad
\int_K \rho \frac{\partial V}{\partial t}  \psi dx d y = 0,
\]
which, from the arbitrariness of $\psi$ and the positivity of $\rho$, implies $\frac{\partial U}{\partial t} = 0$ and $\frac{\partial V}{\partial t}=0$.
Thus, the scheme (\ref{intesemi}) preserves constant velocity solutions.

We now consider the preservation of constant pressure. The energy equation in (\ref{intesemi}) becomes
\[
\frac{d}{dt}\int_K E \psi dx d y =-\int_{\partial K}(E \tilde{\vec{U}}_0\cdot\vec{n}+P_0 \vec{U}_0 \cdot \vec{n})\psi ds
+\int_K(E \tilde{\vec{U}}_0\cdot \nabla\psi +P_0 \vec{U}_0 \cdot \nabla \psi )dxdy,
\]
which can be simplified into
\begin{equation}
\frac{d}{dt}\int_K E \psi dx d y =-\int_{\partial K} E \tilde{\vec{U}}_0\cdot\vec{n}\psi ds
+\int_KE \tilde{\vec{U}}_0\cdot \nabla\psi dxdy .
\label{energy-1}
\end{equation}
Let $\kappa=\frac{1}{\gamma-1}$ and $\chi=\frac{\gamma B}{\gamma-1}$. Then we can rewrite the equation of state (\ref{eos-1}) into
$\rho e = \kappa P  + \chi$. From the definition of $E$, we have $E = \kappa P + \chi + \frac{1}{2} \rho (U^2 + V^2)$.
Inserting this into (\ref{energy-1}), we get
\begin{align*}
& \frac{d}{dt}\int_K \kappa P \psi dx d y + \frac{d}{dt}\int_K \chi \psi dx d y + \frac{d}{dt}\int_K \frac{1}{2} \rho (U^2+V^2) \psi dx d y
\\
= & -\int_{\partial K} \kappa P  \tilde{\vec{U}}_0\cdot\vec{n}\psi ds +\int_K\kappa P \tilde{\vec{U}}_0\cdot \nabla\psi dxdy
\\
& -\int_{\partial K} \chi  \tilde{\vec{U}}_0\cdot\vec{n}\psi ds +\int_K\chi \tilde{\vec{U}}_0\cdot \nabla\psi dxdy
\\
& -\int_{\partial K}  \frac{1}{2} \rho (U^2 + V^2) \tilde{\vec{U}}_0\cdot\vec{n}\psi ds +\int_K \frac{1}{2} \rho (U^2 + V^2) \tilde{\vec{U}}_0\cdot \nabla\psi dxdy.
\end{align*}
Using the assumption that $U = U_0$, $V = V_0$, and $P=P_0$, the fact that $\frac{\partial U}{\partial t} = 0$ and $\frac{\partial V}{\partial t}=0$,
and the equation (\ref{con}), we can derive from the above equation that
\[
\int_K \kappa \frac{\partial P}{\partial t} \psi dx d y = 0, \quad \text{or} \quad \frac{\partial P}{\partial t}  = 0,
\]
provided that
\begin{align}
\label{gam}
\frac{d}{dt}\int_K \kappa \psi d x d y&=-\int_{\partial K}\kappa\tilde{\vec{U}}_0\cdot\vec{n}\psi ds
+\int_K\kappa\tilde{\vec{U}}_0\cdot \nabla\psi dxdy,\\
\label{gamp}
\frac{d}{dt}\int_K \chi \psi d x d y &=-\int_{\partial K}\chi\tilde{\vec{U}}_0\cdot\vec{n}\psi ds+\int_K\chi\tilde{\vec{U}}_0\cdot \nabla\psi dxdy .
\end{align}
In other words, the scheme (\ref{intesemi}) also preserves constant pressure if the above two conditions are satisfied.

We notice from (\ref{kappa-1}) and (\ref{chi-1}) that both $\kappa$ and $\chi$ are linear functions of $Y$. Thus, (\ref{gam}) and (\ref{gamp})
can be satisfied if the following equation for $Y$ holds,
\begin{equation}
\frac{d}{dt}\int_K Y_h \psi d x d y =-\int_{\partial K} Y_h \tilde{\vec{U}}_0\cdot\vec{n}\psi ds
+\int_K Y_h \tilde{\vec{U}}_0\cdot \nabla\psi dxdy .
\label{zufen}
\end{equation}
It is not difficult to see that the above equation is actually a DG discretization of the equation
\[
Y_t+(UY)_x+(VY)_y=0.
\]
Since the species equation (cf. (\ref{eq1})) can be written as
\[
Y_t+UY_x+VY_y=Y_t+(UY)_x+(VY)_y-YU_x-YV_y=0,
\]
we obtain a DG discretization for the species equation as
\begin{align}
\label{zufen1}
\frac{d}{dt}\int_K Y_h \psi d x d y=&-\int_{\partial K}Y_h(\vec{U}-\dot{X})\cdot\vec{n}\psi ds
+\int_K Y_h (\vec{U}-\dot{X}) \cdot \nabla\psi dxdy \notag\\
& + Y_h(\textbf{x}_b)\left (\int_{\partial K}\vec{U}\cdot\vec{n}\psi ds-\int_K\vec{U}\cdot \nabla\psi dxdy\right ),
\quad \forall K\in \mathscr{T}_h,\;\psi\in V_h^k
\end{align}
where $\textbf{x}_b$ is the barycenter of the element $K$. Noticing that the above equation reduces to \eqref{zufen}
for a solution with constant velocity and constant pressure, we know that (\ref{zufen1}), along with (\ref{intesemi}) (or (\ref{vv})), preserves constant velocity and constant pressure. Moreover, \eqref{zufen1} is only quasi-conservative since the last term is non-conservative.

 \subsection{Temporal discretization and limiting}
 \label{sectt}

 Finally, the semi-discrete schemes (\ref{vv}) and (\ref{zufen1}) are discretized in time. Here, we use
an explicit, the third order TVD Runge-Kutta scheme \cite{shu1988}. Casting (\ref{vv}) and (\ref{zufen1}) in the form
$$
\frac{d(\textbf{M}\textbf{U}_h)}{dt}=L_h(\textbf{U}_h,t),
$$
where $\textbf {M}$ is the mass matrix. Noticing that the mesh motion of the method is part of Lagrangian solution, we use the same temporal discretization
to update the physical solution, the mesh coordinates, and other geometrical quantities such as the cell volume and the edge length.
The scheme reads as

Stage 1
\begin{align}
\label{tt1}
\textbf{x}_j^{(1)}&=\textbf{x}_j^{n}+\Delta t_n \dot{X}_j,\notag\\
 \textbf{M}^{(1)}\textbf{U}_h^{(1)}&=\textbf{M}^{n}\textbf{U}_h^n+\Delta t_n L_h(\textbf{U}_h^n,t_n).
\end{align}

Stage 2 
\begin{align}
\textbf{x}_j^{(2)}&=\frac{3}{4}\textbf{x}_j^{n}+\frac{1}{4}(\textbf{x}_j^{(1)}+\Delta t_n \dot{X}_j),\notag\\
 \textbf{M}^{(2)}\textbf{U}_h^{(2)}&=\frac{3}{4}\textbf{M}^{n}\textbf{U}_h^n+\frac{1}{4}(\textbf{M}^{(1)}\textbf{U}_h^{(1)}+\Delta t_n L_h(\textbf{U}_h^{(1)},t_n+\Delta t_n)).
\end{align}

Stage 3 
\begin{align}
\textbf{x}_j^{n+1}&=\frac{1}{3}\textbf{x}_j^{n}+\frac{2}{3}(\textbf{x}_j^{(2)}+\Delta t_n \dot{X}_j),\notag\\
 \textbf{M}^{n+1}\textbf{U}_h^{n+1}&=\frac{1}{3}\textbf{M}^{n}\textbf{U}_h^n+\frac{2}{3}(\textbf{M}^{(2)}\textbf{U}_h^{(2)}+\Delta t_n L_h(\textbf{U}_h^{(2)},t_n+\frac{1}{2}\Delta t_n)).
\end{align}

We need to use a nonlinear limiter to control spurious oscillations in the numerical solution when strong shocks are present.
The process typically contains identification of troubled cells and modification of the DG approximations on those cells.
Here, we employ the minmod-type TVB limiter \cite{cockburn1990,cockburn1989,cockburn1989jcp,cockburn1998,qiu2005,luo2019}
to detect troubled cells and use the multi-resolution WENO limiter developed in \cite{zhu2020,zhu2020a} for the limiting process.
Moreover, we implement this limiting process to the primitive variables component-wisely to keep the pressure non-oscillatory.
Since the process is the same for schemes on a fixed mesh, we omit the detail here and refer the reader to \cite{luo2020}.

% section 3
\section{The computation of the grid velocity}
\label{comvel}
\setcounter{equation}{0}
\setcounter{figure}{0}
\setcounter{table}{0}

We consider the computation of the grid velocity in this section.
To take the advantages of both the Lagrangian meshing and the moving mesh method,
we employ a predictor-corrector strategy to define the grid velocity. More specifically, we first compute the Lagrangian velocity
and obtain a Lagrangian mesh (cf. Section~\ref{lagvel}). Then we use the MMPDE method
\cite{huang1994jcp,huang1994siam,huang2011,luo2019} to improve the quality of the Lagrangian mesh.
The MMPDE method is described in Section~\ref{sec:mmpde}.

\subsection{The Lagrangian meshing}
\label{lagvel}

There are many possible ways to determine the vertex velocity. Here, we adopt a simple yet robust least-squares
algorithm \cite{addessio1990,chen2010}. The first step of the algorithm is to find the flow velocity $U_e^*$ for each edge $e$.
This has been done in \cite{addessio1990,chen2010} by solving a Riemann problem in the normal direction of the edge.
Here we propose to use the idea of the NOK flux (cf. Section~\ref{SEC:NOK}) to compute $U_e^*$. Specifically, from (\ref{NOK-3}) and (\ref{NOK-4})
(and taking $\dot{X} = 0$) we  compute $U_e^*$ as
\begin{equation}
U_e^*=\langle u^1\rangle_++\langle u^1\rangle_-,
\end{equation}
where 
\begin{align}
\label{NOK-5}
&
\langle u^1\rangle_+=(\vec{U}_L\cdot \vec{n}) \langle u^0\rangle_++\frac{1}{2}\frac{e^{-\lambda (\vec{U}_L\cdot \vec{n})^2}}{\sqrt{\pi \lambda}},
\qquad 
\langle u^1\rangle_-=(\vec{U}_R\cdot \vec{n})\langle u^0\rangle_--\frac{1}{2}\frac{e^{-\lambda(\vec{U}_R\cdot \vec{n})^2}}{\sqrt{\pi \lambda}},
\\
\label{NOK-6}
&
\langle u^0\rangle_+=\frac{1}{2}\text{erfc}(-\sqrt{\lambda}(\vec{U}_L\cdot \vec{n})),
\qquad \qquad \quad \;
\langle u^0\rangle_-=\frac{1}{2}\text{erfc}(\sqrt{\lambda}(\vec{U}_R\cdot \vec{n})),
\end{align}
$\lambda=\min\{\frac{1}{c_L^2},\frac{1}{c_R^2}\}$, $c$ is the sound speed, and $\vec{U}_R$ and $\vec{U}_L$ denote the values of $\vec{U}$
from the positive-$\vec{n}$ side and the negative-$\vec{n}$ side of $e$, respectively.

The second step is to obtain the velocity at each vertex
using the least squares method. Specifically, for a given vertex $j$, we would like the normal projection of the vertex velocity
on any of the connecting edges to be equal to the corresponding Riemann velocity, viz.,
\[
\dot{X}_j^L \cdot\vec{n}_{e_i}=U_{e_i}^*, \quad i=1,2,\cdots,n_j
\]
where $n_j$ denotes the number of the edges connected to the vertex $j$. The system can be overdetermined or underdetermined
and thus we solve it by minimizing a weighted sum of the squares of the residuals with respect to $\dot{X}_j^L$,
\[
\sum\limits_i\alpha_{e_i} (\dot{X}_j^L\cdot \vec{n}_{e_i}-U_{e_i}^*)^2,
\]
where $\alpha_{e_i}$ is a weight taken as the mean of the density on each side of the edge in our computation.
Having computed the Lagrangian velocity for all vertices (with proper modifications for boundary vertices so that they stay on the domain),
we can obtain the Lagrangian mesh $\mathscr{T}_h^{n+1,L}$ by updating the location of the vertices as
\[
\textbf{x}_j^L=\textbf{x}_j^n+\Delta t_n\cdot\dot{X}_j^L, \quad j = 1, ..., N_v.
\]

\subsection{The MMPDE moving mesh method}
\label{sec:mmpde}

Generally speaking, a Lagrangian mesh can be severely distorted and thus cannot be used
directly for the numerical solution of PDEs. In this subsection we describe the MMPDE method that is employed
to improve the quality of Lagrangian meshes. 

To start with, we assume a reference computational mesh $\Hat {\mathscr{T}}_c = \{\Hat {\pmb{\xi}}_j\}^{N_v}_{j=1}$
has been chosen, which is a deformation of the physical mesh and is fixed for the whole computation.
Usually we should choose it to be as uniform as possible and can take it as the initial physical mesh.
We also need to use the computational mesh ${\mathscr{T}}_c = \{ {\pmb{\xi}}_j\}^{N_v}_{j=1}$,
which is also a deformation of the physical mesh and will be used as an intermediate variable.

The MMPDE method views any nonuniform mesh as a uniform one in some metric specified by a metric tensor
$\mathbb{M}=\mathbb{M}(\textbf{x})$ that is a symmetric and positive definite matrix for each $\textbf{x}$ and
uniformly positive definite on $\Omega$. Such a uniform mesh $\mathscr{T}_h$ in $\mathbb{M}$ satisfies the so-called equidistribution
and alignment conditions \cite{huang2006}
\begin{align}
& |K|\sqrt{\hbox{det}(\mathbb{M}_K)}=\frac{\sigma_h|K_c|}{|\Omega_c|},\qquad \forall K\in \mathscr{T}_h
\label{eq}
\\
& \frac{1}{d}\text{tr}((F'_K)^{-1}\mathbb{M}_K^{-1}(F'_K)^{-T})=\hbox{det}((F'_K)^{-1}\mathbb{M}_K^{-1}(F'_K)^{-T})^{\frac{1}{d}},
\qquad \forall K\in\mathscr{T}_h
\label{al}
\end{align}
where $d$ is the dimension of the spatial domain ($d=1$ for one dimension and $d=2$ for two dimensions),
$K_c$ is the element in $ \mathscr{T}_c$ corresponding to $K$,
$F_K$ is the affine mapping from $K_c$ to $K$ and $F_K'$ is its Jacobian matrix,
$\mathbb{M}_K$ is the average of $\mathbb{M}$ over $K$, $\hbox{tr}(\cdot)$ and $\hbox{det}(\cdot)$ denote
the trace and determinant of a matrix, respectively, and
\[
|\Omega_c|=\sum\limits_{K_c\in\mathscr{T}_c}|K_c|,\quad
\sigma_h =\sum\limits_{K\in\mathscr{T}_h}|K|\hbox{det}(\mathbb{M}_K)^{\frac{1}{2}} .
\]
The equidistribution condition \eqref{eq} requires that all of the elements have the same size. It determines the size of elements
through the metric tensor $\mathbb{M}$: $|K|$ is smaller when $\text{det}(\mathbb{M}_K)$ is larger and vice versa.
On the other hand, the alignment condition \eqref{al} requires that element $K$, measured in the metric $\mathbb{M}_K$, be similar
to $K_c$. It determines the shape and orientation of $K$ through the combination of $\mathbb{M}_K$ and the shape and orientation of $K_c$.
In the case where $K_c$ is equilateral, the shape and orientation of $K$ is then determined completely by $\mathbb{M}_K$
(through its eigenvalues and eigen-directions).

The objective of the MMPDE method is to generate a mesh satisfying the above two conditions as closely as possible.
This is done by minimizing the energy function
    \begin{align}
    \label{fl}
    I_h(\mathscr{T}_h,\mathscr{T}_c) &= \sum\limits_{K\in\mathscr{T}_h}|K|\sqrt{\hbox{det}(\mathbb{M}_K)}(\hbox{tr}((F'_K)^{-1}\mathbb{M}^{-1}_K(F'_K)^{-T}))^{\frac{3 d}{4}} \notag\\
     &\qquad + d^{\frac{3 d}{4}}\sum\limits_{K\in\mathscr{T}_h} |K|\sqrt{\hbox{det}(\mathbb{M}_K)}\left (\frac{|K_c|}{|K|\sqrt{\hbox{det}(\mathbb{M}_K)}}\right )^{\frac{3}{2}},
    \end{align}
which is a Riemann sum of a continuous functional \cite{huang2001} based on equidistribution and alignment for variational mesh adaptation.
Here we use the $\xi$-formulation where $\mathscr{T}_h$ is taken as a currently available physical mesh (denoted by $\tilde{\mathscr{T}}_h$,
its choice will be discussed later)
and $I_h(\tilde{\mathscr{T}}_h,\mathscr{T}_c)$
is minimized with respect to $\{\pmb{\xi}_j\}$. The minimization is realized by solving the gradient system (i.e., the mesh equation)
\begin{align}
\label{MM}
\frac{d\pmb{\xi}_j}{dt}=-\frac{\text{det}(\mathbb{M}(\textbf{x}_j))^{\frac{1}{4}}}{\tau}\left ( \frac{\partial I_h}{\partial \pmb{\xi}_j}\right )^T,\quad j=1,2,\cdots,N_v
\end {align}
where $\frac{\partial I_h}{\partial \pmb{\xi}_j}$ is considered as a row vector and $\tau>0$ is a parameter used
to adjust the response time of mesh movement to changes in $\mathbb{M}$.
The analytical formula for $\frac{\partial I_h}{\partial \pmb{\xi}_j}$ can be found in \cite{HK2015}.

The mesh equation (\ref{MM}) (with proper modifications for boundary vertices) can be integrated from $t_n$ to ${ t_{n+1}}$, starting  with the reference computational mesh $\Hat {\mathscr{T}}_c$ as an initial mesh, to obtain a new computational mesh $\mathscr{T}_c^{n+1}$. This new mesh forms
a correspondence $\Psi_h$ with the physical mesh $\tilde{\mathscr{T}}_h$ having the property $\tilde{\textbf{x}}_j=\Psi_h(\pmb{\xi}_j^{n+1}),\; j=1,2,\cdots,N_v$.
The new physical mesh $\mathscr{T}_h^{n+1}$ is then defined as $\textbf{x}_j^{n+1}=\Psi_h(\Hat {\pmb{\xi}}_j),\; j=1,2,\cdots,N_v$,
which can be computed using linear interpolation.

In the context of the MMPDE method, it is common to take $\tilde{\mathscr{T}}_h$ as
the physical mesh at the current step, i.e., $\tilde{\mathscr{T}}_h = \mathscr{T}_h^{n}$. In our current situation, however, we choose $\tilde{\mathscr{T}}_h
= \mathscr{T}_h^{n+1,L}$ to take the advantages of the Lagrangian mesh $\mathscr{T}_h^{n+1,L}$.
In this setting, the MMPDE method plays a role of improving the quality of the Lagrangian mesh.
It has been shown in \cite{HK2018} analytically and numerically that
the $x$-formulation of the MMPDE method, where $I_h(\mathscr{T}_h,\hat{\mathscr{T}}_c)$ (with $\mathscr{T}_c = \hat{\mathscr{T}}_c$)
is minimized by solving its gradient system with respect to $\{\textbf{x}_j\}$, produces non-singular (tangling-free) meshes
when the metric tensor is bounded and the initial mesh is non-singular.
Although the $\xi$-formulation (that we use in this work) is different from the $x$-formulation,
they are numerical approximations to the same continuous optimization problem and thus we can expect the $\xi$-formulation also leads
to non-singular meshes especially when the mesh is sufficiently fine. Indeed, our experience shows that the $\xi$-formulation is
as robust as the $x$-formulation in terms of generating non-singular meshes.
Moreover, with a proper metric tensor (see its definition below), the MMPDE method can be expected to concentrate mesh points
near regions of shocks and material interfaces. Thus, we hope that the use of the MMPDE method to improve quality of the Lagrangian mesh
will not (at least not significantly) move points away from shocks and material interfaces. In this sense, the MMPDE maintains
a level of mesh concentration of the Lagrangian mesh.

Having obtained the new physical mesh $\mathscr{T}_h^{n+1}$, we can compute the nodal grid velocity using (\ref{gridve}) and then
$\dot{X}=(U_g,V_g)$ as the piecewise linear interpolant of the nodal grid velocity.

A key to the success for the MMPDE moving mesh method is the selection of a proper metric tension $\mathbb{M}$ that
provides the information needed to control the size, shape and orientation of the mesh elements throughout the domain.
In this work, we employ a metric tensor that is known to be optimal for the $L^2$ norm of linear interpolation error \cite{huang2003}
and has been used for conservation laws in \cite{luo2019}, i.e.,
\begin{align}
\label{mer}
\mathbb{M}=\hbox{det}(\mathbb{I}+|H(u_h^n)|)^{-\frac{1}{d+4}}(\mathbb{I}+|H(u^n_h)|),
\end {align}
where $u$ denotes a physical variable,
$\mathbb{I}$ is the $d\times d$ identity matrix, $H(u_h^n)$ is a recovered Hessian from the numerical solution $u^n_h$,
and $|H(u^n_h)|=Q\hbox{diag}(|\lambda_1|,\cdots,|\lambda_d|)Q^T$ with $Q\hbox{diag}(\lambda_1,\cdots,\lambda_d)Q^T$
being the eigen-decomposition of $H(u^n_h)$. The Hessian is recovered using the quadratic least squares fitting \cite{zhang2005}.
For a smoother mesh, a low-pass-filter smoothing algorithm \cite{huang2011} is applied to the metric tensor several sweeps every time it is computed.

For multi-component flows, we take $u$ to be the quantity
\begin{align}
\label{ent}
S=1+\beta_1\left (\frac{\rho}{\| \rho\|_\infty}\right )^2+\beta_2\left (\frac{P}{\| P\|_\infty}\right )^2+\beta_3\left (\frac{Y}{\| Y\|_\infty}\right )^2,
\end{align}
where $\beta_i, i=1,2,3$ are positive parameters and $\| \cdot \|_\infty$ is the maximum norm.
In the computation of the metric tensors, the nodal values of $\rho$, $P$, and $Y$ are needed. They are calculated as the volume
average of the cell average of these variables in the neighboring cells of each vertex.

%at $\textbf{x}_j\; (j=1,2,\cdots,N_v)$ are computed as
%\begin{align}
%\label{ent}
%%&S_j=1+\beta_1(\frac{\rho_j}{\max\limits_{1\leqslant m\leqslant N_v}(\rho_m)})^2
%%+\beta_2(\frac{P_j}{\max\limits_{1\leqslant m\leqslant N_v}(P_m)})^2+\beta_3(\frac{Y_j}{\max\limits_{1\leqslant m\leqslant N_v}(Y_m)})^2,\\
%&{\rho}_j=\frac{\sum\limits_{K\in\omega_j}|K|{\rho}_K}{\sum\limits_{K\in\omega_j}|K|},\qquad 
% {P_j}=\frac{\sum\limits_{K\in\omega_j}|K|{P}_K}{\sum\limits_{K\in\omega_j}|K|},\notag\qquad
% {Y_j}=\frac{\sum\limits_{K\in\omega_j}|K|{Y}_K}{\sum\limits_{K\in\omega_j}|K|},\notag
%\end{align}
%where $\omega_j$ is the element patch associated with $\textbf{x}_j$. 

% section 4
\section{Numerical examples}
\label{secnum}
\setcounter{equation}{0}
\setcounter{figure}{0}
\setcounter{table}{0}

In this section we present numerical results obtained with the DG-ALE method described
in the previous sections for a selection of one- and two-dimensional examples.
Recall that the method has been described in two dimensions. Its implementation in one dimension
is similar. The CFL number in time step selection is set to be 0.3 for $P^1$ elements, 0.15 for
$P^2$ elements. 
The parameter $\tau$ in \eqref{MM} is taken as $0.1$ for accuracy test problems,
$10^{-3}$ and $10^{-4}$ for one- and two-dimensional { systems} with discontinuities, respectively.
The parameter $\beta_i$'s in \eqref{ent} are taken as 1.0, 1.0, and 0 for all the computations unless otherwise stated.
The results obtained with the DG-ALE method will be marked as ``ALEMM" in figures.
%Periodic boundary conditions are used for all but Examples~\ref{examgasgas} and \ref{exambw}.

\subsection{One-dimensional examples}

\begin{exam}{\em
\label{examorder}
To verify the convergence order of the new method, we first consider a sine wave composed of two components with
$\gamma_1=1.4$, $\gamma_2=1.9$, $B_1=1$, and $B_2=0$.
The initial condition is given by
\begin{align*}
\rho(x,0)=1+0.2\sin(\pi x),\quad v(x,0)=1,\quad P(x,0)=1,\quad Y(x,0)=0.5+0.5\sin(\pi x).
\end{align*}
A periodic boundary condition is used. We take the computational domain as (0, 2) and compute the solution up to $t = 0.5$. 
 
 The error of the density is listed in Table \ref{exorder}, which shows the second-order convergence for $P^1$ elements
 and the third-order convergence for $P^2$ elements for the DG-ALE method.

%  \begin{table}
%\caption{Example~\ref{examorder}: The solution error in the density at $t=0.5$.}
%\renewcommand{\multirowsetup}{\centering}
%\begin{center}
%\begin{tabular}{|c|c|c|c|c|c|c|c|c|c|c|c|c|}
%\hline
%$k$ & $N$   &40           &80            & 160 & 320 & 640 &1280 \\
%\hline
%\multirow{6}{1cm}{1}
%%\hline
% & $L^1$   &1.141e-4      & 2.851e-5 & 7.186e-6 & 1.820e-6  & 4.613e-7 & 1.165e-7\\
%%\hline
% &  Order      & \quad  & 2.001       & 1.988      &  1.981       & 1.980 &1.984 \\
% %\hline
% &$L^2$   & 1.469e-4      & 3.684e-5 & 9.354e-6 &  2.391e-6  & 6.110e-7 & 1.552e-7\\
%%\hline
% &Order   & \quad   & 1.995     & 1.978      &   1.968      & 1.969 & 1.977 \\
%%\hline
% &$L_{\infty}$ & 4.451e-4 & 1.105e-4 & 2.888e-5  &  7.539e-6 & 1.931e-6 & 4.881e-7\\
%%\hline
% &Order    & \quad  & 2.010      & 1.936     &   1.938      & 1.965 & 1.984   \\
% \hline
% \multirow{6}{1cm}{2}
%
% & $L^1$  & 9.594e-6      &1.255e-6 & 1.525e-7 &  1.810e-8 & 2.166e-9 & 2.634e-10 \\
%%\hline
% &  Order      & \quad & 2.934    & 3.041     & 3.075    & 3.063 & 3.039  \\
%% \hline
% &$L^2$  & 1.361e-5       & 1.827e-6 & 2.228e-7 &  2.598e-8 & 3.043e-9 & 3.640e-10\\
%%\hline
% &Order       & \quad &2.897     & 3.036       & 3.100    & 3.094 & 3.063 \\
%%\hline
% &$L_{\infty}$ & 5.101e-5 & 7.733e-6 & 1.021e-6  & 1.186e-7 & 1.332e-8 & 1.491e-9 \\
%%\hline
% &Order    & \quad       &2.722    & 2.921     & 3.105   & 3.155 &3.160   \\
% \hline 
%\end{tabular}%\\
%\end{center}
%\label{exorder}
%\end{table}

 \begin{table}
\caption{Example~\ref{examorder}: The solution error in the density at $t=0.5$.}
\renewcommand{\multirowsetup}{\centering}
\begin{center}
\begin{tabular}{|c|c|c|c|c|c|c|c|c|c|c|c|c|}
\hline
$k$ & $N$   &40           &80            & 160 & 320 & 640 &1280 \\
\hline
\multirow{6}{1cm}{1}
%\hline
 & $L^1$   &1.128e-4      & 2.826e-5 & 7.150e-6 & 1.818e-6  & 4.615e-7 & 1.166e-7\\
%\hline
 &  Order      & \quad  & 1.997       & 1.983      &  1.975       & 1.978 &1.985 \\
 %\hline
 &$L^2$   & 1.443e-4      & 3.636e-5 & 9.296e-6 &  2.388e-6  & 6.111e-7 & 1.552e-7\\
%\hline
 &Order   & \quad   & 1.989     & 1.968      &   1.961      & 1.966 & 1.977 \\
%\hline
 &$L_{\infty}$ & 4.521e-4 & 1.125e-4 & 2.925e-5  &  7.592e-6 & 1.938e-6 & 4.889e-7\\
%\hline
 &Order    & \quad  & 2.006      & 1.944     &   1.946      & 1.970 & 1.987   \\
 \hline
 \multirow{6}{1cm}{2}

 & $L^1$  & 9.595e-6      &1.225e-6 & 1.498e-7 &  1.801e-8 & 2.168e-9 & 2.640e-10 \\
%\hline
 &  Order      & \quad & 2.969    & 3.032     & 3.056    & 3.054 & 3.038  \\
% \hline
 &$L^2$  & 1.380e-5       & 1.784e-6 & 2.175e-7 &  2.579e-8 & 3.049e-9 & 3.656e-10\\
%\hline
 &Order       & \quad &2.951     & 3.036       & 3.076    & 3.080 & 3.060 \\
%\hline
 &$L_{\infty}$ & 5.748e-5 & 7.204e-6 & 9.491e-7  & 1.167e-7 & 1.347e-8 & 1.524e-9 \\
%\hline
 &Order    & \quad       &2.996    & 2.924     & 3.023   & 3.115 &3.144   \\
 \hline 
\end{tabular}%\\
\end{center}
\label{exorder}
\end{table}

}\end{exam}

% interface only problem
\begin{exam}{\em
\label{examinterfaceonly}
To demonstrate the non-oscillatory property for the pressure and velocity fields, in this example
we consider the interface only problem with the initial condition
\begin{equation*}
(\rho,v,P,\gamma,B)=
\begin{cases}
(1,1,1,1.4,1),  \;&x\leqslant 0\\
(0.125,1,1,1.9,0), \;&x>0.
\end{cases}
\end{equation*}
The computational domain is taken as (-5,5) and the inflow/outflow boundary conditions are employed at the ends.

The numerical results obtained with $N = 100$ at $T=2$
are plotted in Fig \ref{figinterface}. The figure shows that the DG-ALE method preserves the constant velocity and pressure
and appears to be free of oscillations near the material interface.
Moreover, the results show that the method is able to concentrate mesh points in the regions of the material interface, which
is a strong indication that the incorporation of the Lagrangian meshing with the MMPDE moving mesh method works well.
Furthermore, the close-up of the density at the interface demonstrates that
$P^2$-DG gives more accurate results than $P^1$-DG.
   
\begin{figure}[hbtp]
\centering
\subfigure[Density]{
\includegraphics[width=0.44\textwidth]{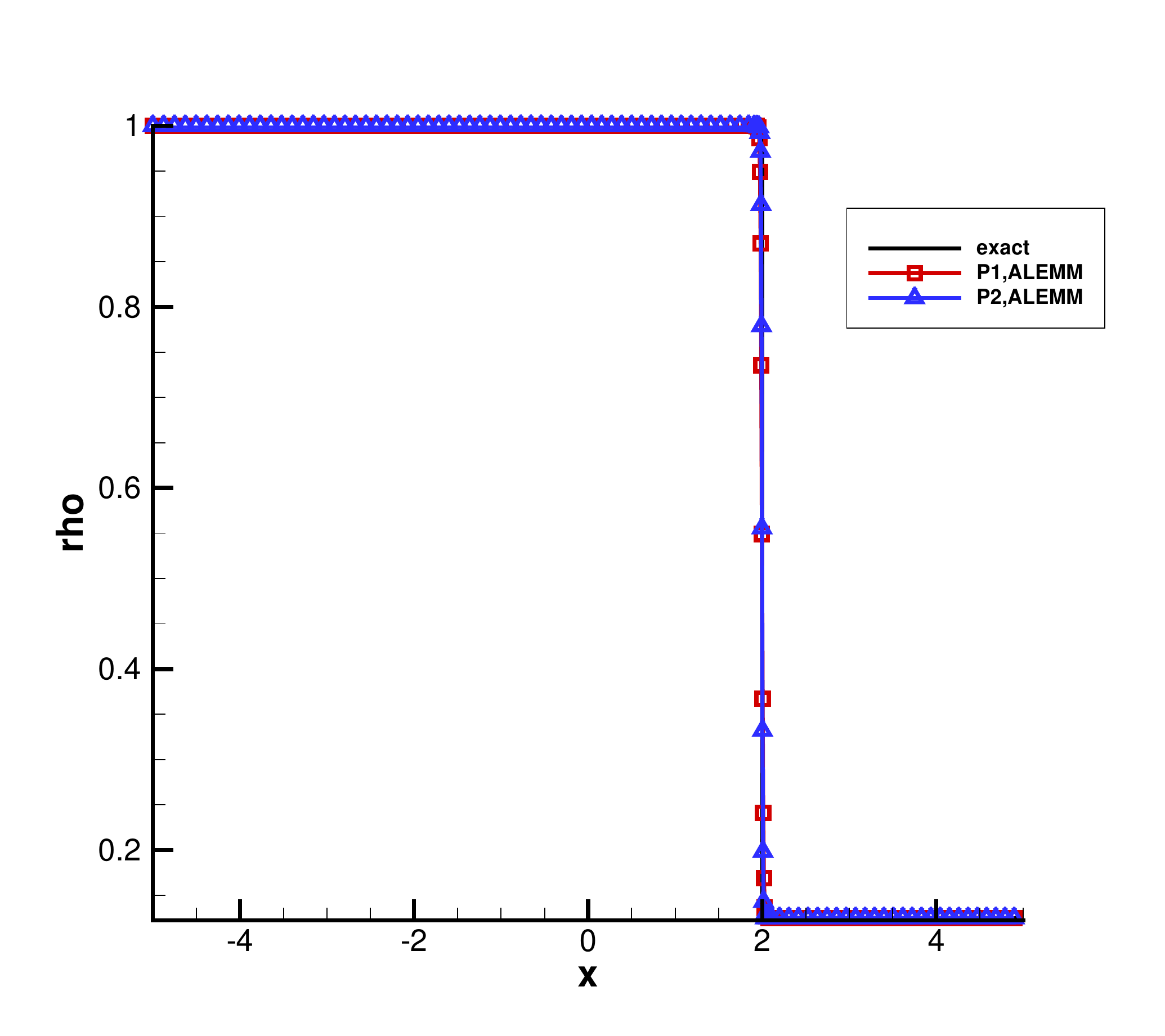}}
\subfigure[Close view of (a)]{
\includegraphics[width=0.44\textwidth]{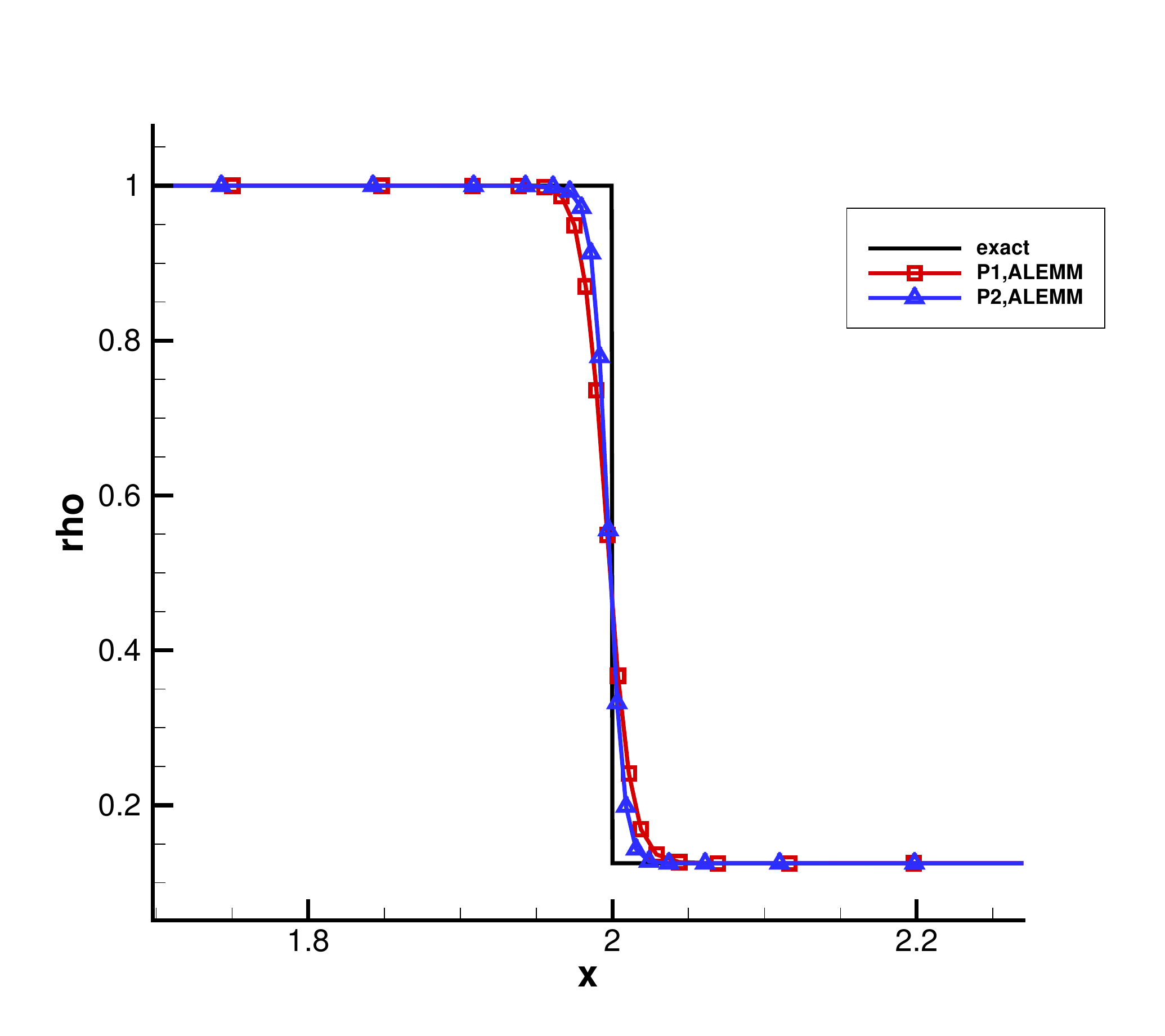}}
\subfigure[Velocity]{
\includegraphics[width=0.44\textwidth]{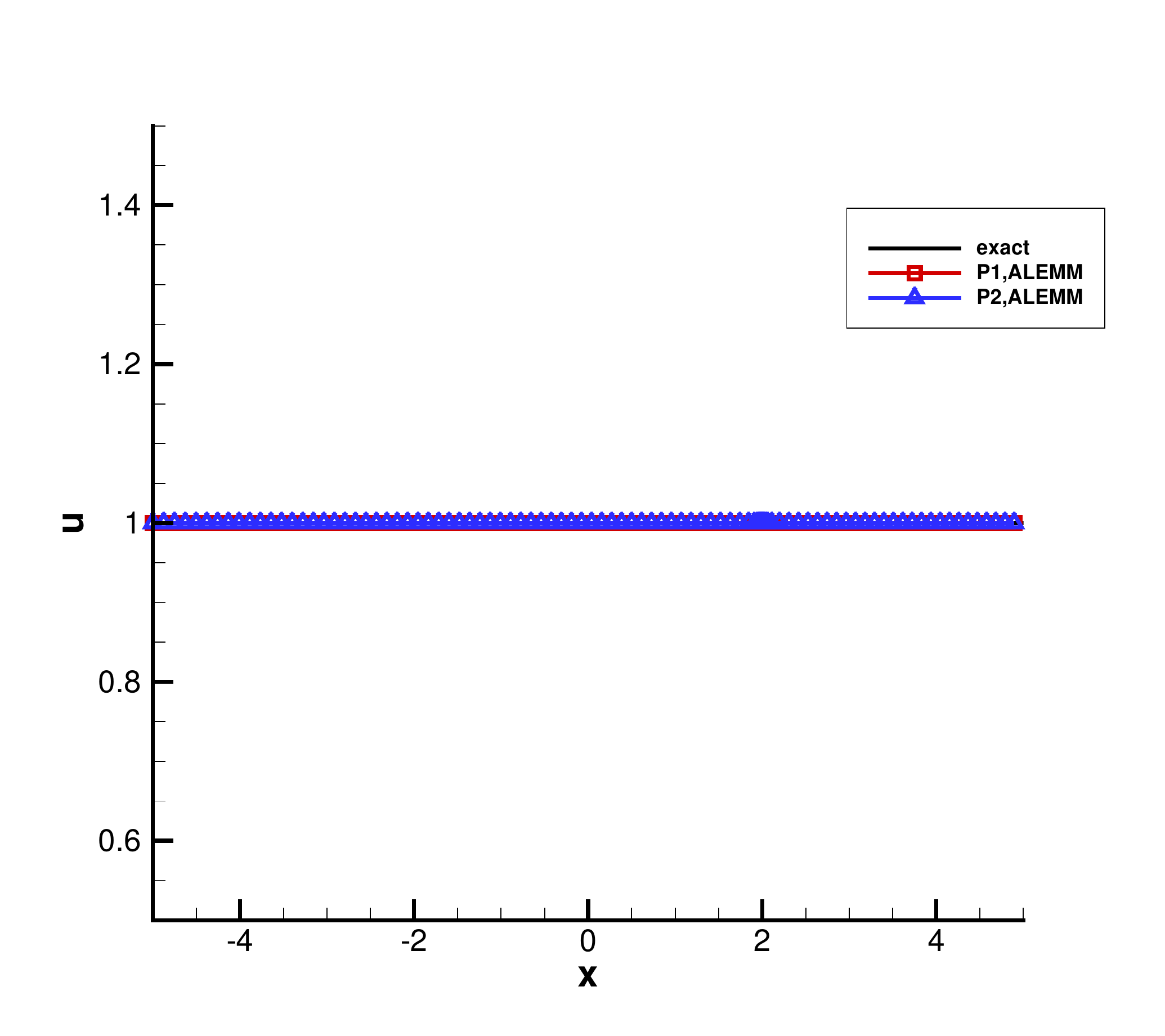}}
\subfigure[Pressure]{
\includegraphics[width=0.44\textwidth]{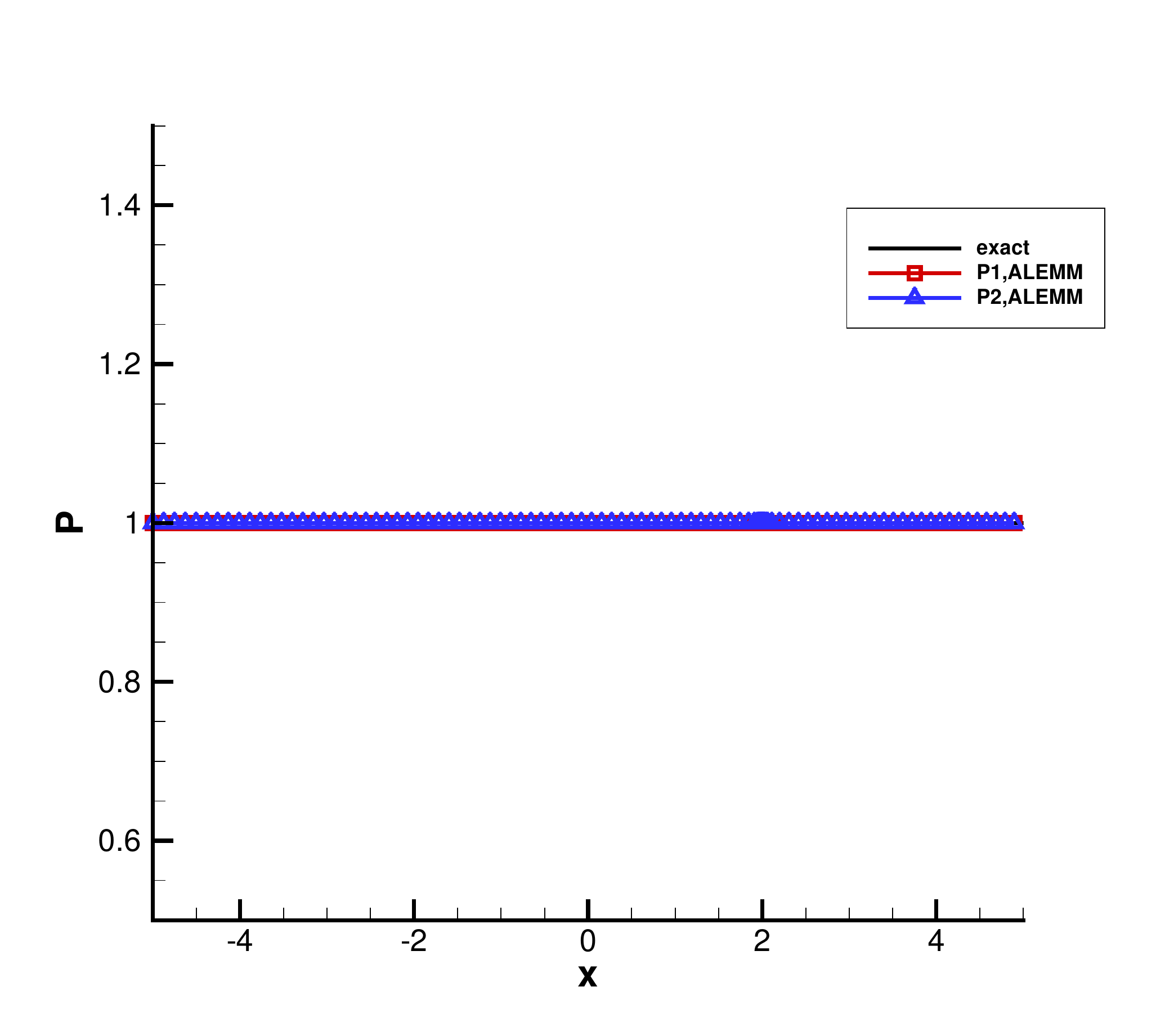}}
\subfigure[Mesh trajactories with $P^1$ elements]{
\includegraphics[width=0.44\textwidth]{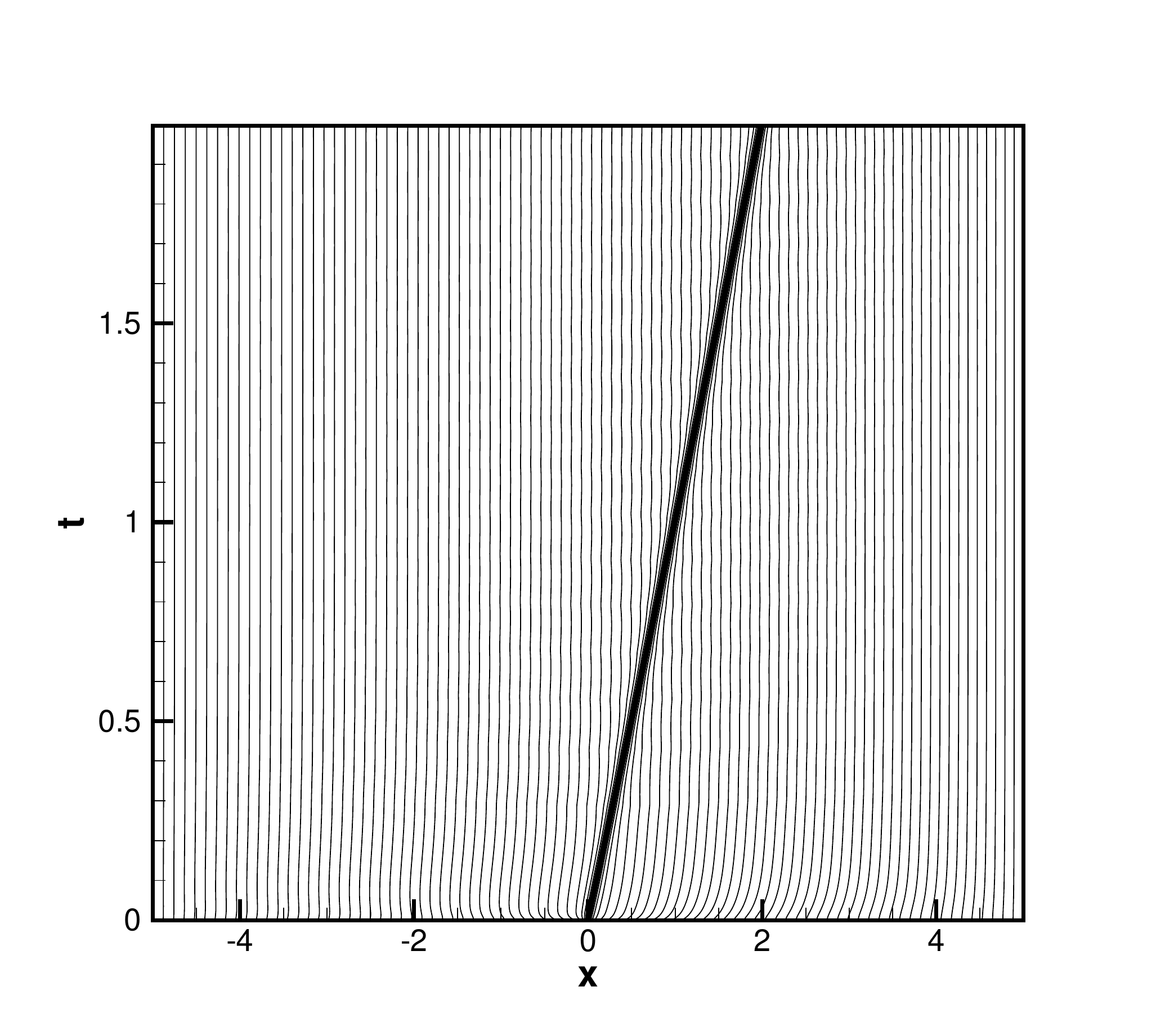}}
\subfigure[Mesh trajactories with $P^2$ elements]{
\includegraphics[width=0.44\textwidth]{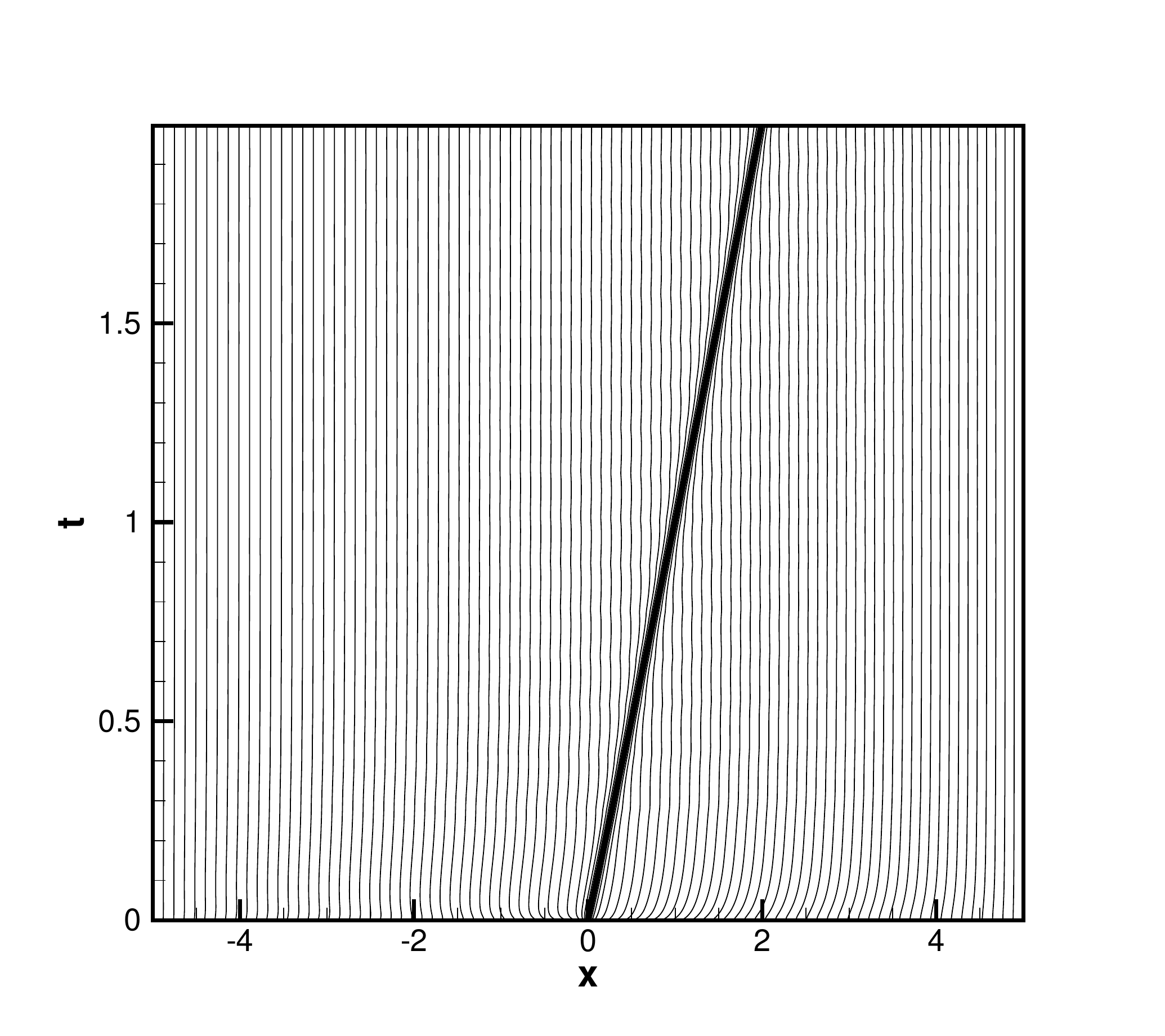}}
\caption{Example~\ref{examinterfaceonly} The DG-ALE method with $N=100$.}
\label{figinterface}
\end{figure}

}\end{exam}

  % Shu-Osher problem
  \begin{exam}{\em
\label{examshuosher}
To test the numerical viscosity of the method, the Shu-Osher problem, which was first introduced for a single material in \cite{shu1989},
is considered for two components  \cite{liu2017} in this example (where each of $\gamma$ and $B$ takes two different values 
on the domain).
This problem contains both shocks and complex smooth region structures. We solve the Euler equations with a moving shock interacting
with a sine wave in density. The initial data is
\begin{equation*}
(\rho,v,P,\gamma,B)=
\begin{cases}
(3.857143,2.629369,\frac{31}{3},1.4,1),  \;&x\leqslant -4\\
(1+0.2\sin(5x),0,1,1.9,0), \;&x>-4.
\end{cases}
\end{equation*}
The inflow and outflow boundary conditions are employed at the ends. The physical domain is taken as (-5, 5).

The computed density obtained with $N = 150$ at $t=1.8$ is plotted in Fig. \ref{figshuoshersmo40} where the reference solution is
obtained by the second-order DG method with 2000 uniform points. The figure shows that both the shocks and the complex structures
in the smooth regions are resolved well by the DG-ALE scheme while $P^2$-DG is more accurate than $P^1$-DG.
The mesh trajectories Fig. \ref{figshuoshersmo40} also show that the mesh points are concentrated near the moving shocks and
smooth region structures for all time in the time range of the computation.
   
\begin{figure}[hbtp]
\centering
\subfigure[Density]{
\includegraphics[width=0.44\textwidth]{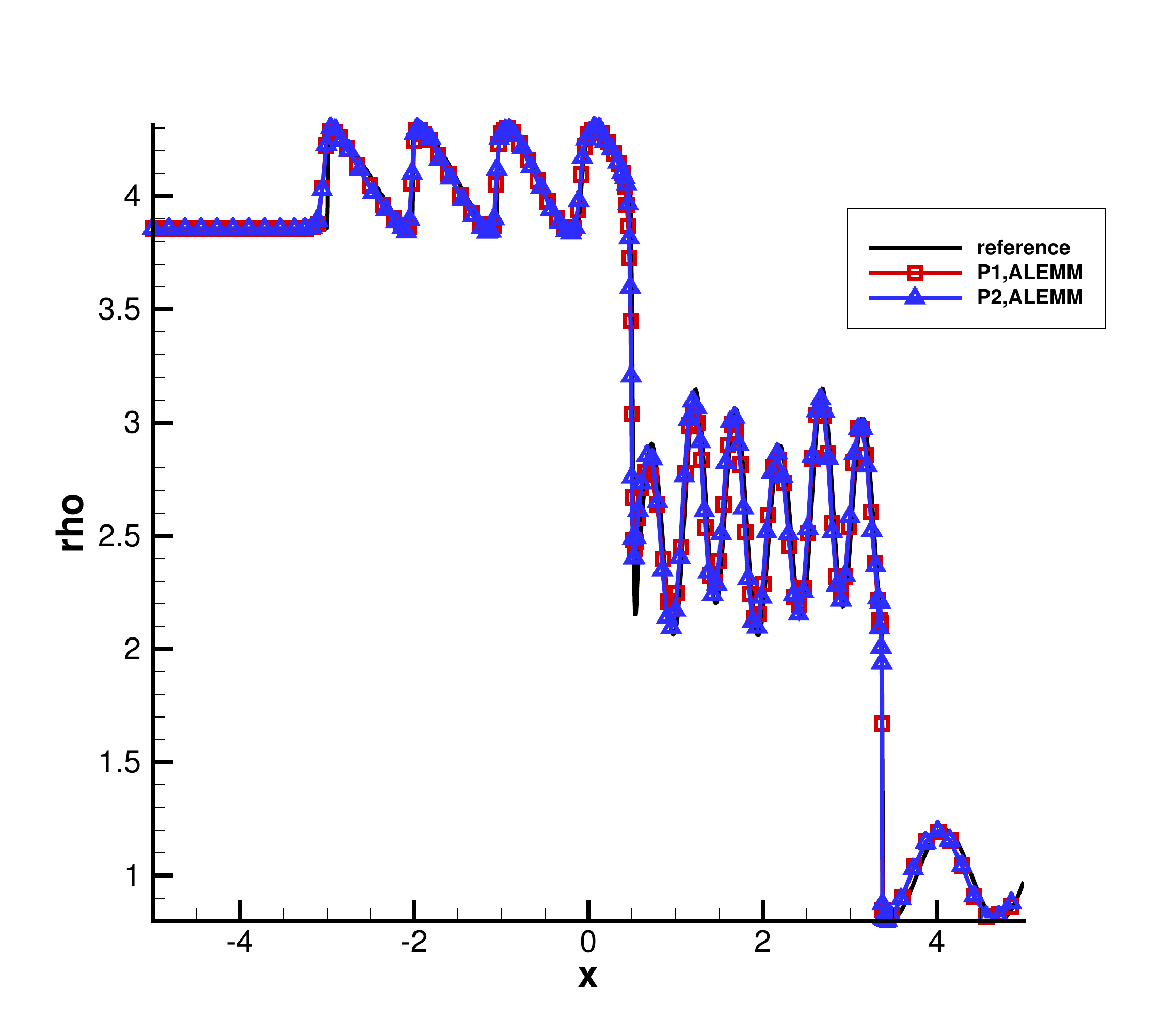}}
\subfigure[Close view of (a)]{
\includegraphics[width=0.44\textwidth]{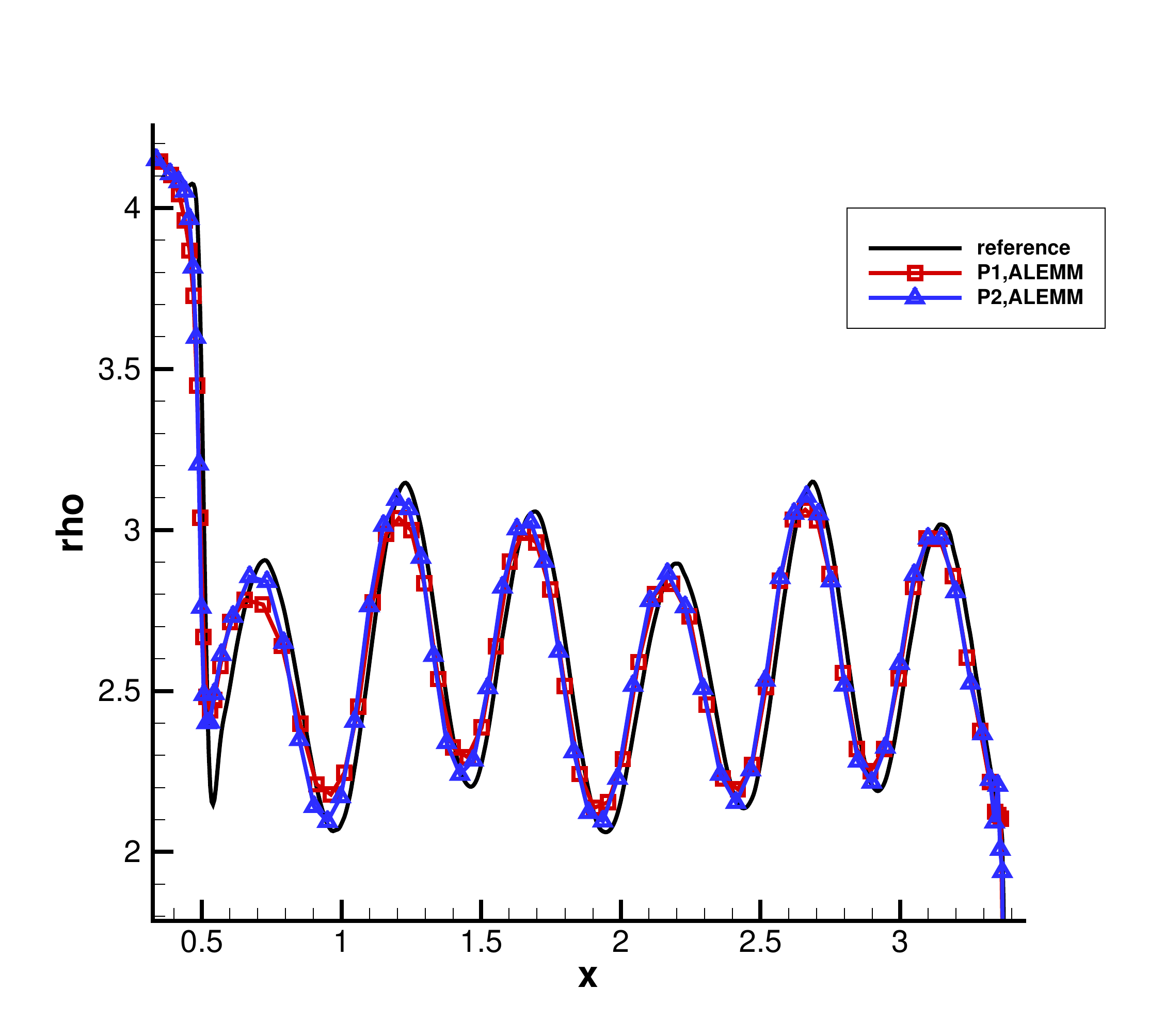}}
\subfigure[Mesh trajactories with $P^1$ elements]{
\includegraphics[width=0.44\textwidth]{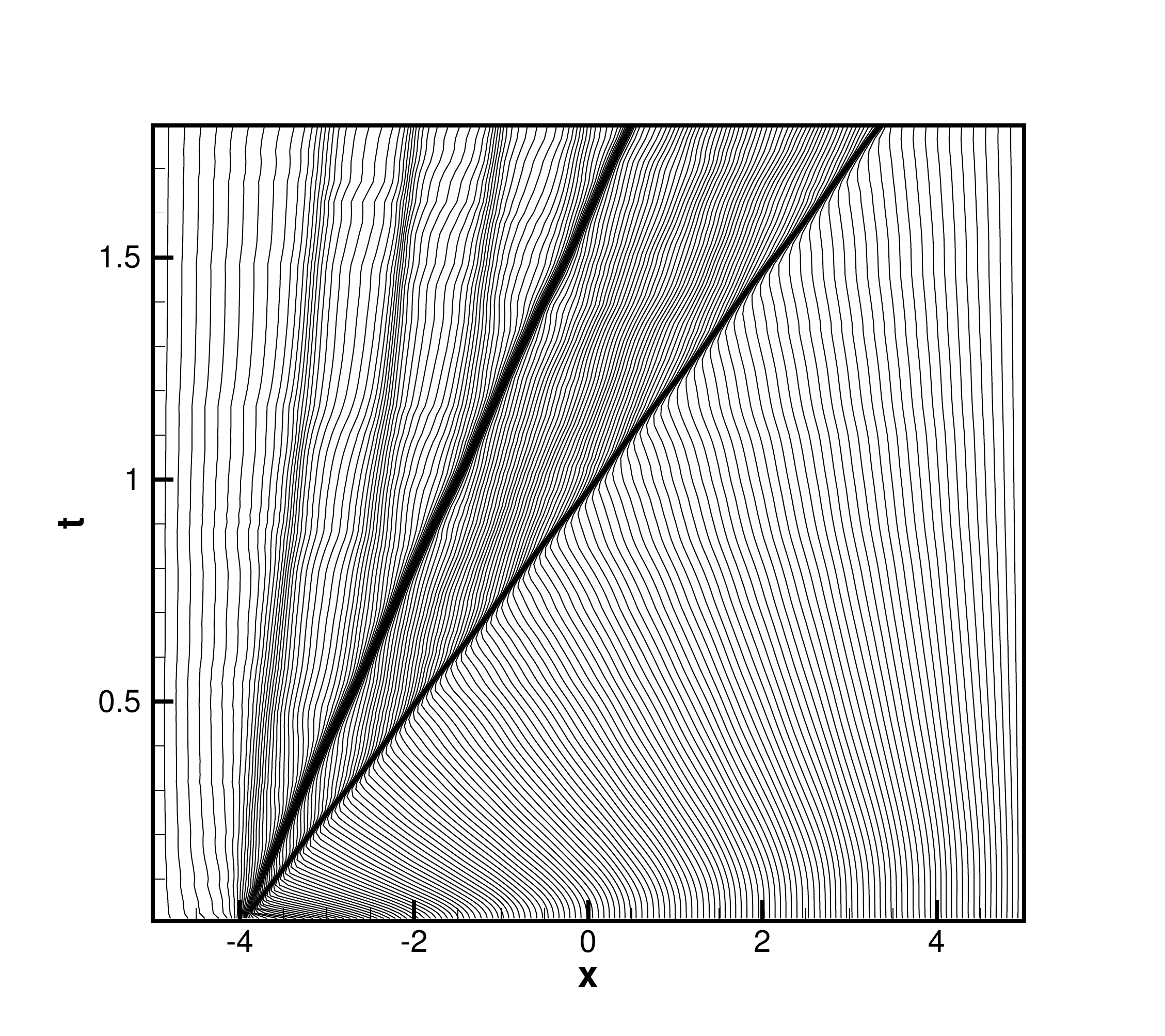}}
\subfigure[Mesh trajactories with $P^2$ elements]{
\includegraphics[width=0.44\textwidth]{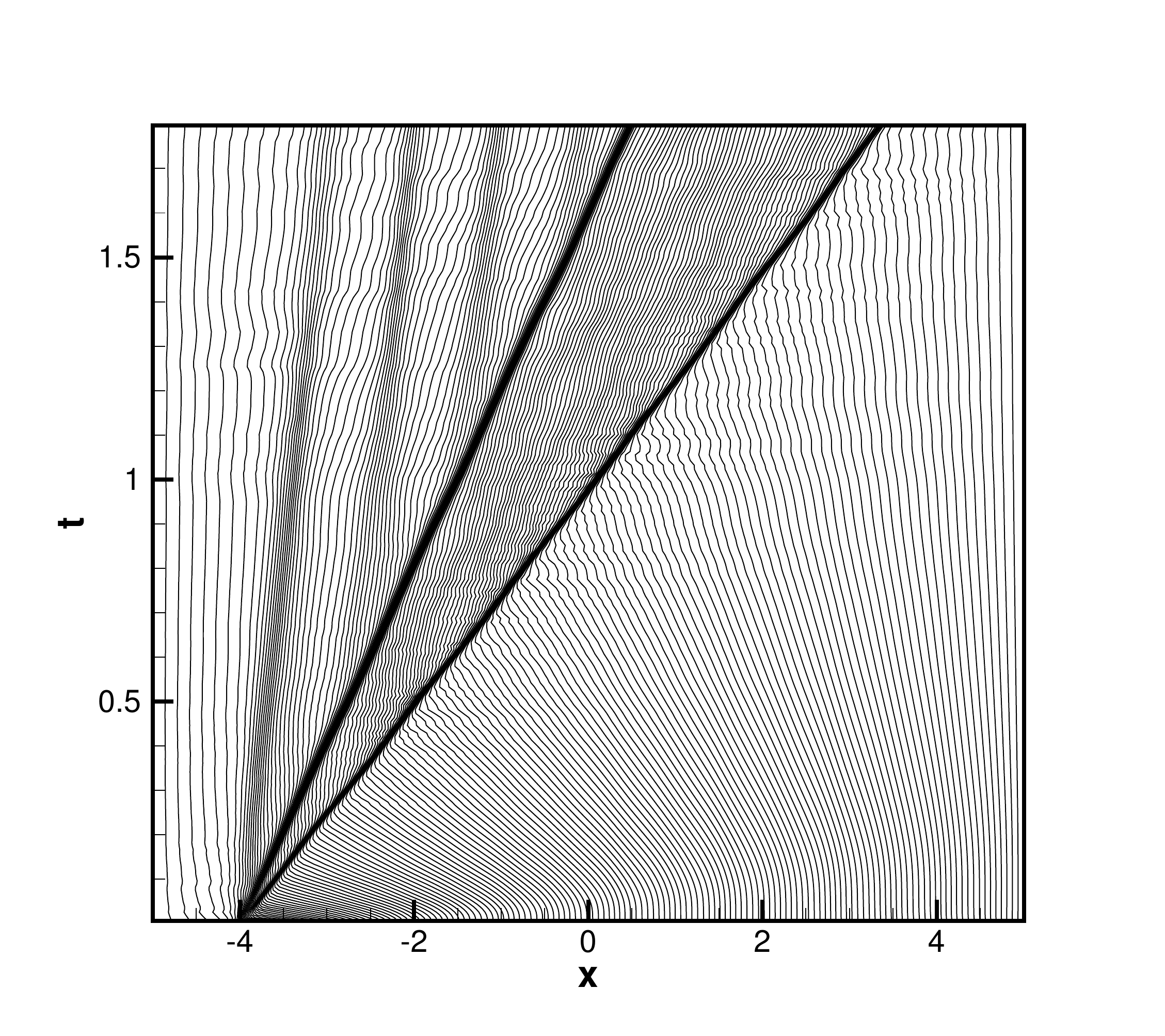}}
\caption{Example~\ref{examshuosher} The DG-ALE method with $N=150$.}
\label{figshuoshersmo40}
\end{figure}

}
\end{exam}

%%%%%%%%%gas-liquid
\begin{exam}{\em
\label{examgasliquid}
In this example the gas-liquid shock tube test with a strong shock wave is considered. This test is challenging since the shock
and the material interface are close to each other and the pressure ratio is excessively high. The initial condition is
\begin{equation*}
(\rho,v,P,\gamma,B)=
\begin{cases}
(10^3,0,10^9,4.4,6\times 10^8),  \; & x\leqslant 0.5\\
(50,0,10^5,1.4,0), \;  &x>0.5
\end{cases}
\end{equation*}
and the computational domain is (-0.2, 1). The inflow and outflow boundary conditions are employed at the ends. The final time is $t=0.0002$.
For this problem, an analytical solution is available {\cite{razo2016}}.

The numerical results with $N = 2000$ are plotted in Fig. \ref{figgasliquid}. One can see that both the shock and the material interface are resolved
by the scheme. Moreover,  $P^2$-DG produces more accurate solutions than $P^1$-DG.
%The trajectories of the mesh with 1000 points are plotted in Fig. \ref{figgasliquid}.

To show the advantage of the current scheme, it is compared with the pure Lagrangian DG method (marked by ``LAG"),
the ALE using the Winslow method (marked using ``ALEWIN") to smooth and improve the Lagrangian mesh,
and the moving mesh DG method (marked by ``MM") that is the DG-ALE method without using the Lagrangian mesh
in Figs. \ref{figgasliquidwinmm}, \ref{figgasliquidmmlag} and \ref{figgasliquidmmmm} , respectively.
From these figures, we can see that the current DG-ALE method gives more accurate results than ALEWIN, MM, and LAG.
Once again, this demonstrates that the incorporation of the Lagrangian meshing with the MMPDE moving mesh method used
in the current DG-ALE method works really well to concentrate mesh points in regions of shocks and material interfaces.
Moreover, the solutions obtained with LAG are oscillatory at the discontinuities. 
 
%%%%%%%%%%%%%ALE P1P2
\begin{figure}[hbtp]
\centering
\subfigure[Density]{
\includegraphics[width=0.44\textwidth]{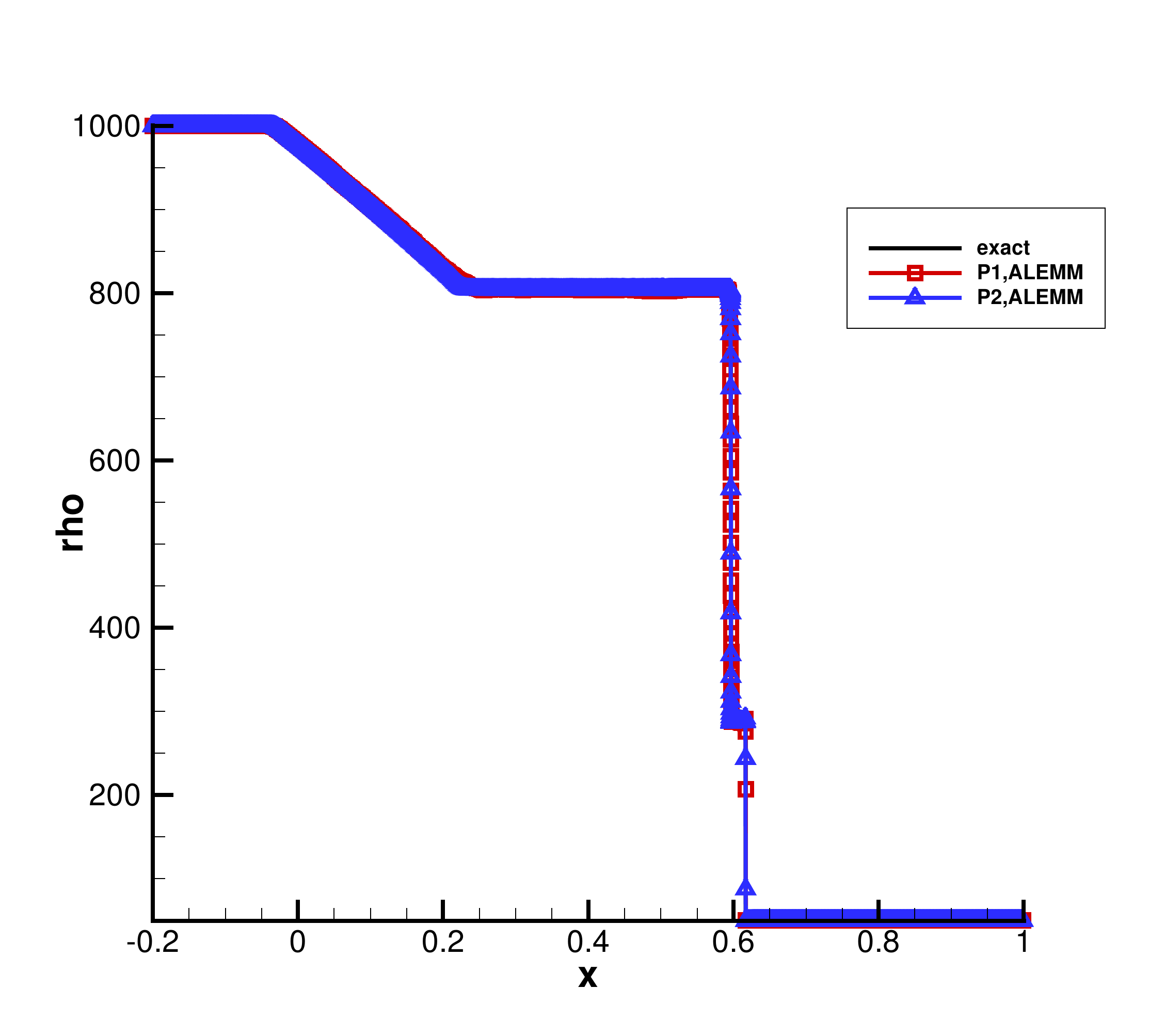}}
\subfigure[Close view of (a)]{
\includegraphics[width=0.44\textwidth]{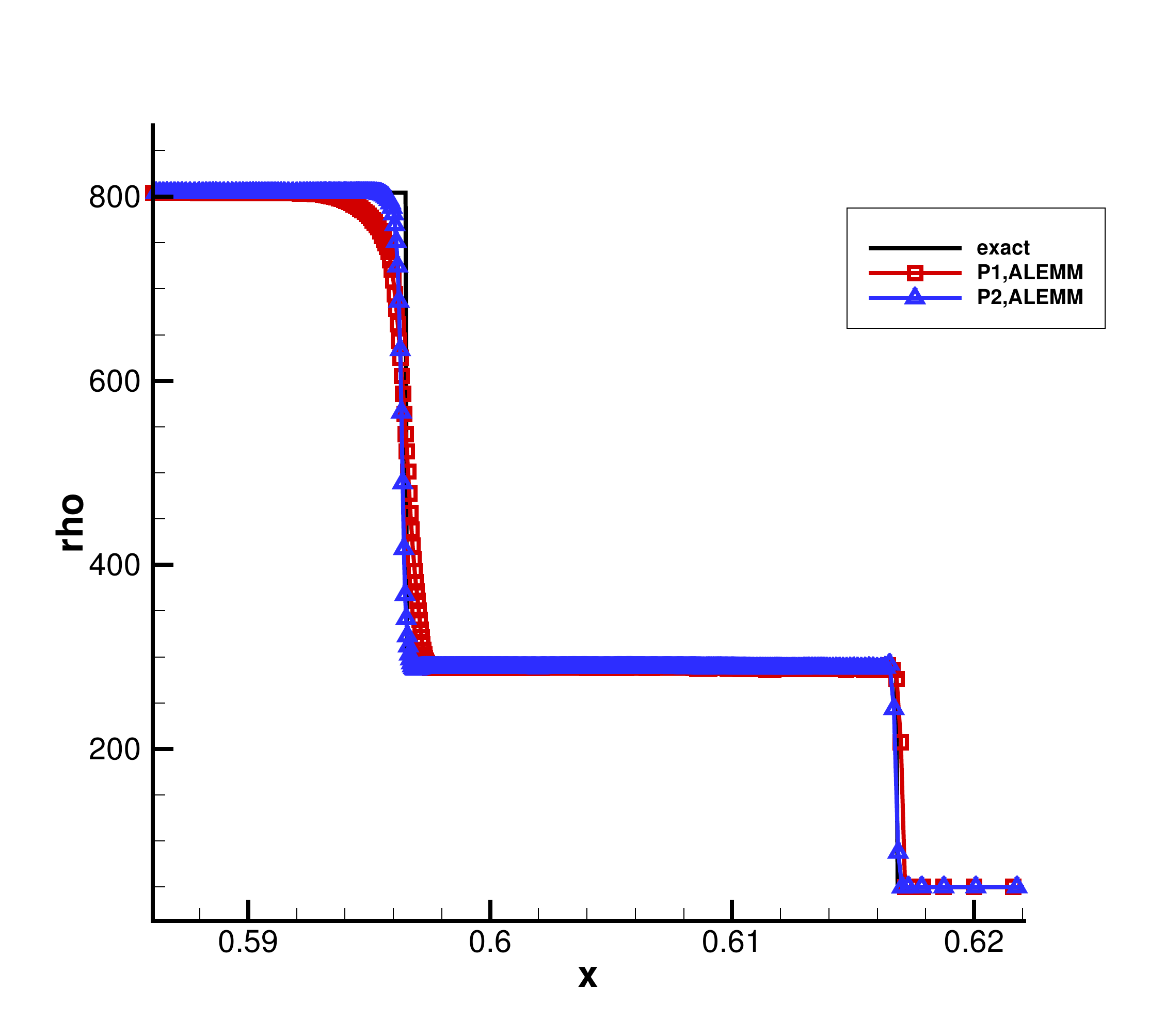}}
\subfigure[Mesh trajactories with $P^1$ elements]{
\includegraphics[width=0.44\textwidth]{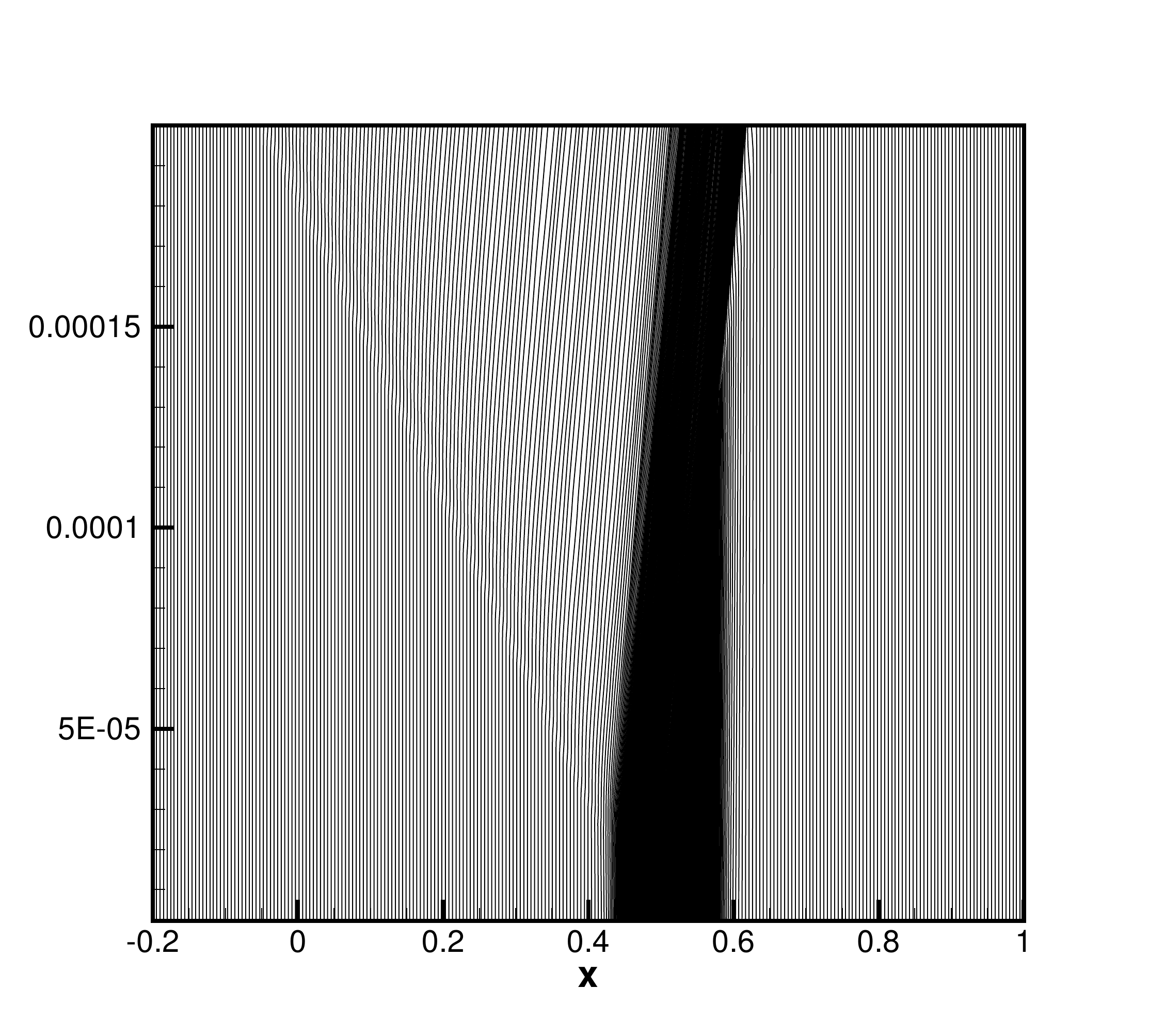}}
\subfigure[Mesh trajactories with $P^2$ elements]{
\includegraphics[width=0.44\textwidth]{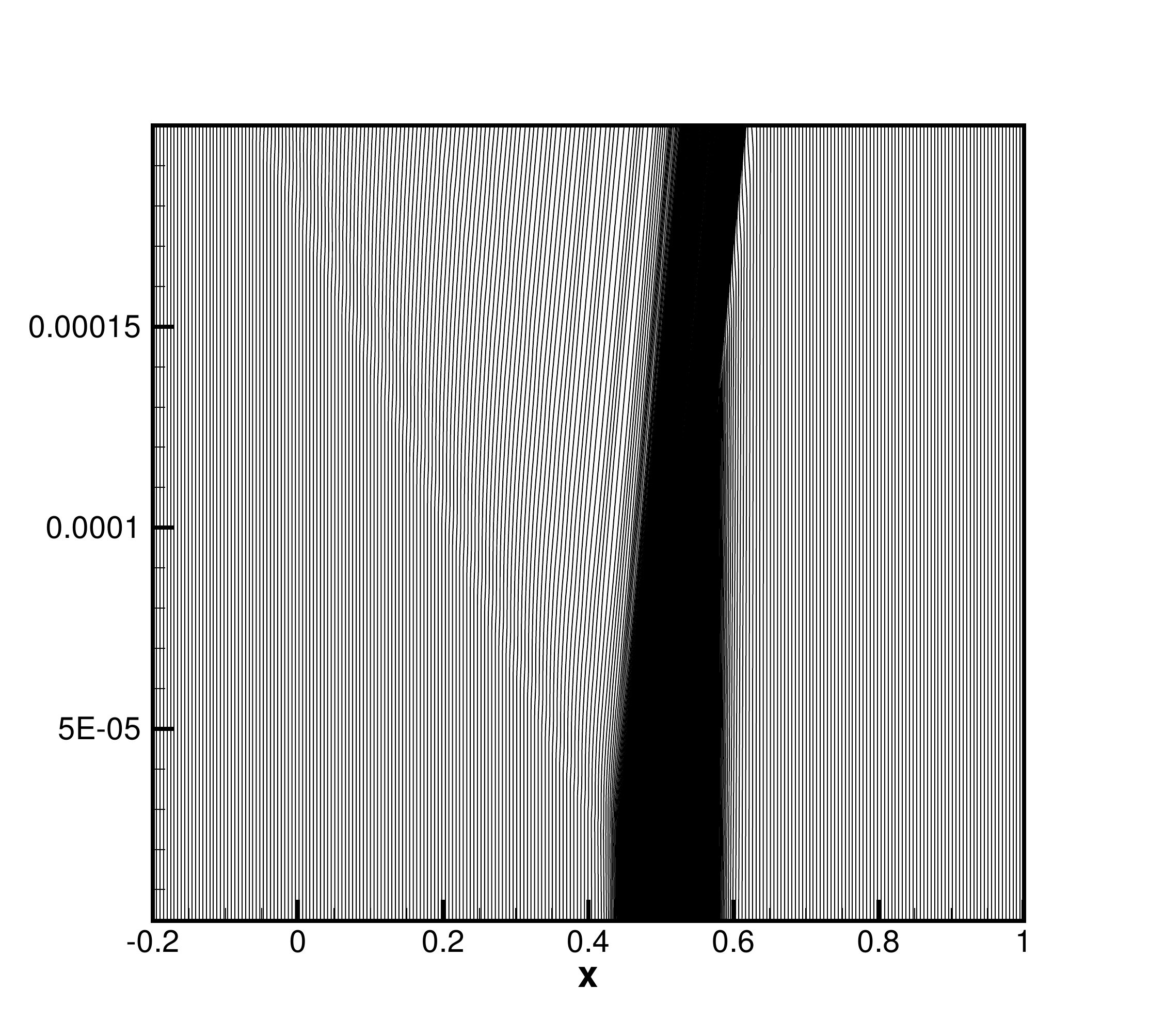}}
\caption{Example~\ref{examgasliquid} The DG-ALE method with $N=2000$.}
\label{figgasliquid}
\end{figure}
%%%%%%%%%%%%%%%

%%%%%%%%%%%%ALEMM and ALEWIN
\begin{figure}[hbtp]
\centering
\subfigure[Density, $P^1$-DG]{
\includegraphics[width=0.44\textwidth]{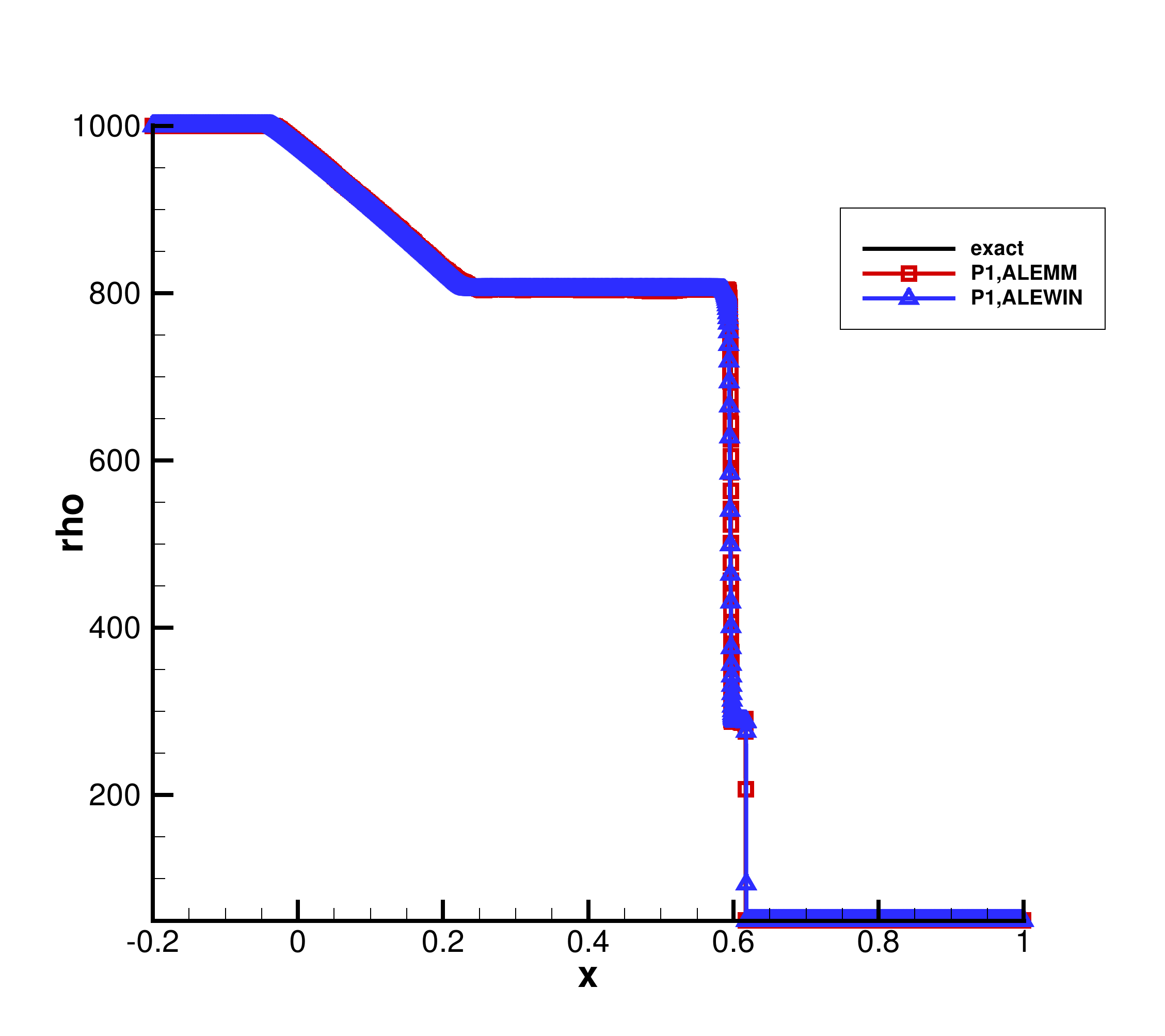}}
\subfigure[Close view of (a)]{
\includegraphics[width=0.44\textwidth]{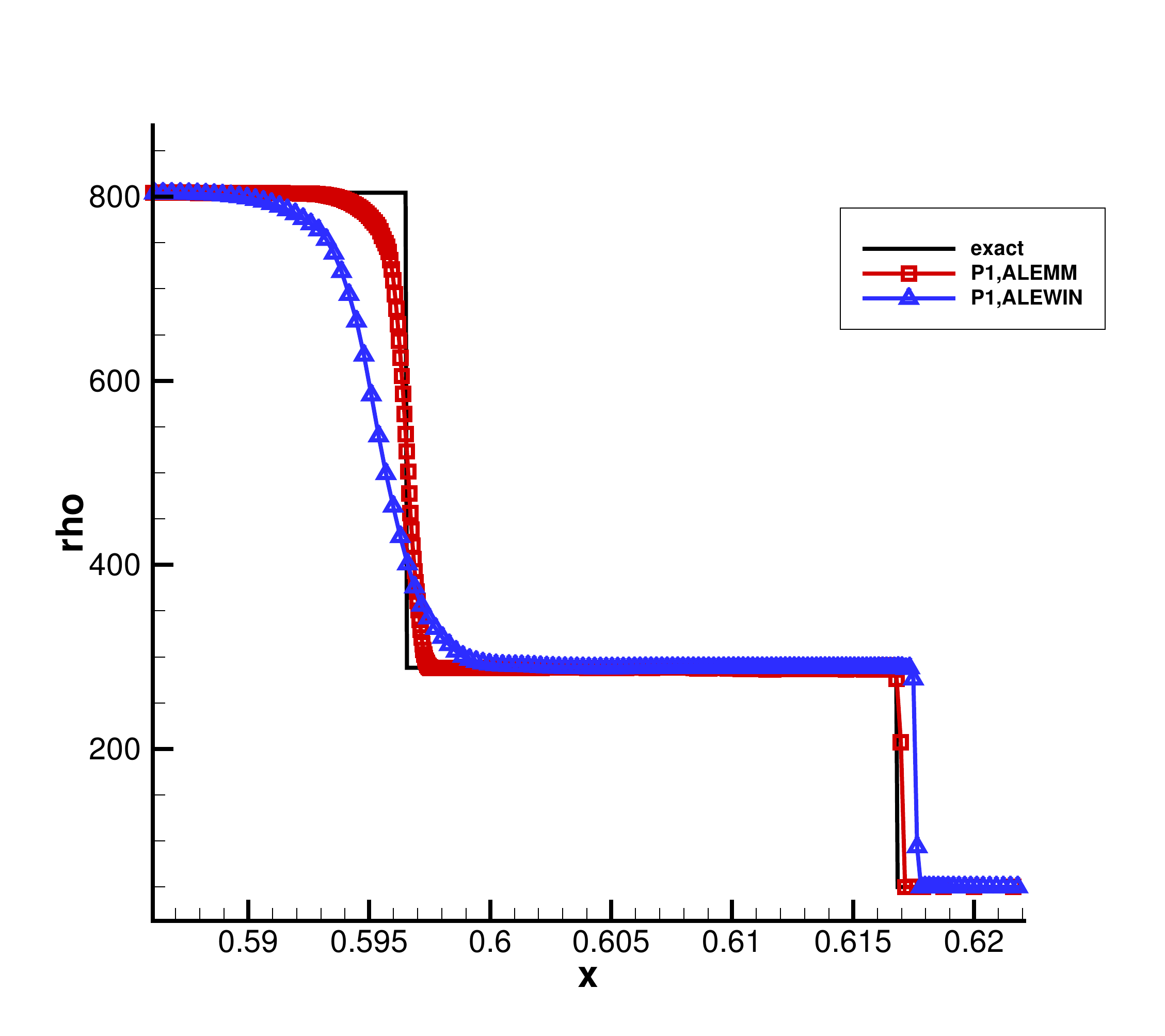}}
\subfigure[Density, $P^2$-DG]{
\includegraphics[width=0.44\textwidth]{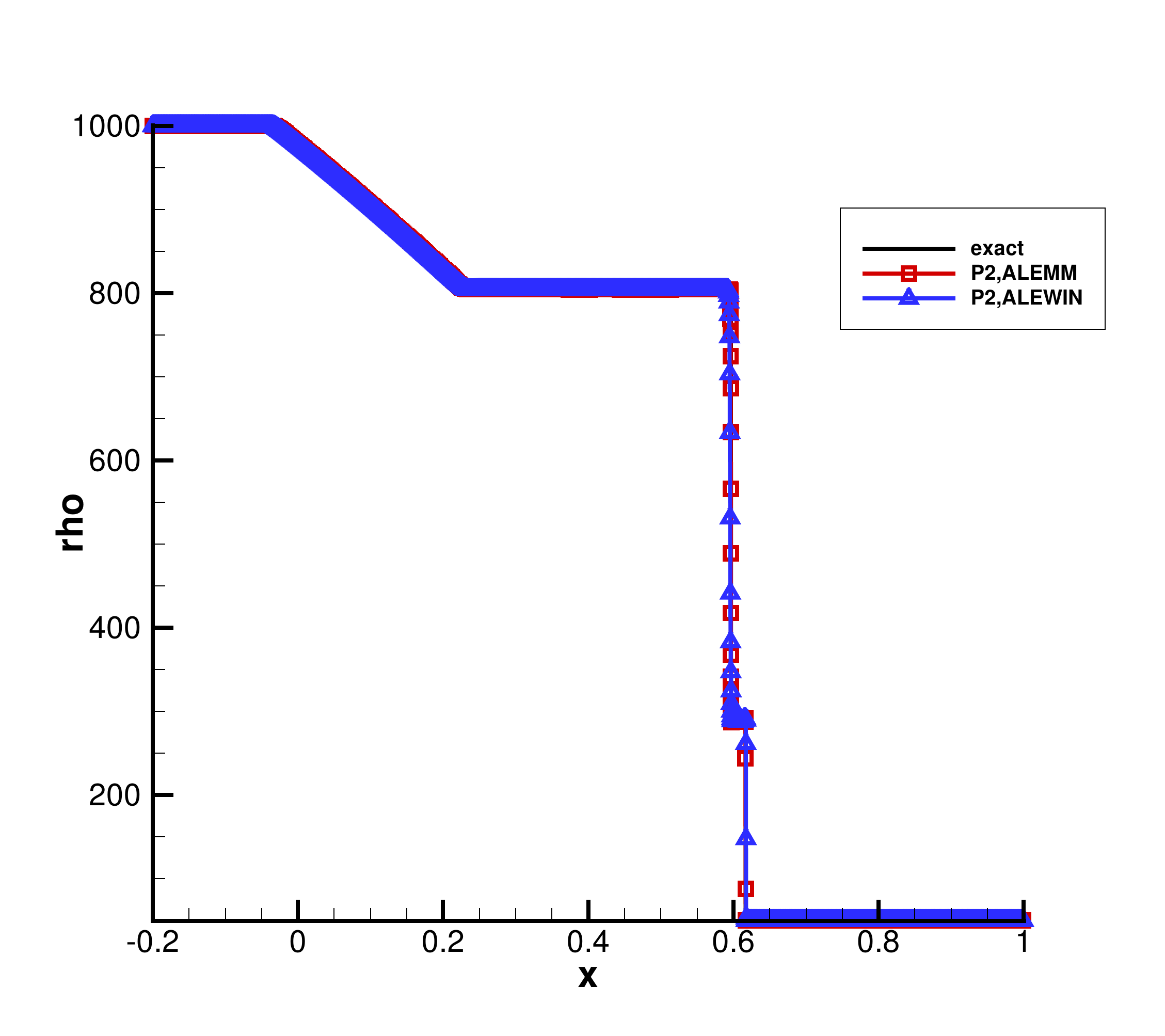}}
\subfigure[Close view of (c)]{
\includegraphics[width=0.44\textwidth]{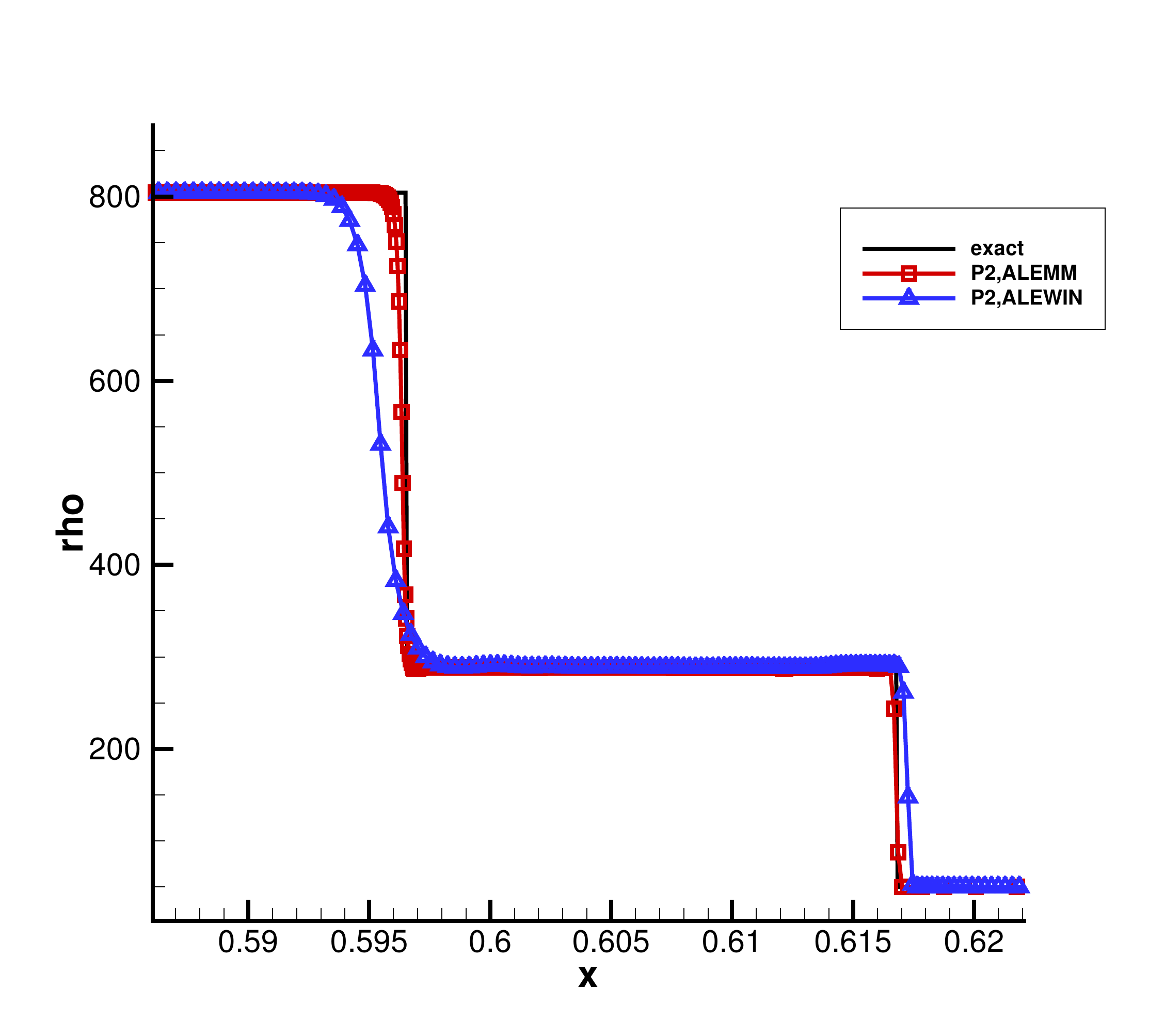}}
\caption{Example~\ref{examgasliquid} The DG-ALE method with moving mesh is compared with DG-ALE with Winslow smoothness. $N=2000$.}
\label{figgasliquidwinmm}
\end{figure}
%%%%%%%%%%end

%%%%%%%%%%ALEMM and Lag
\begin{figure}[hbtp]
\centering
\subfigure[Density, $P^1$-DG]{
\includegraphics[width=0.44\textwidth]{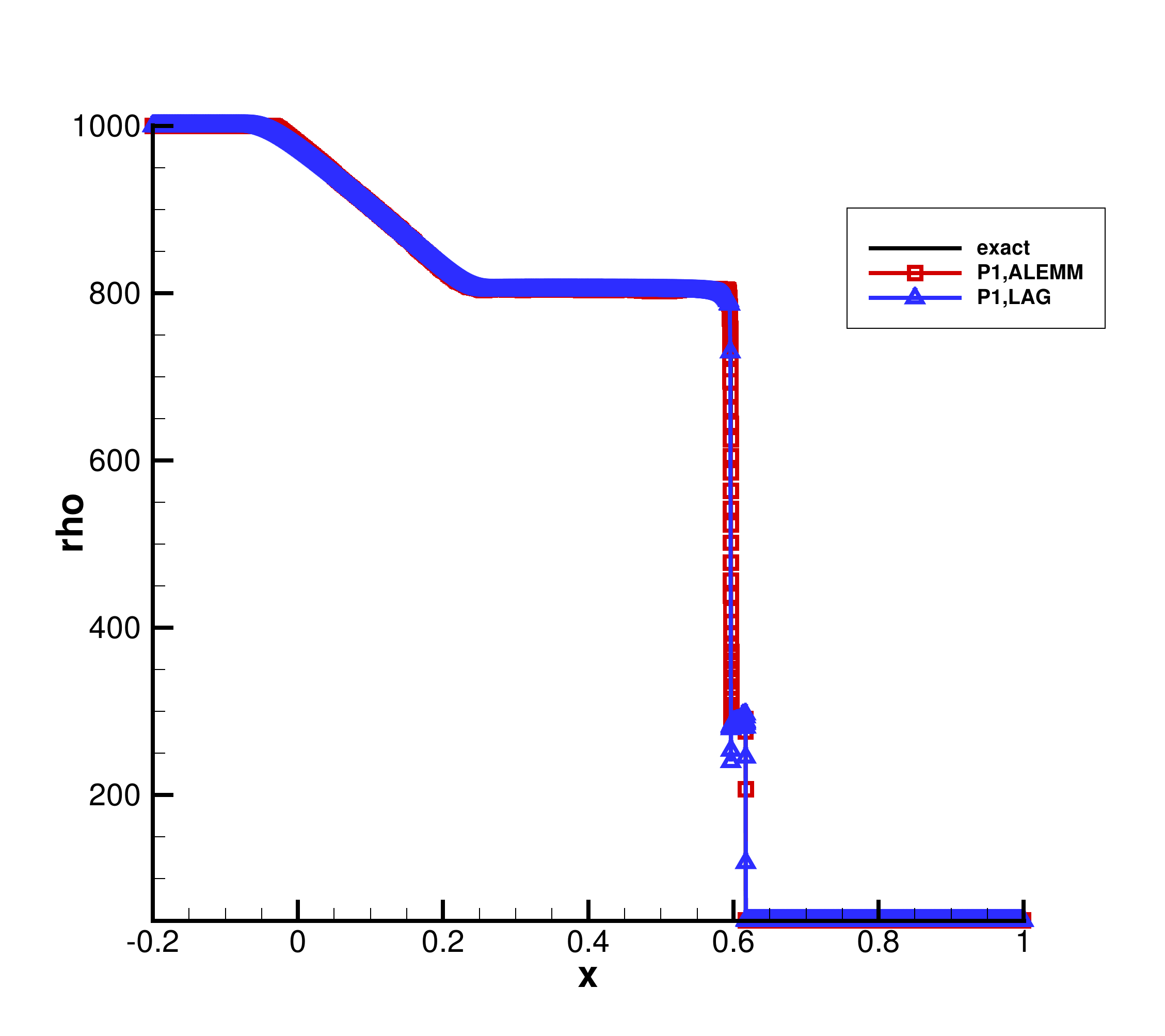}}
\subfigure[Close view of (a)]{
\includegraphics[width=0.44\textwidth]{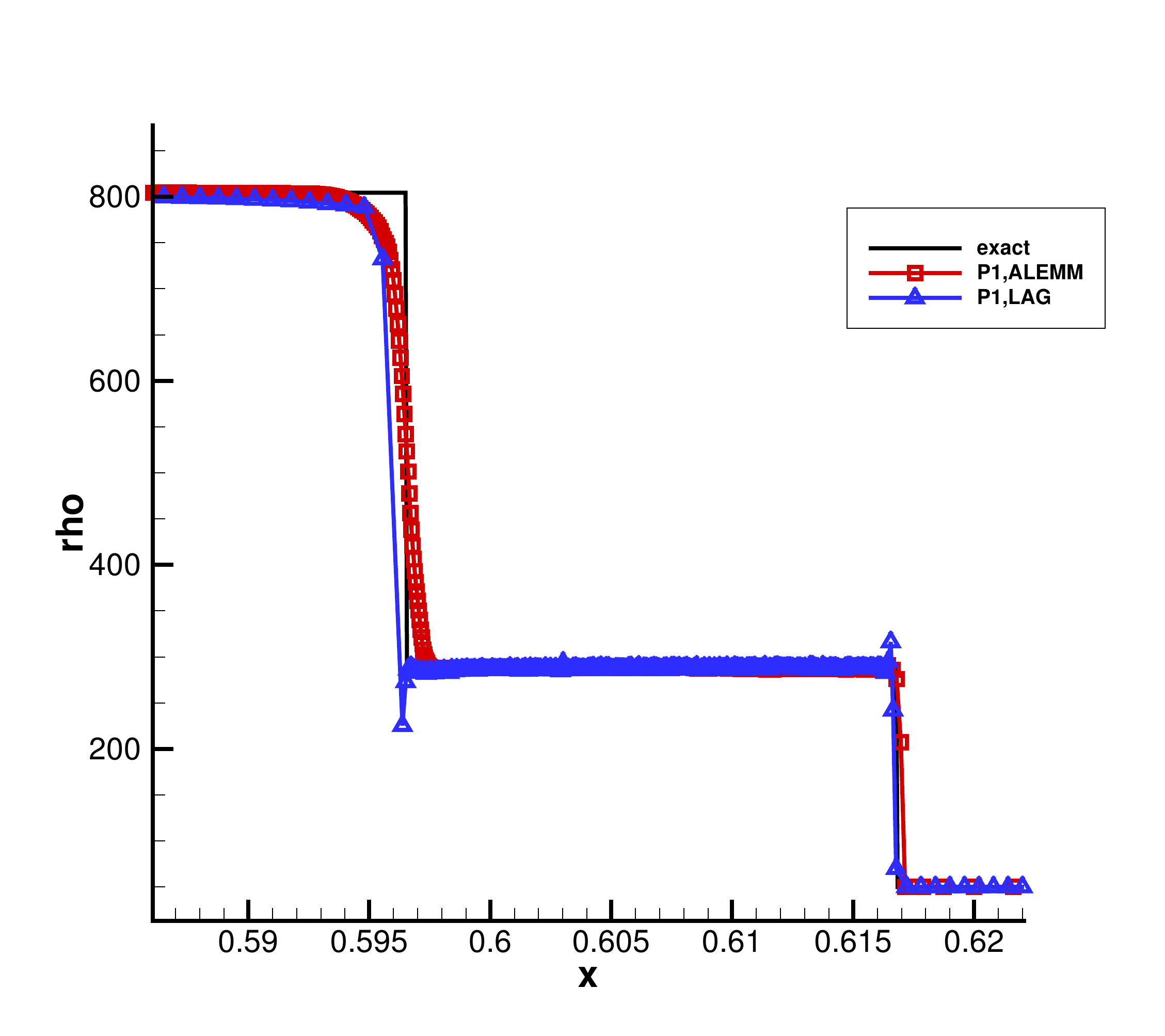}}
\subfigure[Density, $P^2$-DG]{
\includegraphics[width=0.44\textwidth]{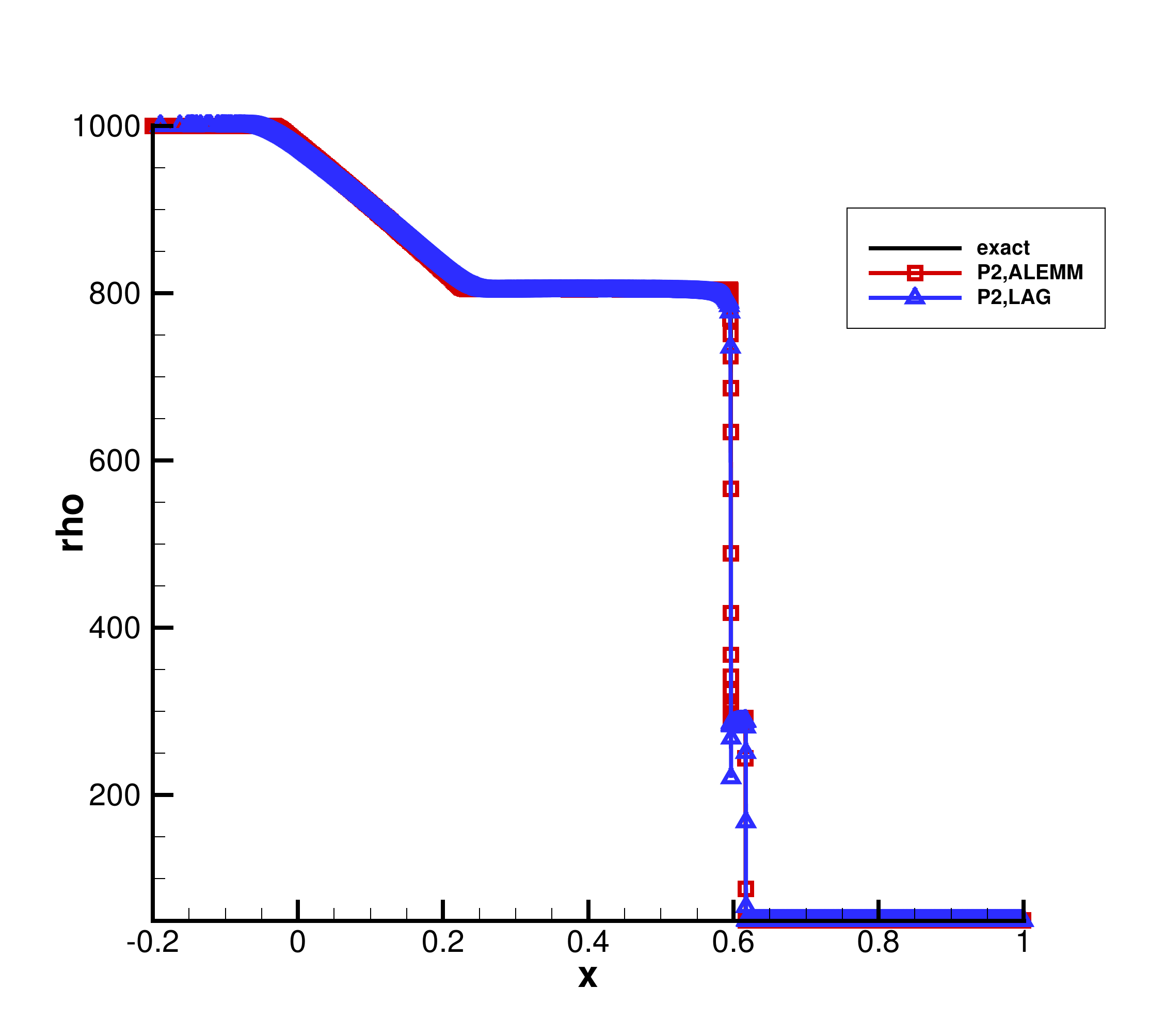}}
\subfigure[Close view of (c)]{
\includegraphics[width=0.44\textwidth]{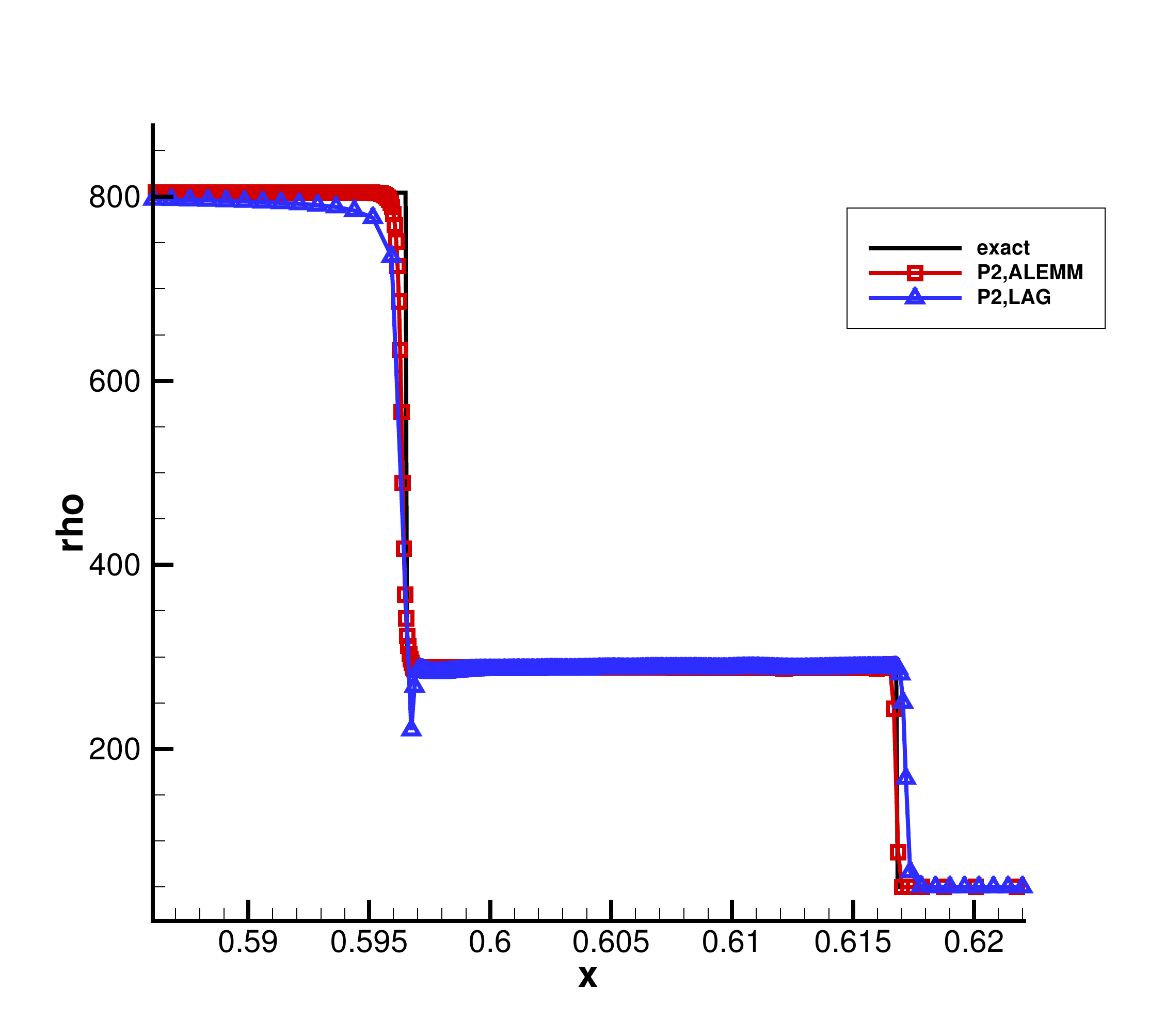}}
\caption{Example~\ref{examgasliquid} The DG-ALE method with moving mesh is compared with the Lagrangian DG method. $N=2000$.}
\label{figgasliquidmmlag}
\end{figure}
%%%%%%%%%%%%%%%%%end

%%%%%%%%%%ALEMM and MM
\begin{figure}[hbtp]
\centering
\subfigure[Density, $P^1$-DG]{
\includegraphics[width=0.44\textwidth]{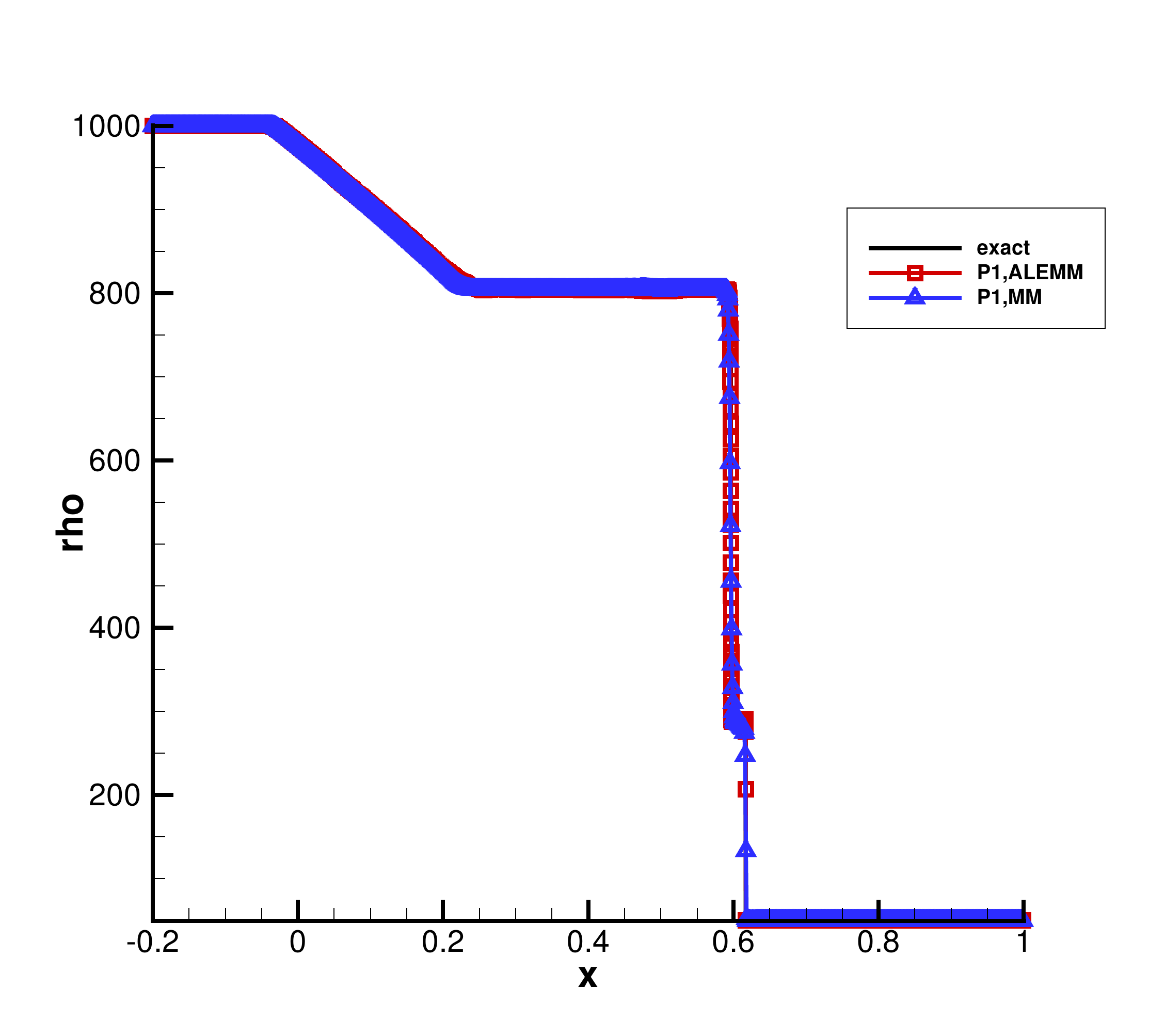}}
\subfigure[Close view of (a)]{
\includegraphics[width=0.44\textwidth]{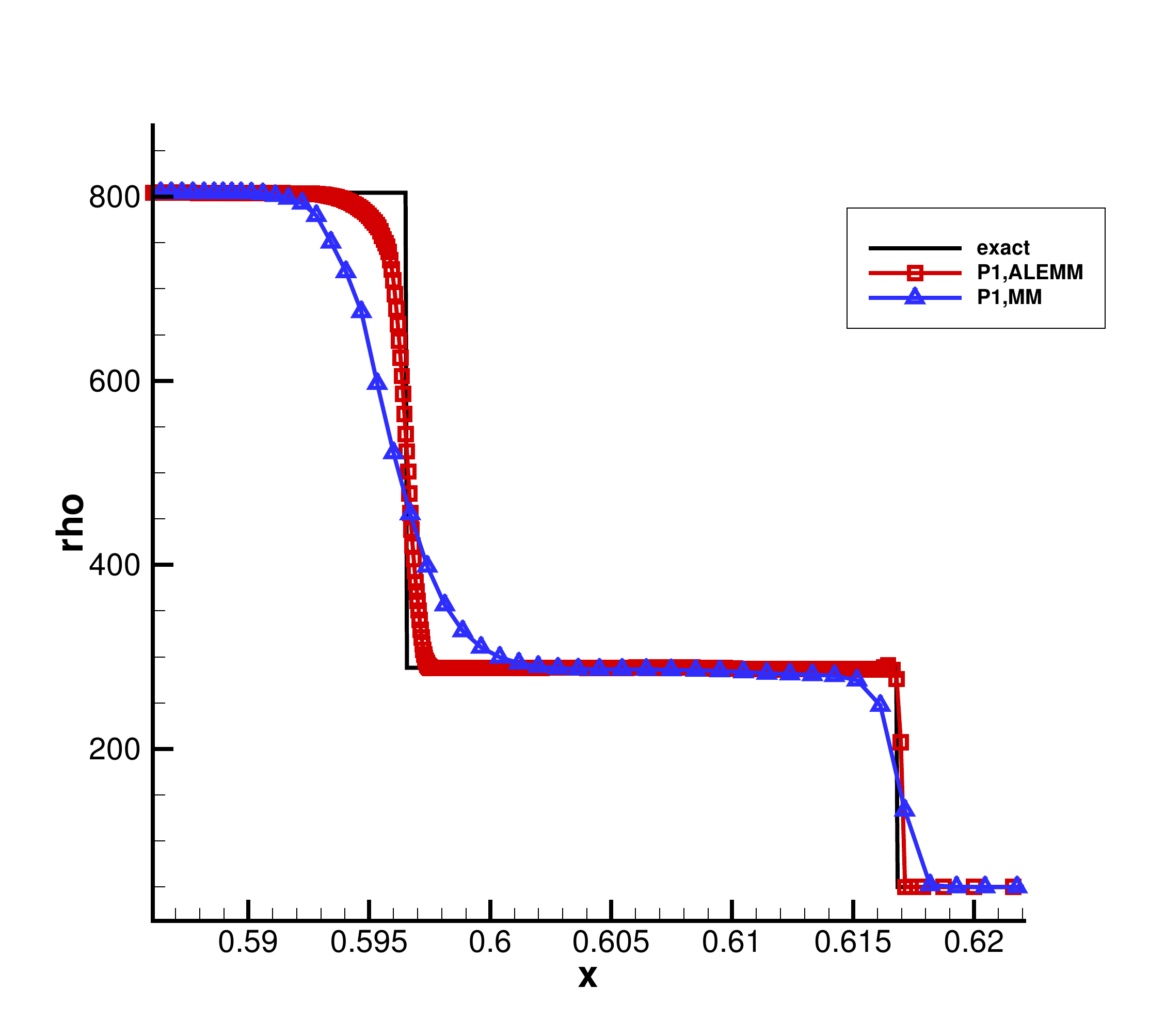}}
\subfigure[Density, $P^2$-DG]{
\includegraphics[width=0.44\textwidth]{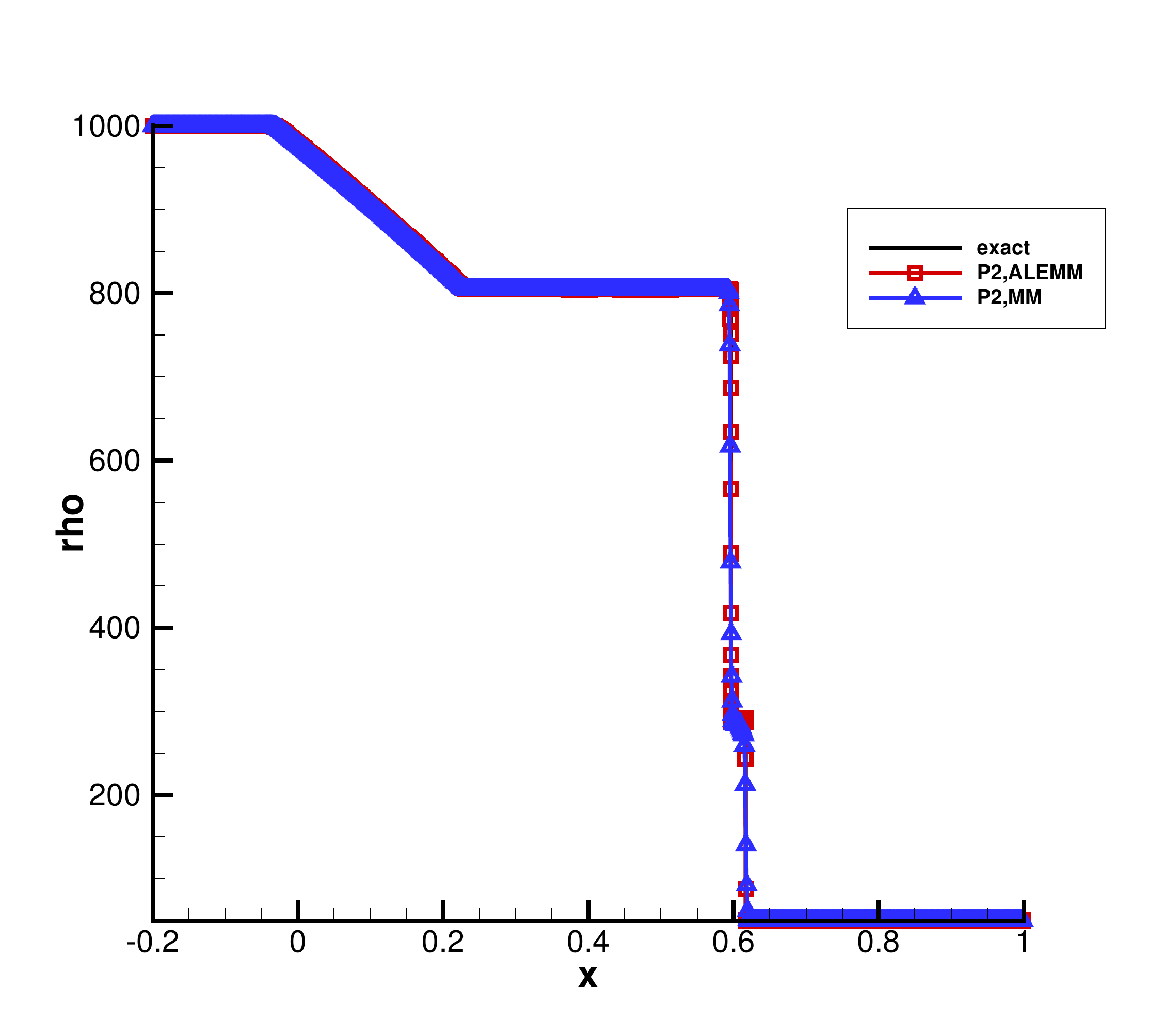}}
\subfigure[Close view of (c)]{
\includegraphics[width=0.44\textwidth]{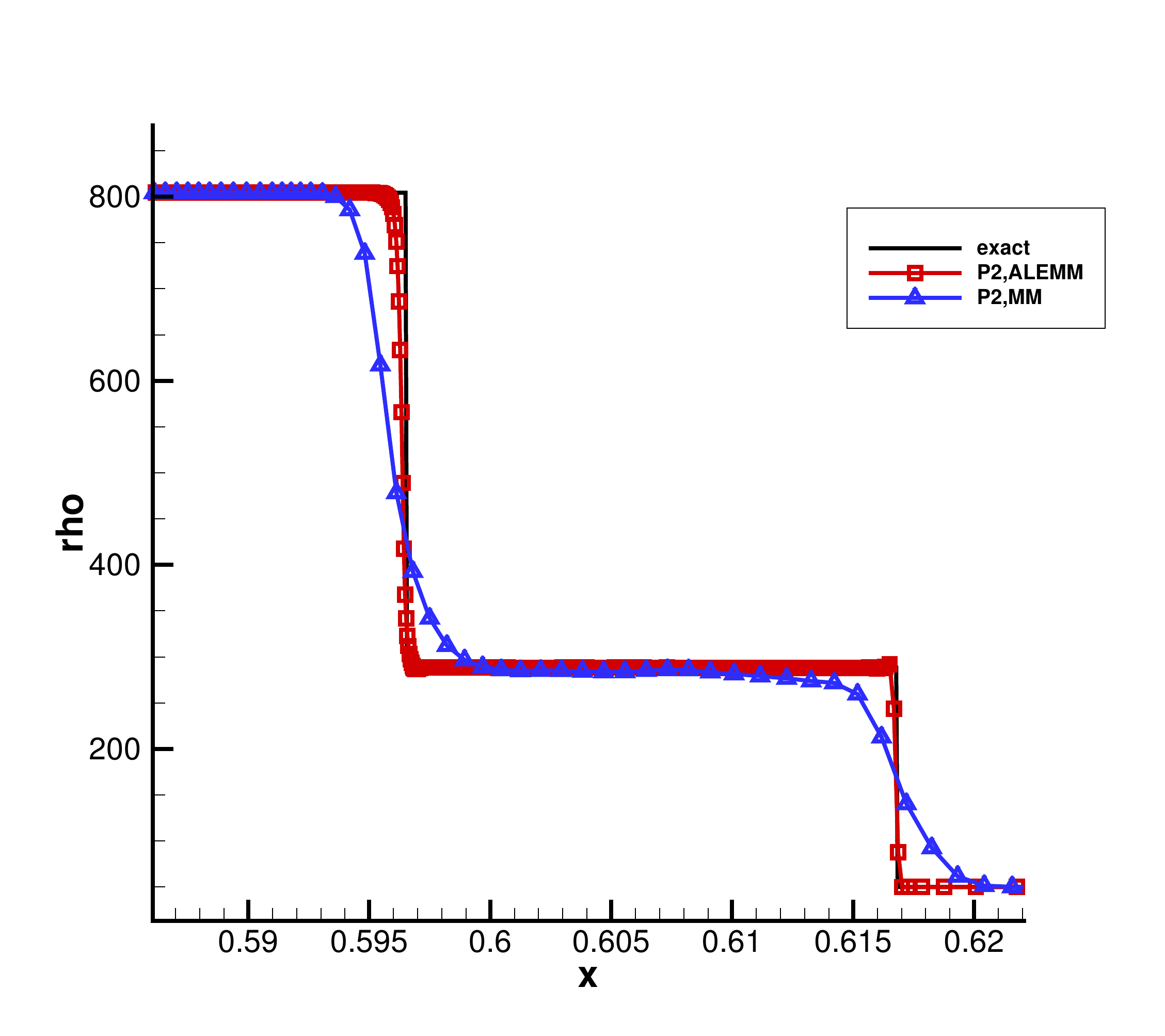}}
\caption{Example~\ref{examgasliquid} The DG-ALE method with moving mesh is compared with the moving mesh DG method
(i.e., the DG-ALE method without using the Lagrangian mesh). $N=2000$.}
\label{figgasliquidmmmm}
\end{figure}
%%%%%%%%%%%%%%%%%end
}
\end{exam}

\subsection{Two-dimensional examples}

\begin{figure}[hbtp]
              \begin{center}
              {\includegraphics[scale=0.3]{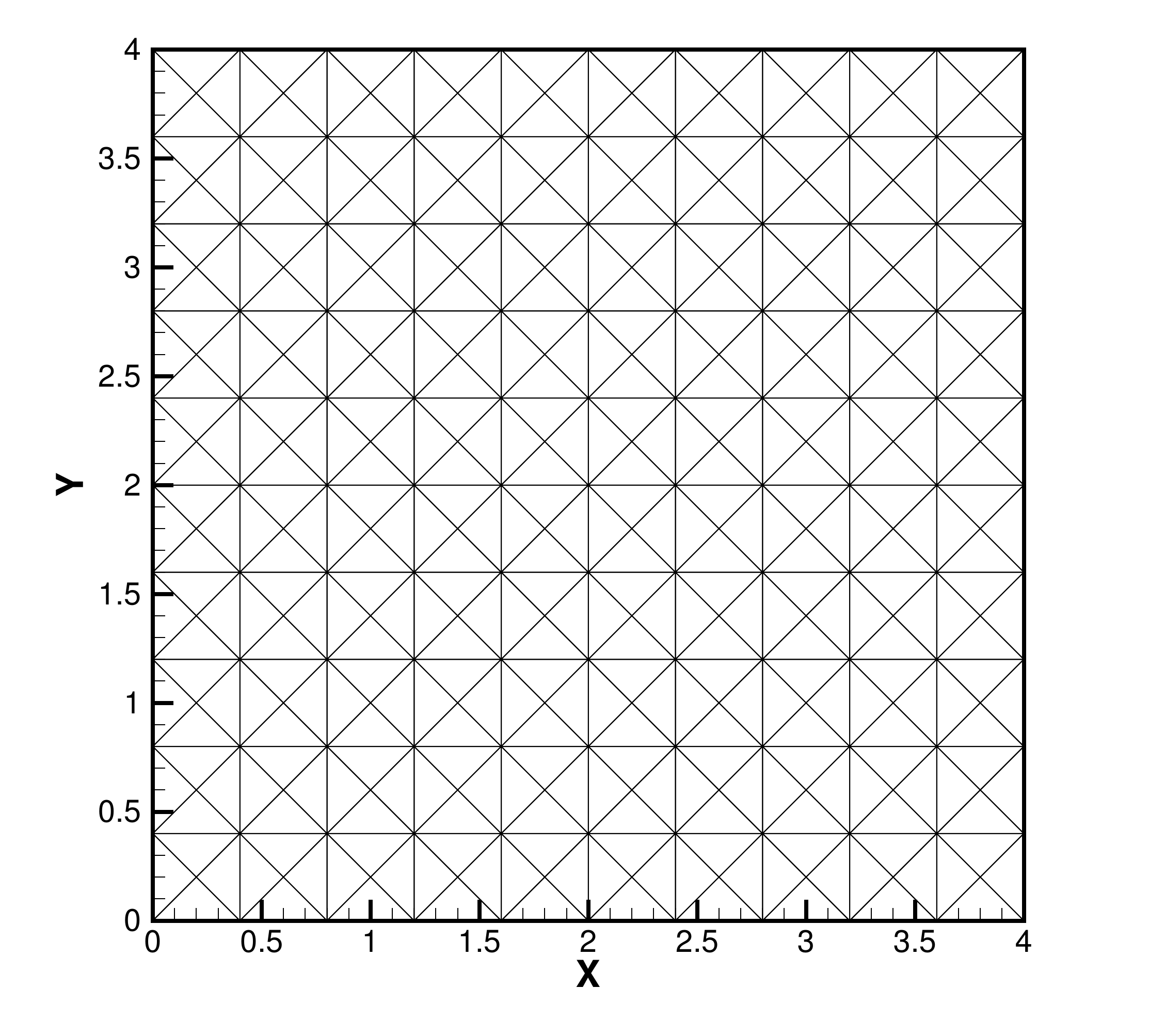}}
                \caption{An initial triangular mesh used in two dimensional computation. A moving mesh associated with
                this mesh is denoted by $N = 10\times 10 \times 4$.}
                \label{sample}
                \end{center}
                \end{figure}

            For two-dimensional examples, an initial triangular mesh is obtained by dividing any rectangular element into
            four triangular elements; see Fig. \ref{sample} for an example. A moving mesh associated with a mesh
            in Fig. \ref{sample} will be denoted by $N = 10\times 10 \times 4$.
            
\begin{exam}{\em
\label{examorder2d}
To verify the convergence order of the DG-ALE method in two dimensions, we consider a sine-wave problem and choose the parameters 
$\gamma_1=1.4$, $\gamma_2=1.9$, $B_1=1$, $B_2=0$.
The initial condition is given by
\begin{align*}
\rho(x,y,0)&=1+0.2\sin(\pi (x+y)),\quad U(x,y,0)=1,\quad V(x,y,0)=1, \\
P(x,y,0)&=1,\quad Y(x,y,0)=0.5+0.5\sin(\pi (x+y)) .
\end{align*}
Periodic boundary conditions are used in both directions. We take the computational domain as $(0,2)\times(0,2)$. 
 
 We compute the solution up to $t= 1$. 
 The error in the density is listed in Table~\ref{exorder2d}, which shows the second-order convergence for $P^1$ elements
 and the third-order convergence for $P^2$ elements for the DG-ALE method.

\begin{table}
\caption{Example~\ref{examorder2d}: The error in the density at $t=1$.}
\renewcommand{\multirowsetup}{\centering}
\begin{center}
\begin{tabular}{|c|c|c|c|c|c|c|c|c|c|c|c|c|}
\hline
$k$ & $N$  &$4\times 4\times 4$      & $8\times 8\times 4$ & $16\times16\times 4$ & $32\times32\times 4$ &$64\times64\times 4$&$128\times128\times 4$ \\
\hline
\multirow{6}{1cm}{1}
%\hline
 & $L^1$ &1.498e-2      & 3.159e-3 & 7.358e-4       & 1.767e-4  & 4.454e-5& 1.140e-5 \\
%\hline
 &  Order      & \quad  & 2.245           & 2.102          &  2.058        & 1.988 & 1.966\\
 %\hline
 &$L^2$   & 1.959e-2      & 4.149e-3  & 9.762e-4    &  2.325e-4  & 5.966e-5 & 1.553e-5\\
%\hline
 &Order   & \quad           & 2.239         & 2.088       &   2.070      & 1.962& 1.942 \\
%\hline
 &$L_{\infty}$ & 5.348e-2 & 1.540e-2 & 3.419e-3  &  8.166e-4  & 2.135e-4 & 5.421e-5\\
%\hline
 &Order           & \quad     & 1.796       & 2.171     &   2.066     & 1.936   & 1.977  \\
 \hline
 \multirow{6}{1cm}{2}

 & $L^1$  & 2.057e-3      &3.360e-4 & 5.516e-5 &  1.112e-5 & 1.685e-6 & 1.876e-7 \\
%\hline
 &  Order      & \quad & 2.614    & 2.607     & 2.311    & 2.722  &3.168 \\
% \hline
 &$L^2$  & 2.709e-3       & 4.020e-4 & 6.957e-5 &  1.511e-5 & 2.468e-6& 2.938e-7 \\
%\hline
 &Order       & \quad &2.752     & 2.531       & 2.203    & 2.614& 3.070\\
%\hline
 &$L_{\infty}$ & 8.688e-3 & 1.102e-3 & 2.707e-4  & 7.037e-5 & 1.163e-5  &1.491e-6\\
%\hline
 &Order    & \quad       &2.980    & 2.025     & 1.944   & 2.597   &2.963 \\
 \hline 
\end{tabular}%\\
\end{center}
\label{exorder2d}
\end{table}

}\end{exam}

%%%%%%%%%%%%%interface only%%%%%%%%
%
\begin{exam}{\em
\label{examinterface2d}
 We test a 2D interface only problem with the initial condition
\begin{equation*}
(\rho,\mu,\nu,P,\gamma,B)=
\begin{cases}
(2,1,1,1,4.4,1),  \; & (x-0.2)^2+(y-0.2)^2\leqslant 0.01\\
(1,1,1,1,1.4,0), \;  &(x-0.2)^2+(y-0.2)^2>0.01.
\end{cases}
\end{equation*}
For this example, nonreflecting boundary conditions are used for the top and bottom boundary and inflow and outflow boundary conditions
are used for the left and right boundary. The computational domain is taken as $(0,1)\times(0,1)$ and the numerical computation is stopped at $t=0.3$.

The results obtained with $N = 100\times 100\times 4$ are plotted in Fig.~\ref{figinterface2dp1}.
One can observe that the quasi-conservative DG-ALE method can preserve the constant pressure and velocity and seems to be free of oscillations
near the material interfaces. The density along the line $y=0.5$ and the final adaptive meshes for both $P^1$-DG and $P^2$-DG
are plotted in Fig.~\ref{figinterface2dtraj}. It is clear that the mesh points are concentrated near the material interface and
the results with $P^2$-DG have better resolution than the ones with $P^1$-DG.
\begin{figure}[hbtp]
\begin{center}
\mbox{%\subfigure[Euler: $N_t=75\times 25\times 4$]
{\includegraphics[width=0.35\textwidth]{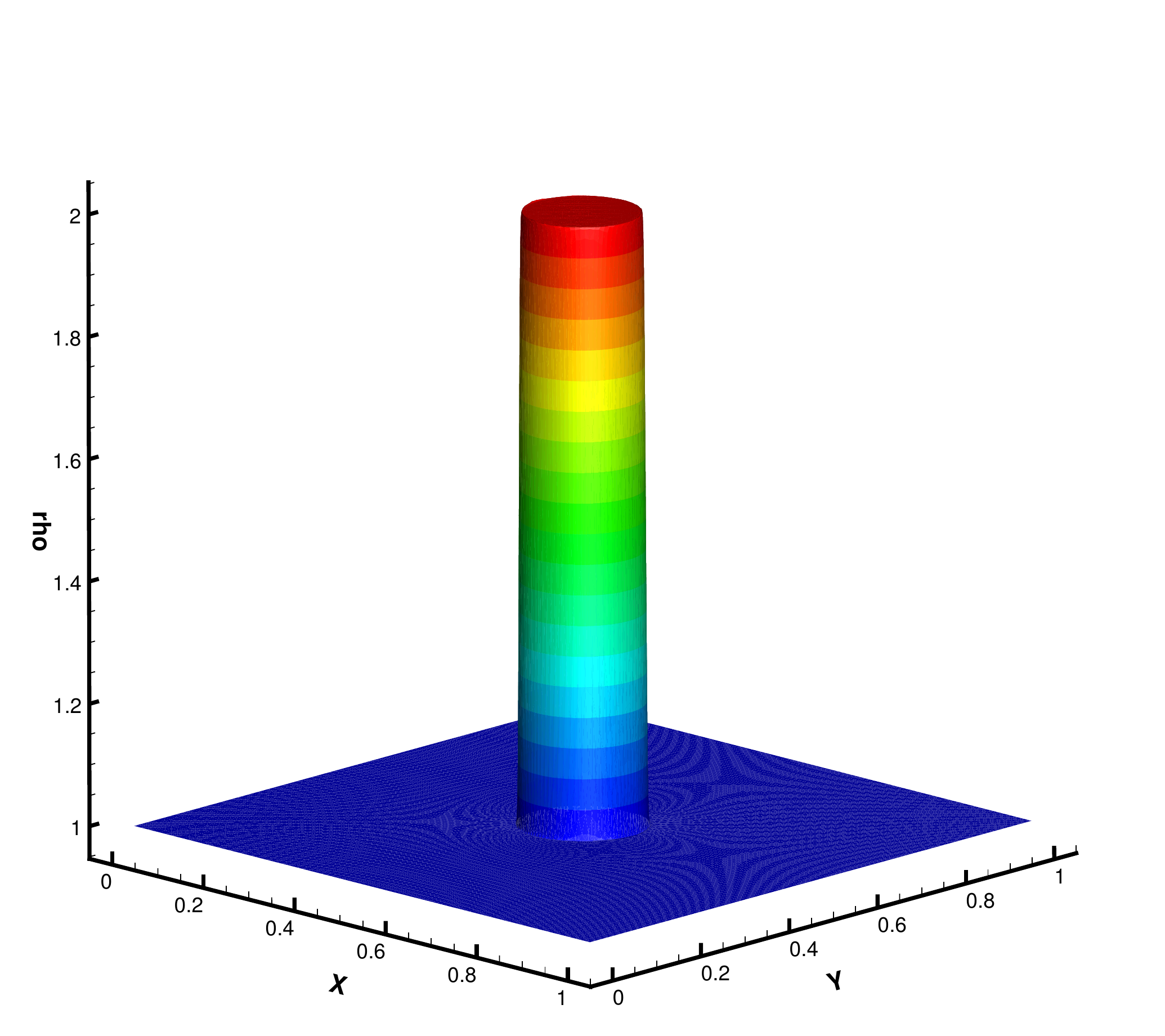}}
%\subfigure[Euler: $N_t=75\times 25\times 4$]
{\includegraphics[width=0.35\textwidth]{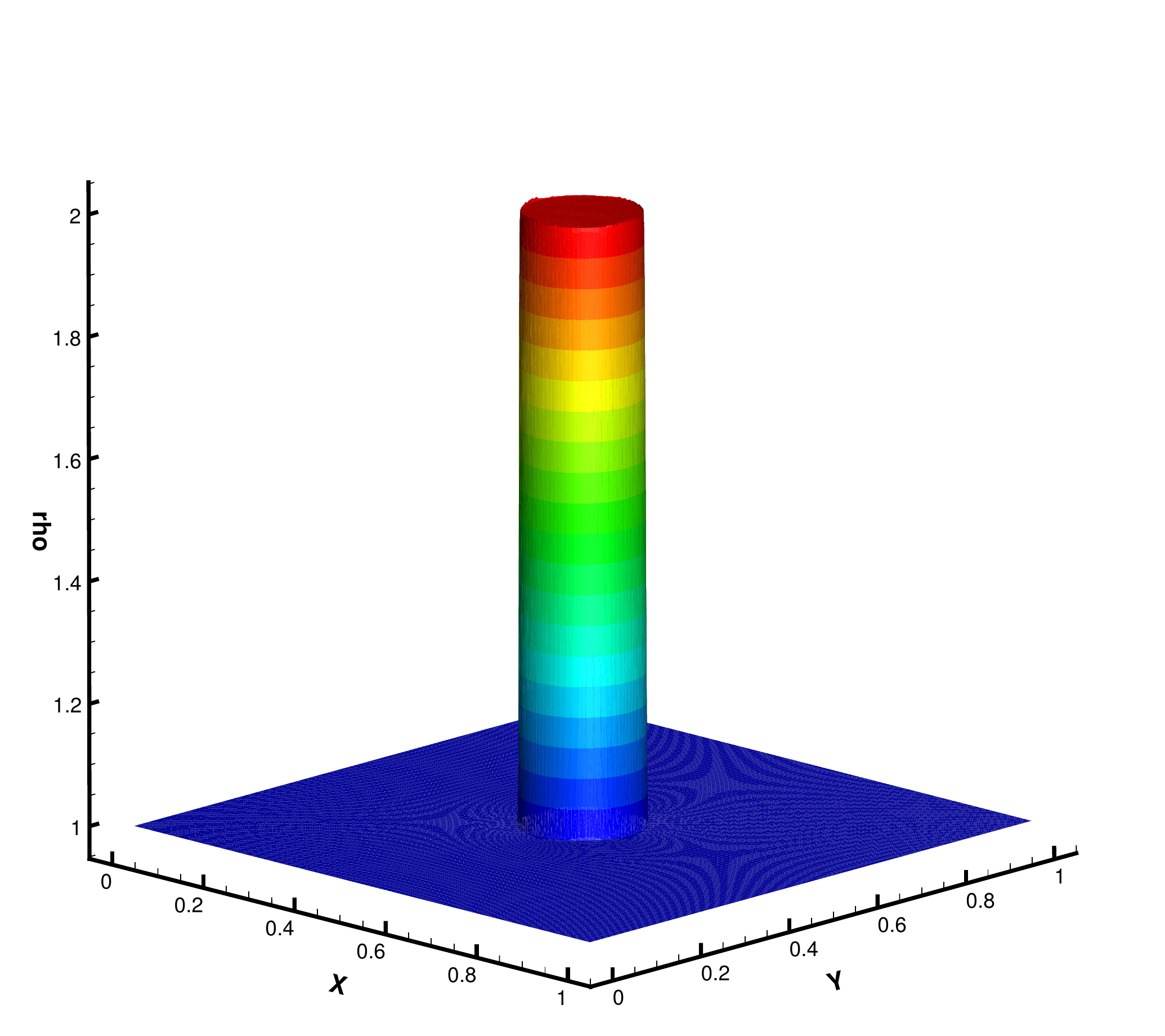}}
}

\mbox{%\subfigure[Euler: $N_t=75\times 25\times 4$]
{\includegraphics[width=0.35\textwidth]{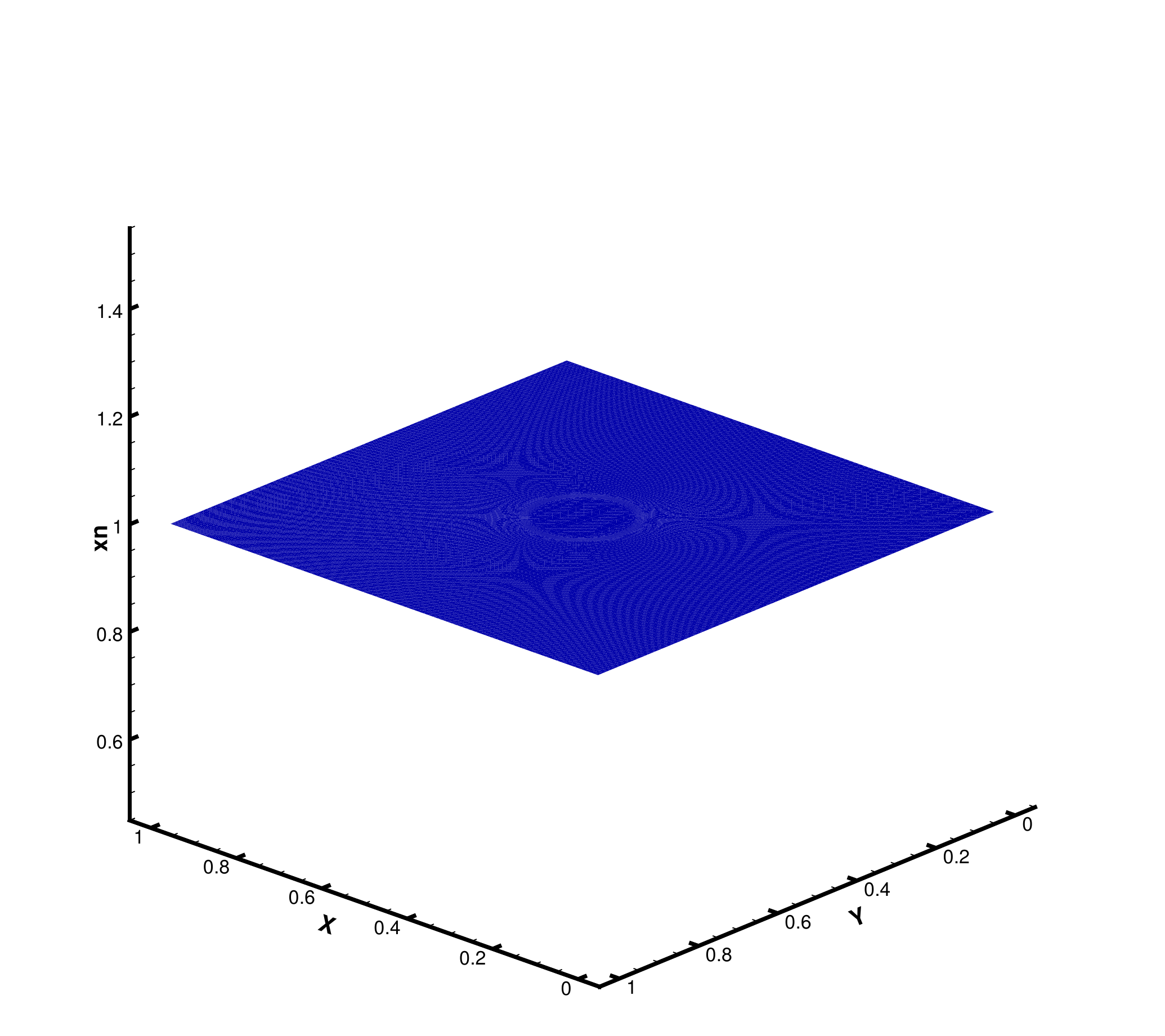}}
%\subfigure[Euler: $N_t=75\times 25\times 4$]
{\includegraphics[width=0.35\textwidth]{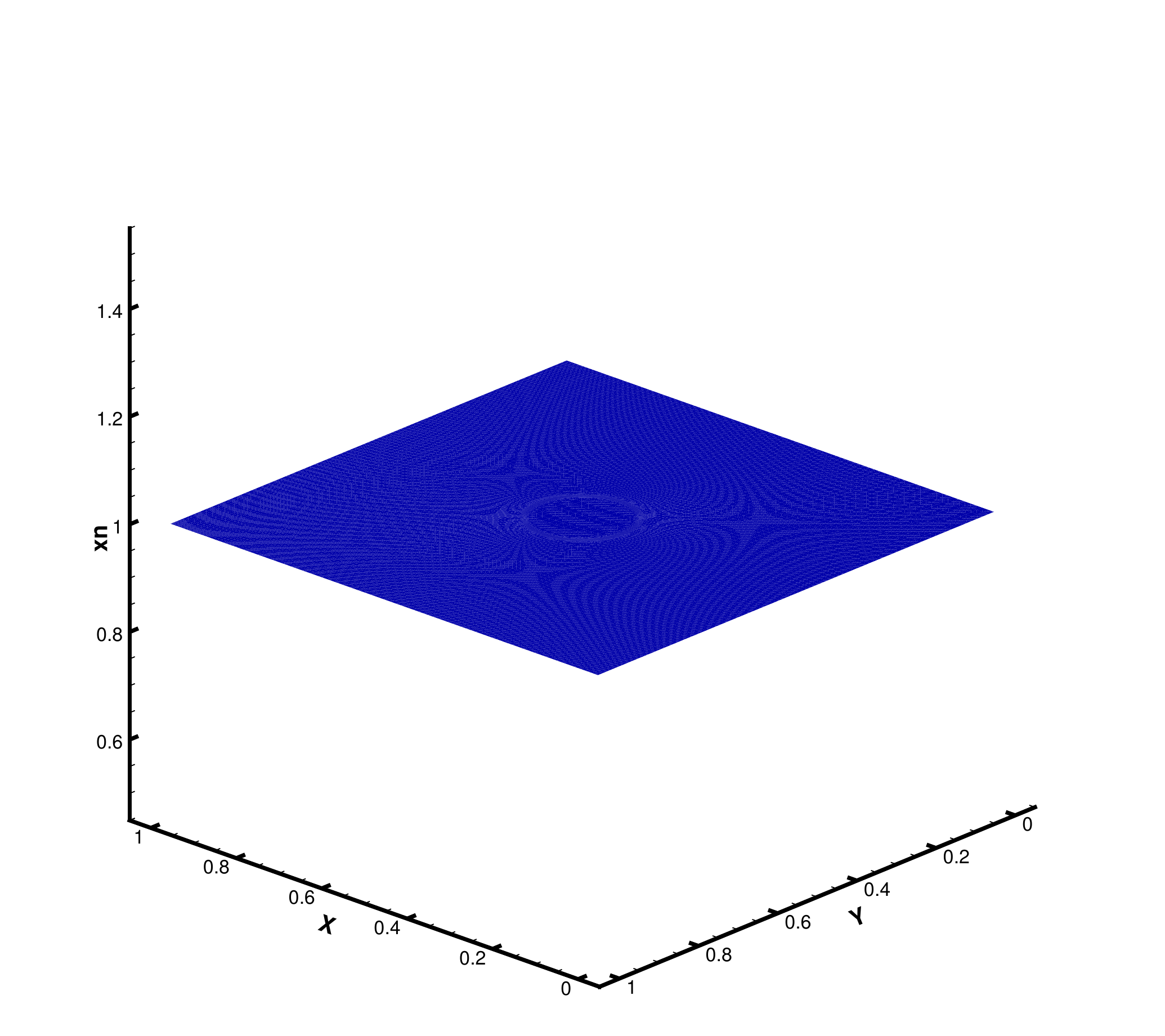}}
}

\mbox{%\subfigure[ALE: $N_t=75\times 25\times 4$]
{\includegraphics[width=0.35\textwidth]{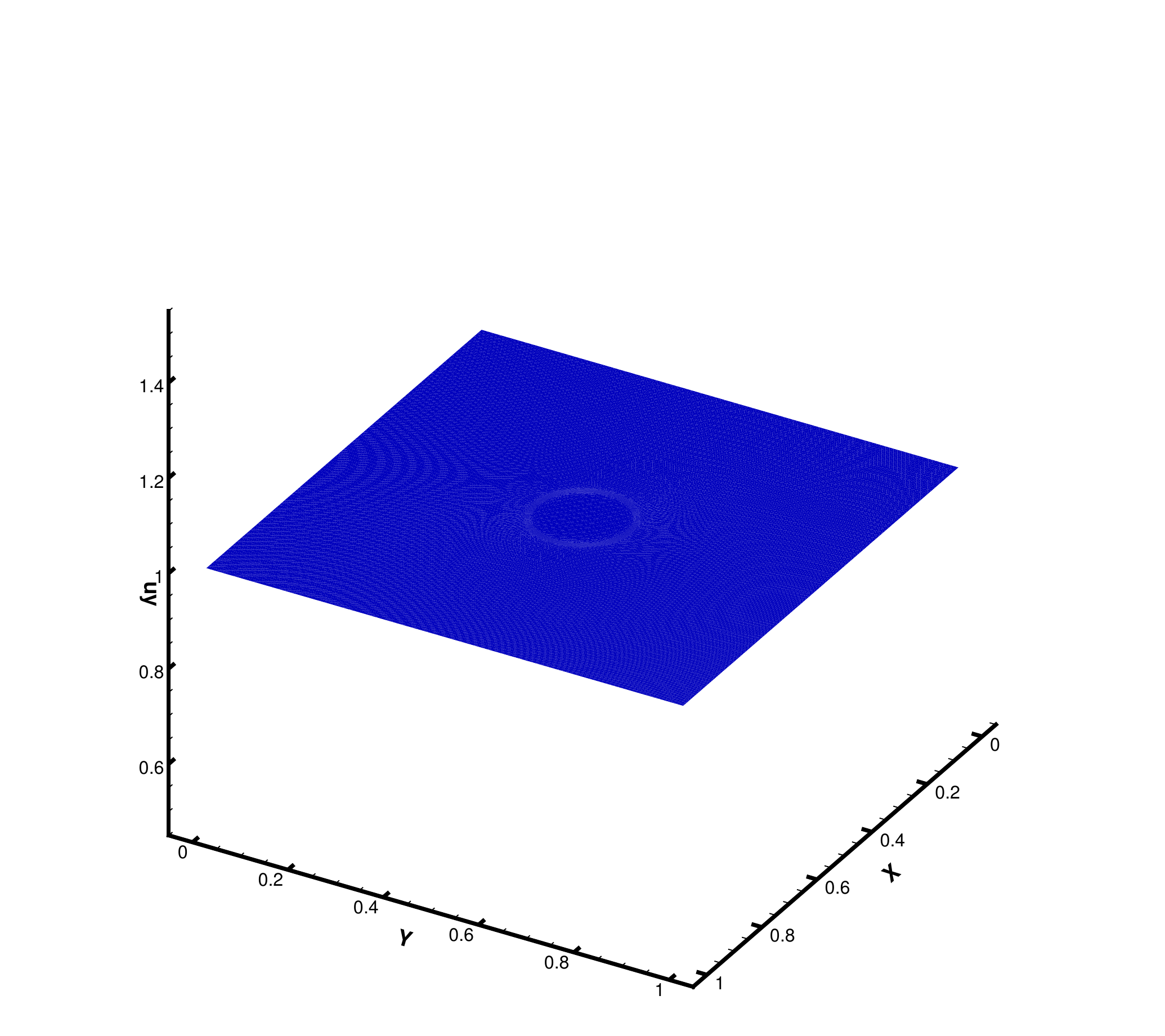}}\quad
%\subfigure[Mesh]
{\includegraphics[width=0.35\textwidth]{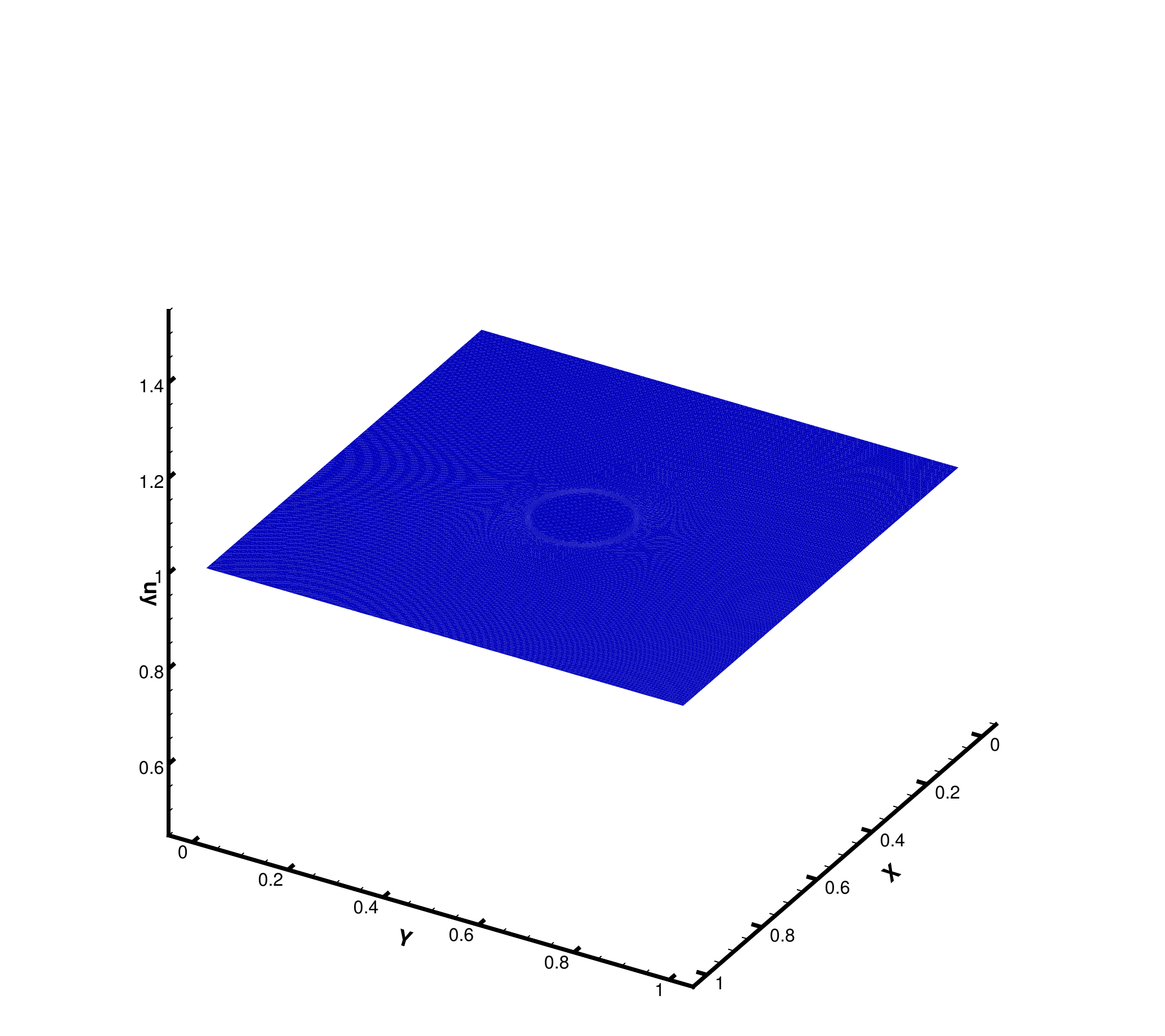}}
}

\mbox{%\subfigure[ALE: $N_t=75\times 25\times 4$]
{\includegraphics[width=0.35\textwidth]{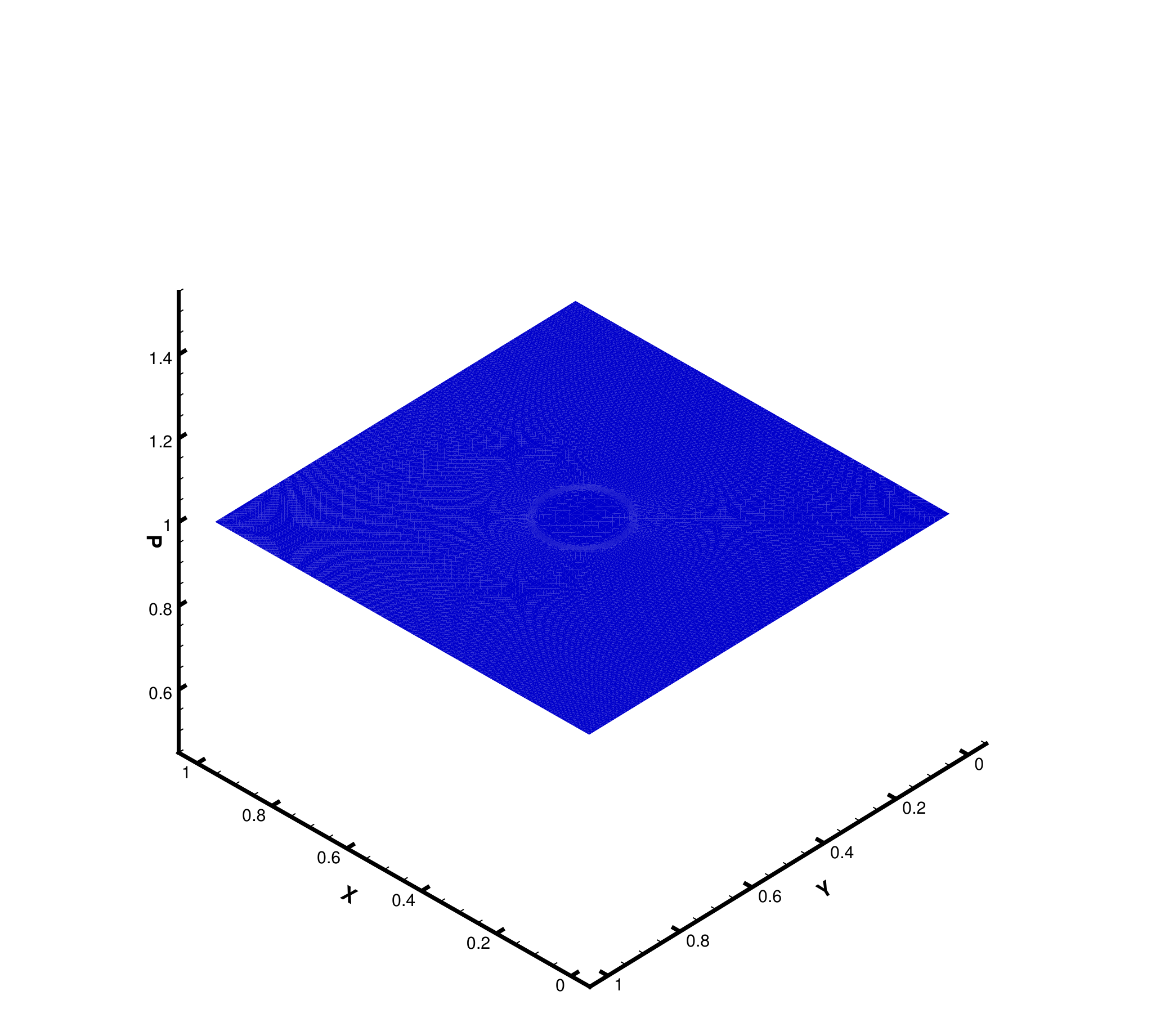}}\quad
%\subfigure[Mesh]
{\includegraphics[width=0.35\textwidth]{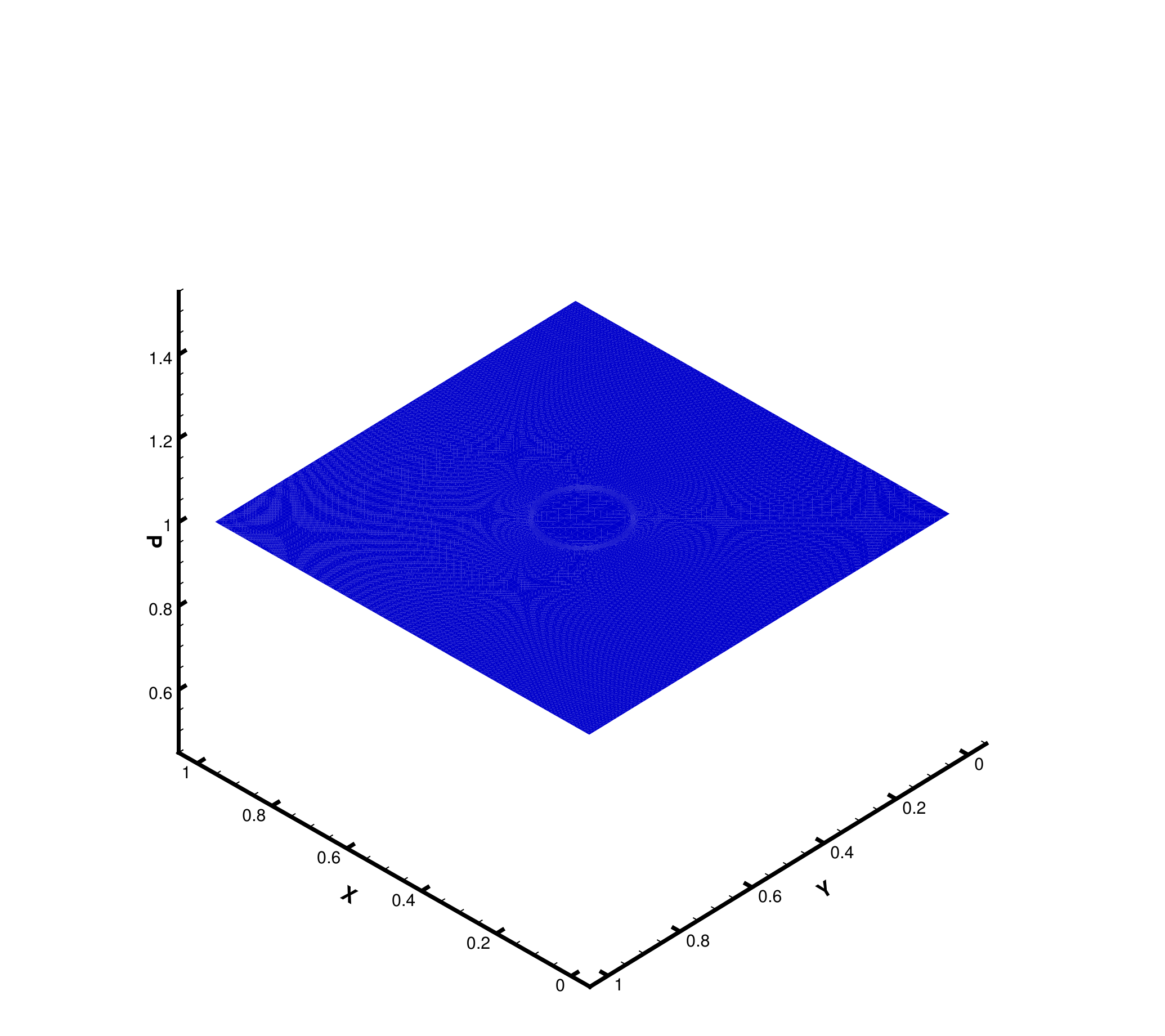}}
}

\caption{Example~\ref{examinterface2d}. From top row to bottom row: density, velocity in the $x$-direction, velocity in the $y$-direction,
and pressure. Left column: $P^1$-DG; Right column: $P^2$-DG. $N = 100\times 100\times 4$}
\label{figinterface2dp1}
\end{center}
\end{figure}

\begin{figure}[hbtp]
\begin{center}

\mbox{\subfigure[The density along the line $y=0.5$]
{\includegraphics[width=0.4\textwidth]{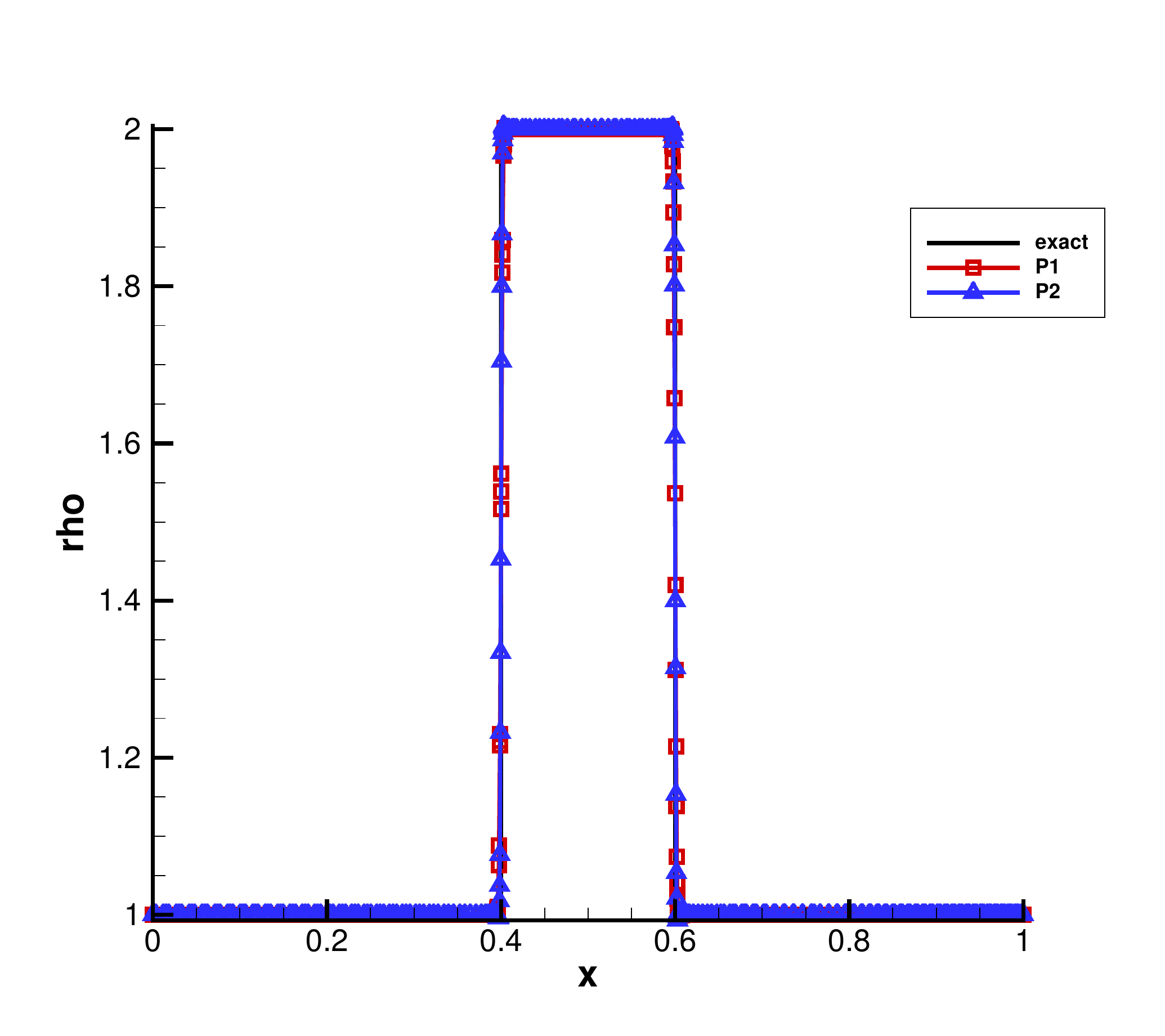}}\quad
\subfigure[Zoom of (a)]
{\includegraphics[width=0.4\textwidth]{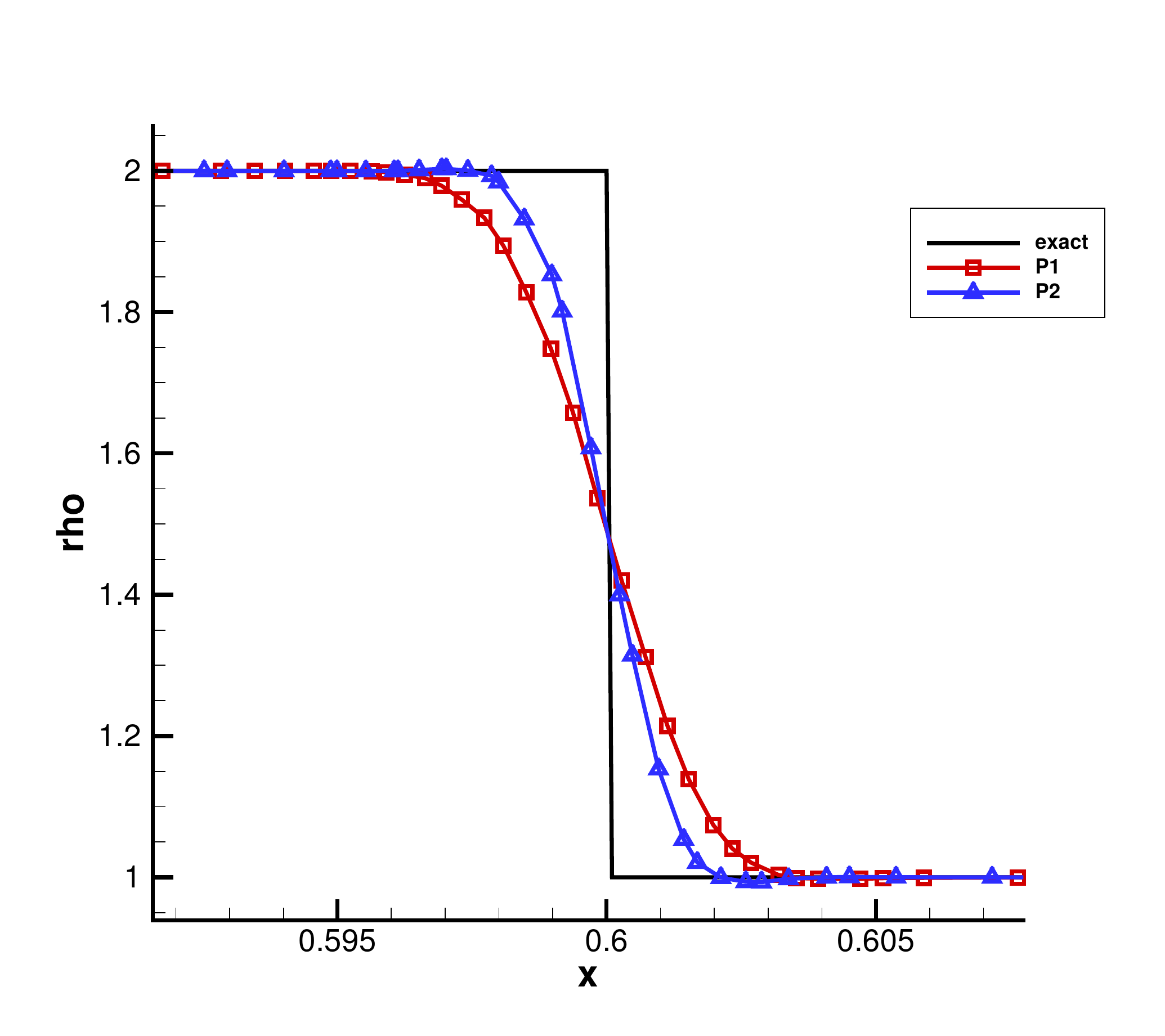}}
}

\mbox{\subfigure[Mesh, $P^1$-DG]
{\includegraphics[width=0.4\textwidth]{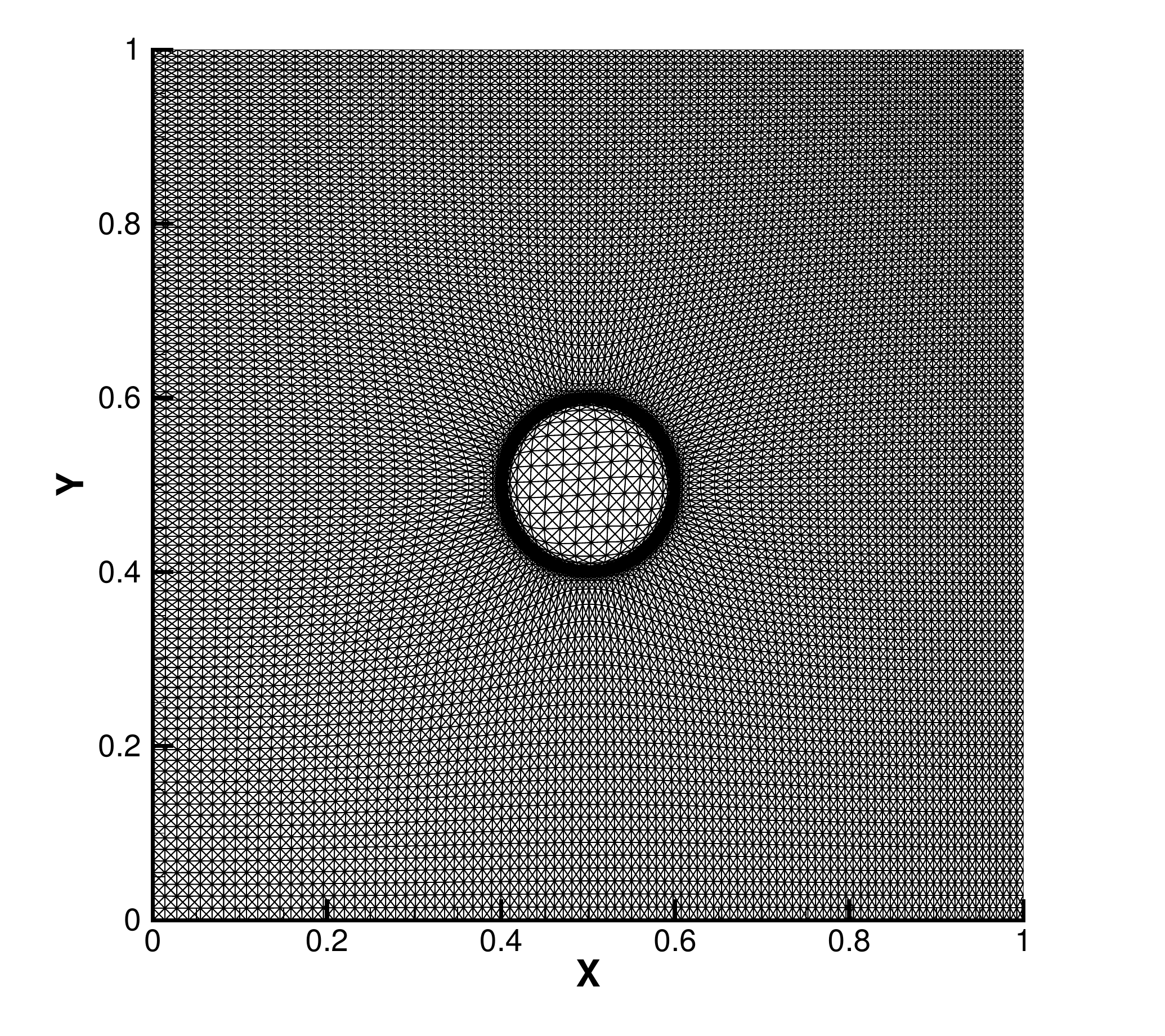}}\quad
\subfigure[Mesh, $P^2$-DG]
{\includegraphics[width=0.4\textwidth]{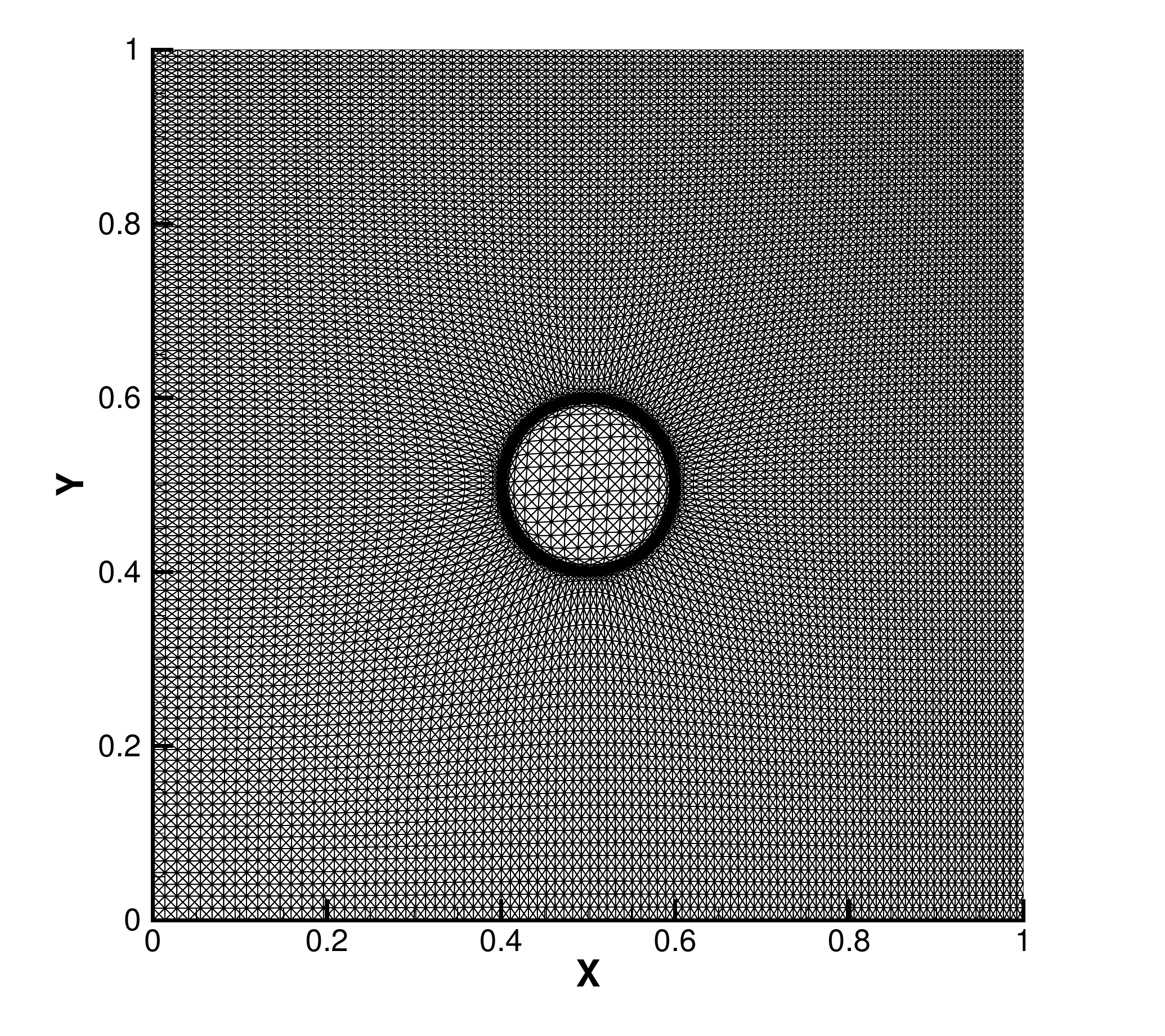}}
}

\caption{Example~\ref{examinterface2d}. The density along $y=0.5$ and the final meshes of $N = 100\times 100\times 4$. }
\label{figinterface2dtraj}
\end{center}
\end{figure}
}
\end{exam}

%%%%%%%%%%%%%%%%%%gas-gas interface%%%%%%55
\begin{exam}{\em
\label{examgasgas}
We consider an air shock impacting on a helium bubble \cite{qiu2008}. 
The domain of this problem is $(-3,4)\times(-3,3)$ and the initial condition is
\begin{equation*}
(\rho,U,V,P,\gamma,B)=
\begin{cases}
(0.138,0,0,1,\frac{5}{3},0),  \; & x^2+y^2 \leqslant 1\\
(1.3764,0.394,0,1.5698,1.4,0), \;   & x< -1.2\\
(1,0,0,1,1.4,0), \; & x\geqslant -1.2  \text{ and } x^2+y^2 > 1.
\end{cases}
\end{equation*}
Reflective boundary conditions are used for the top and bottom boundary and inflow and outflow boundary conditions are adopted
for the left and right boundary. The final time is $t=4$. The parameters $\beta_i$'s in \eqref{ent} are all taken as 1 in this example.

We plot contours of the density and the adaptive mesh of $N=70\times 60\times 4$ at $t=0.5$, 1, 2, and 4 in Figs.~\ref{figgasgas}
and \ref{figgasgastraj}. One can see that the mesh points are concentrated at the material interface, an indication that the combination
of the Lagrangian meshing and the MMPDE moving mesh method works well.
Moreover, $P^2$-DG produces sharper contours of the density near the material interface than $P^1$-DG.

\begin{figure}[hbtp]
 \begin{center}
 \mbox{%\subfigure[t=0]
 {\includegraphics[width=0.35\textwidth]{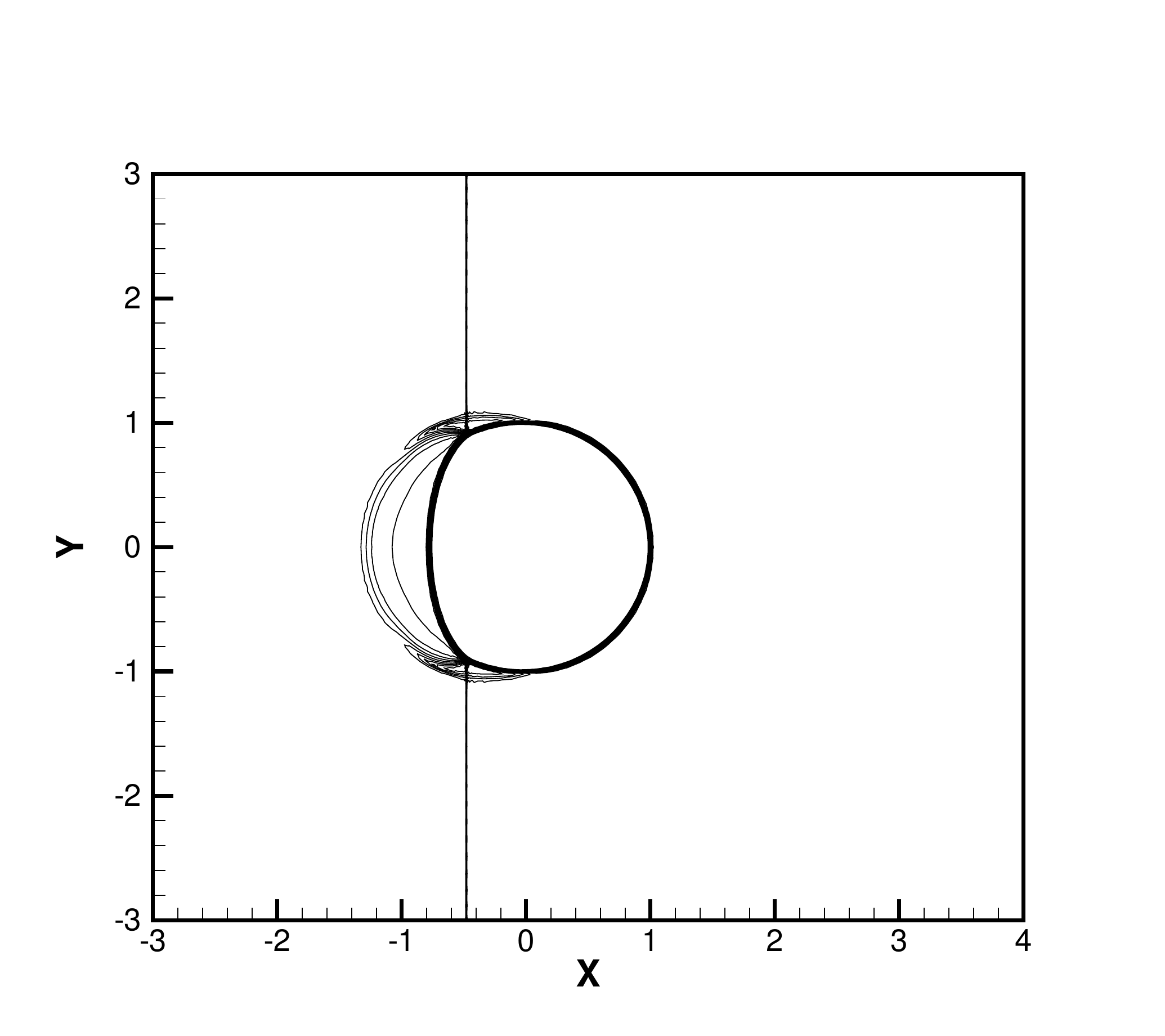}}\quad
  % \subfigure[close view of (a) near shock]
   {\includegraphics[width=0.35\textwidth]{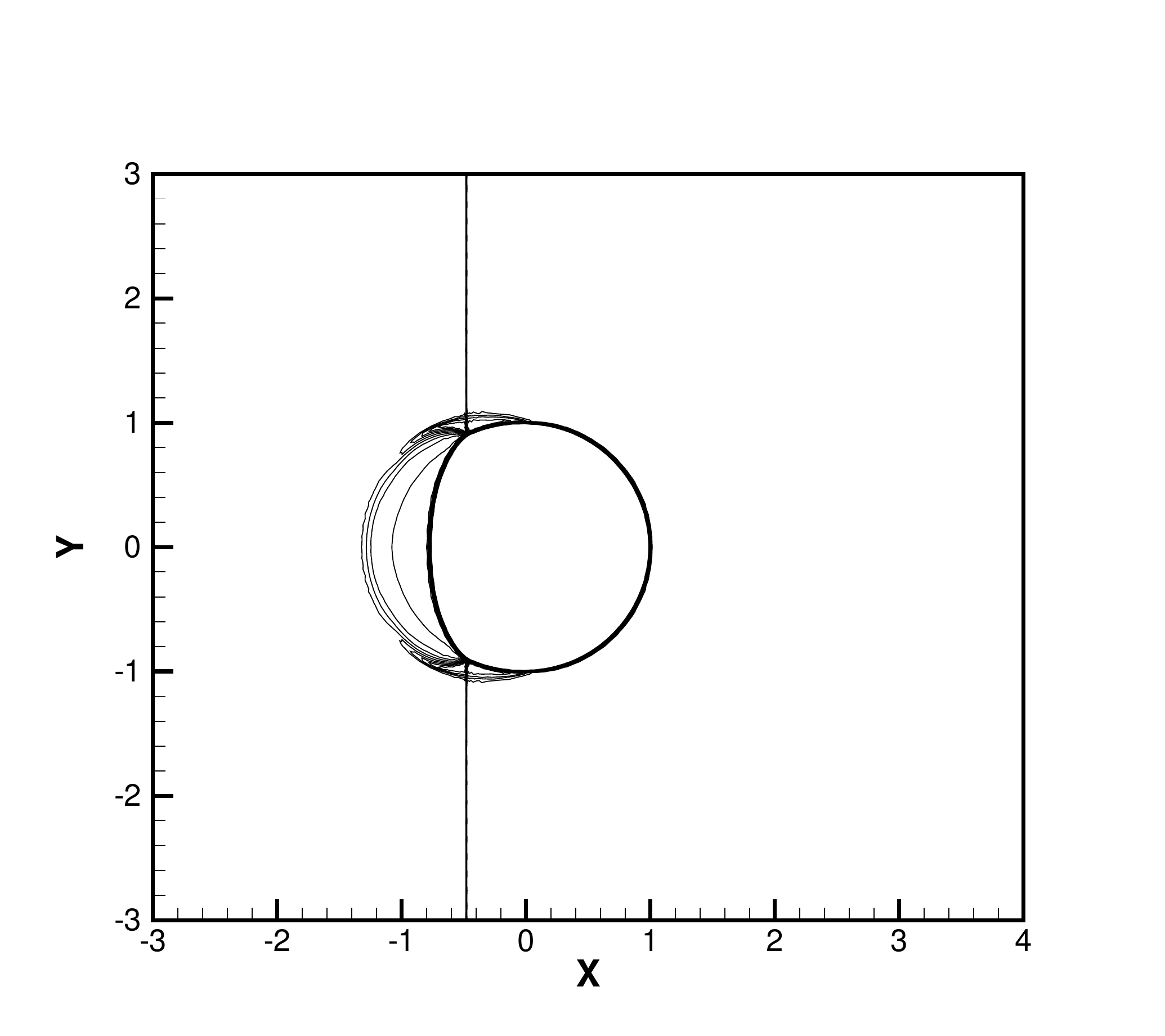}}
   }
   
 \mbox{%\subfigure[t=0]
 {\includegraphics[width=0.35\textwidth]{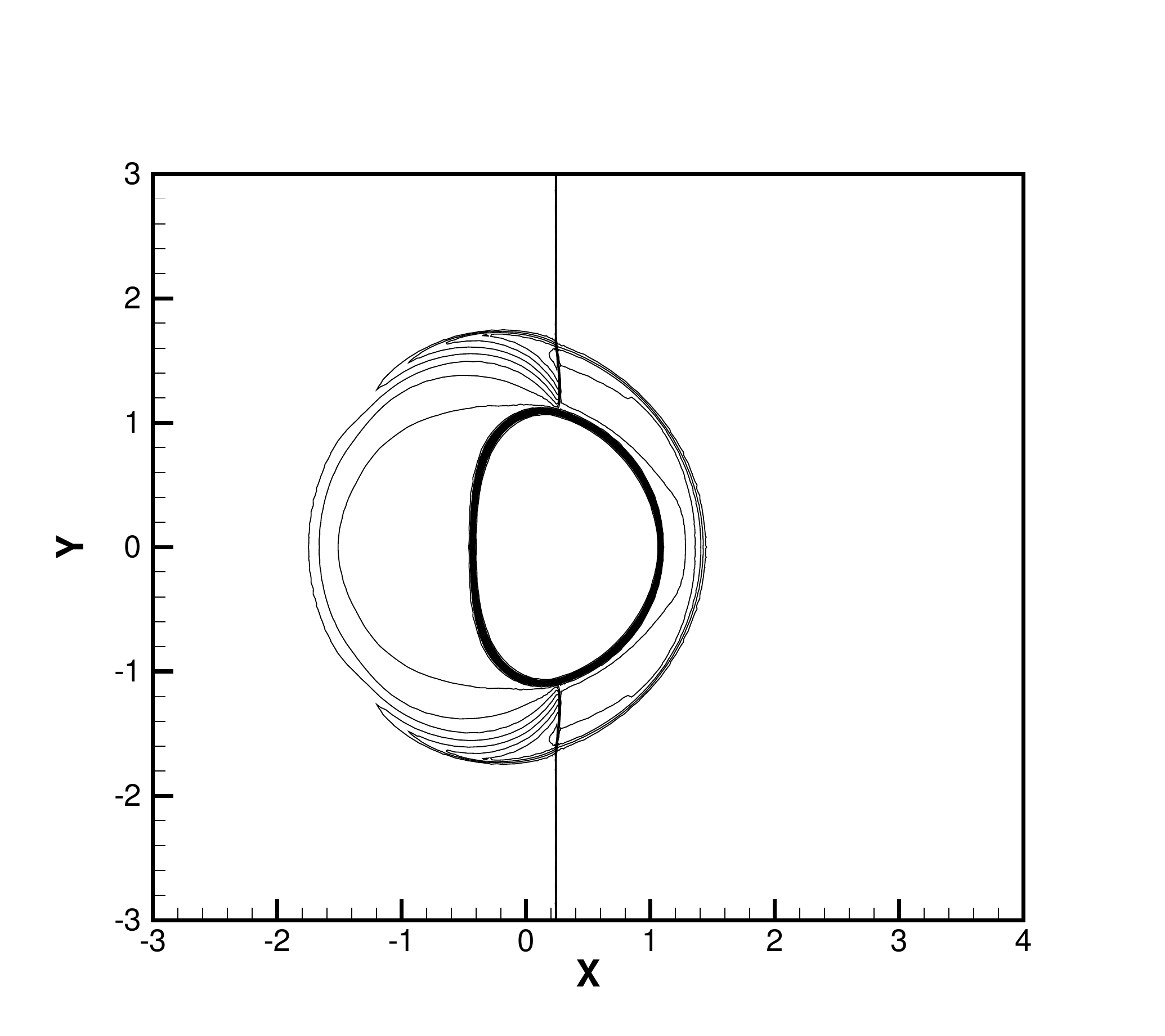}}\quad
  % \subfigure[close view of (a) near shock]
   {\includegraphics[width=0.35\textwidth]{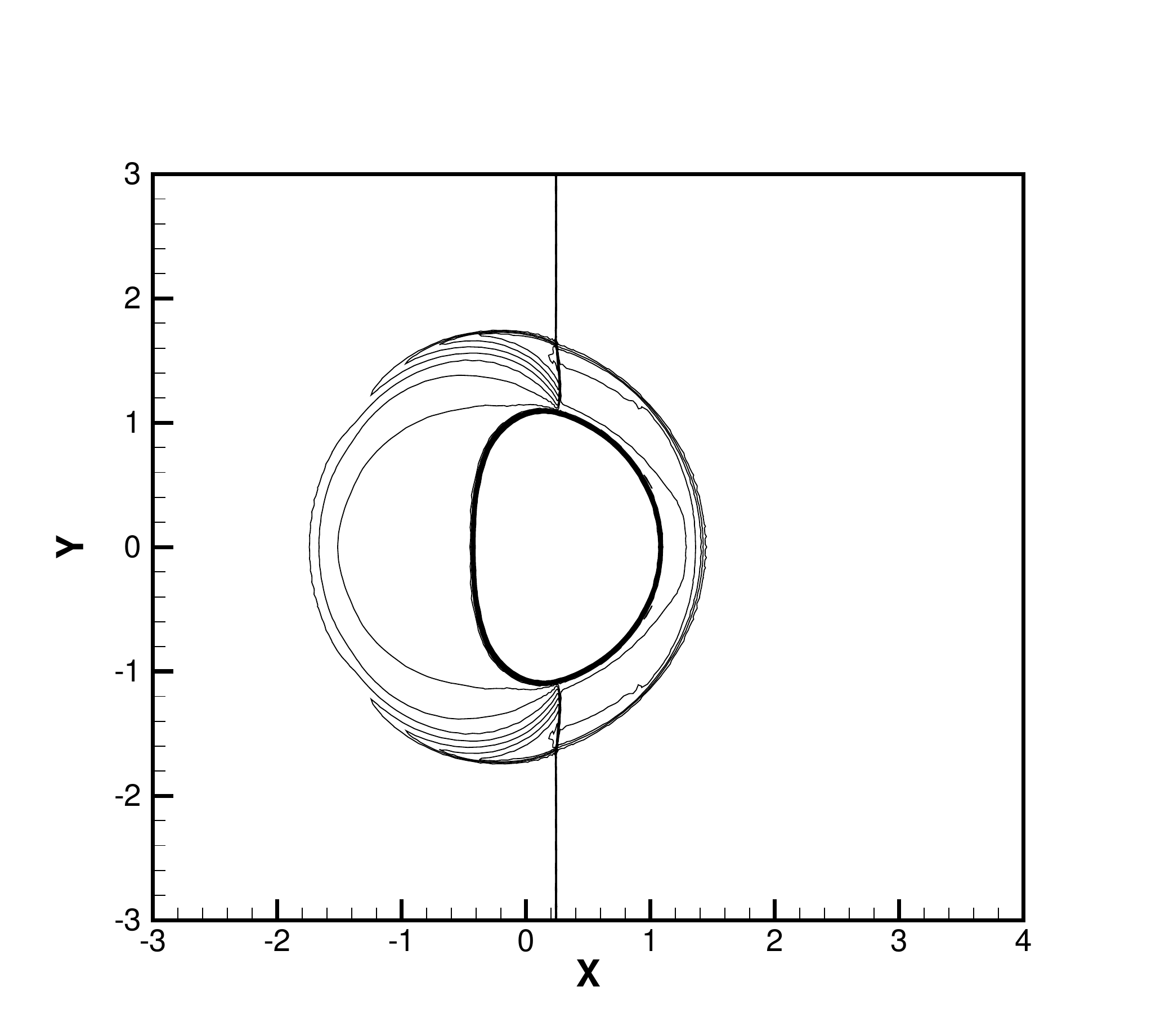}}
   }
 \mbox{%\subfigure[t=0]
 {\includegraphics[width=0.35\textwidth]{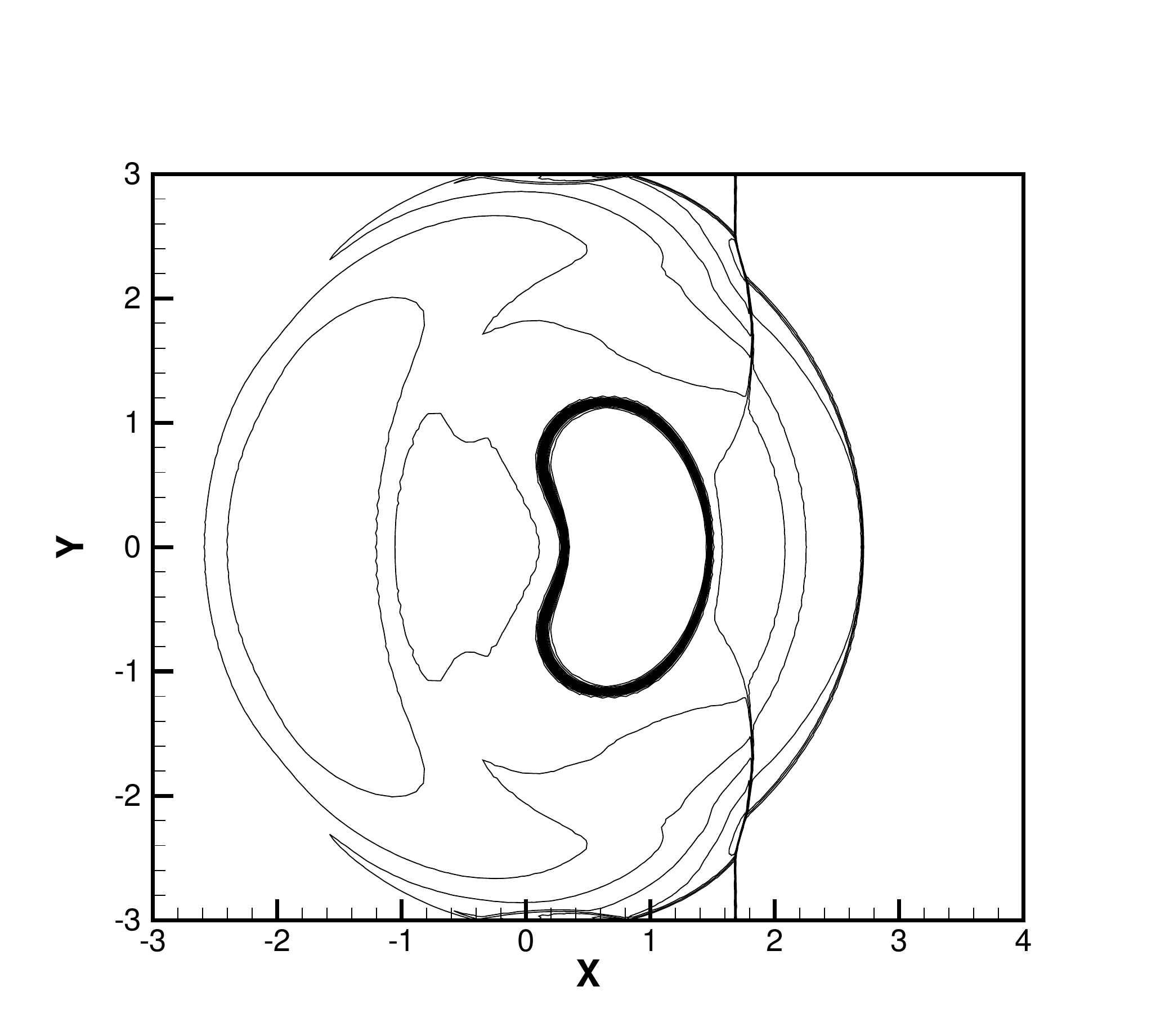}}\quad
  % \subfigure[close view of (a) near shock]
   {\includegraphics[width=0.35\textwidth]{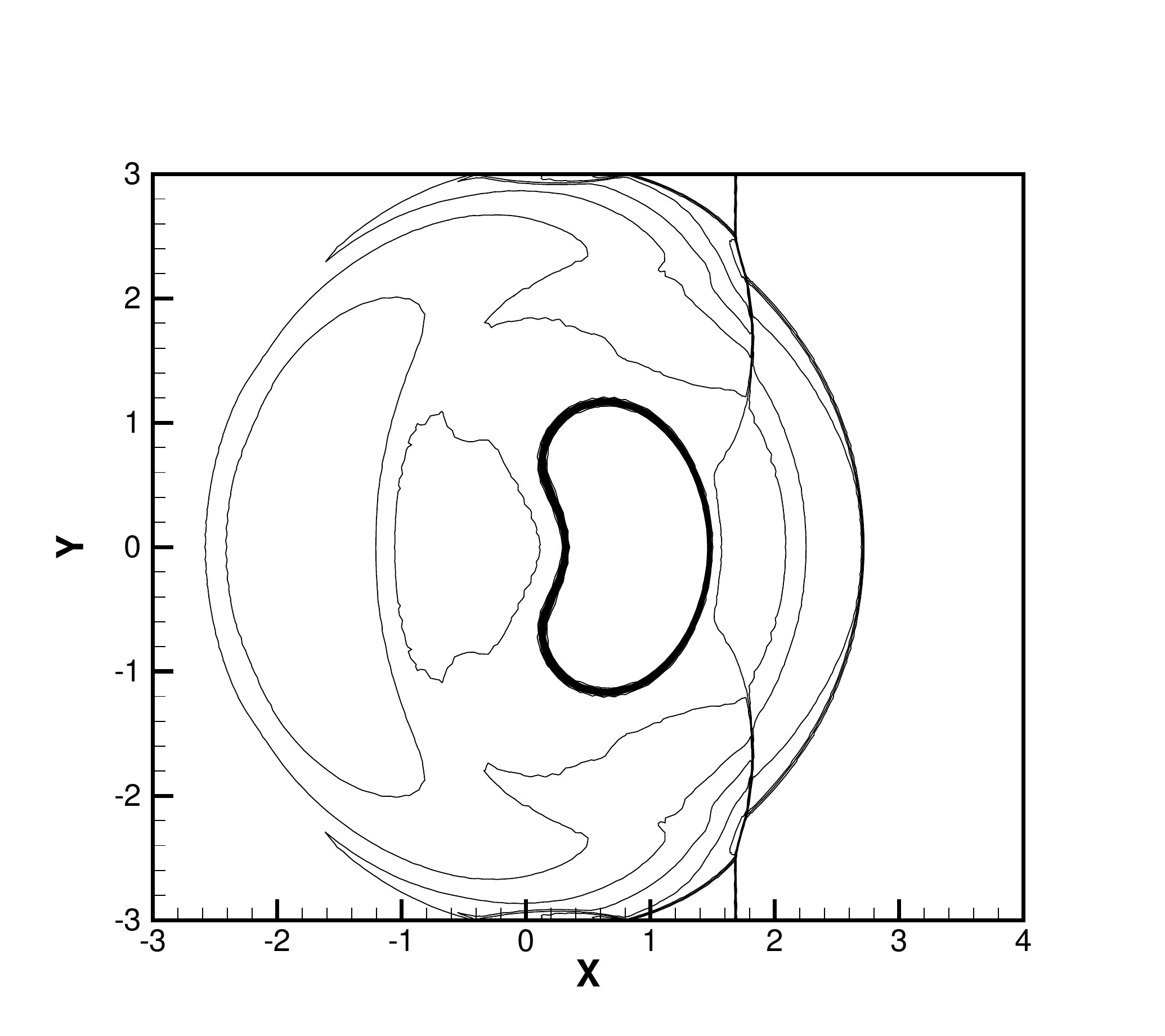}}
   }

\mbox{%\subfigure[t=0]
 {\includegraphics[width=0.35\textwidth]{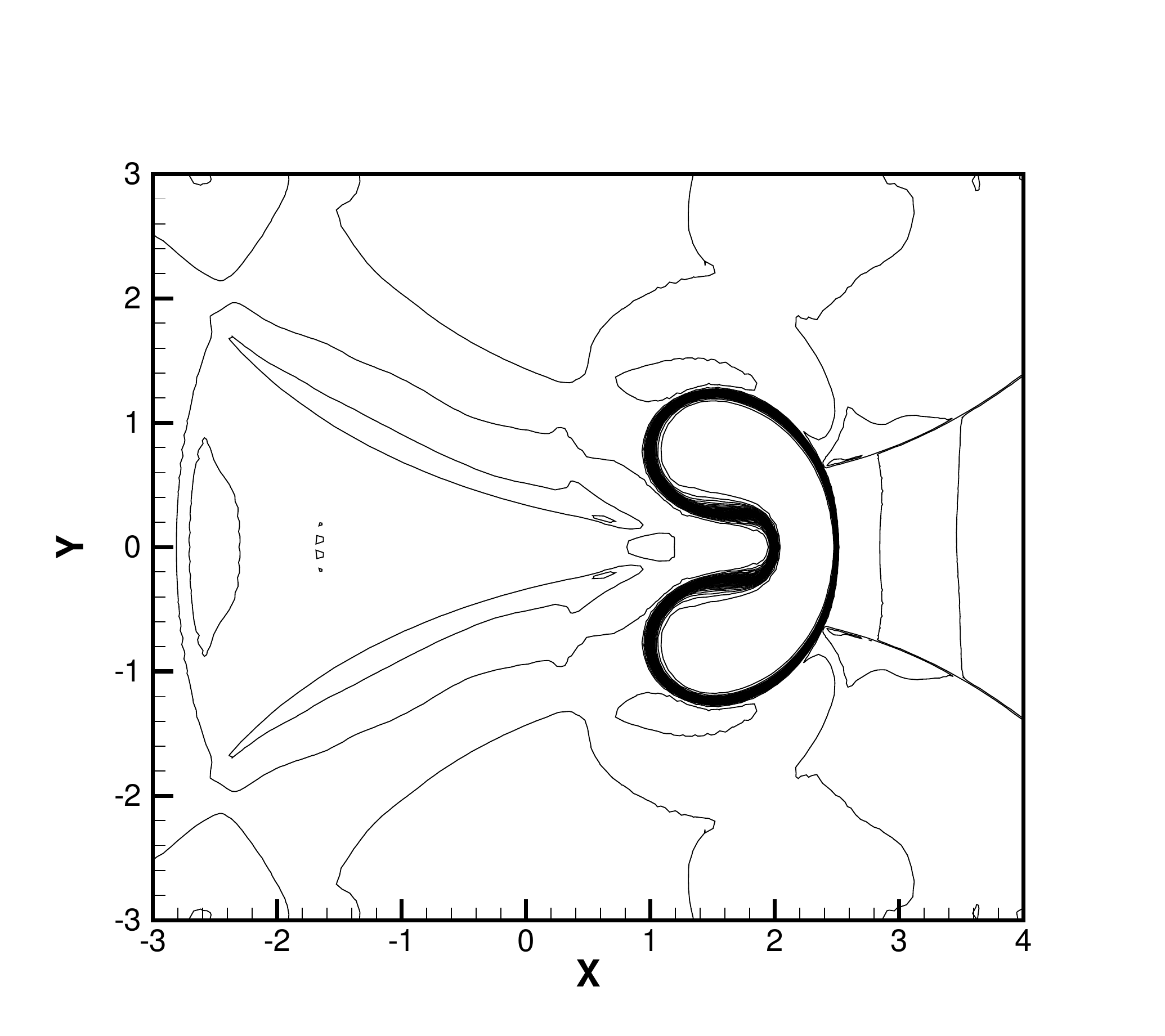}}\quad
  % \subfigure[close view of (a) near shock]
   {\includegraphics[width=0.35\textwidth]{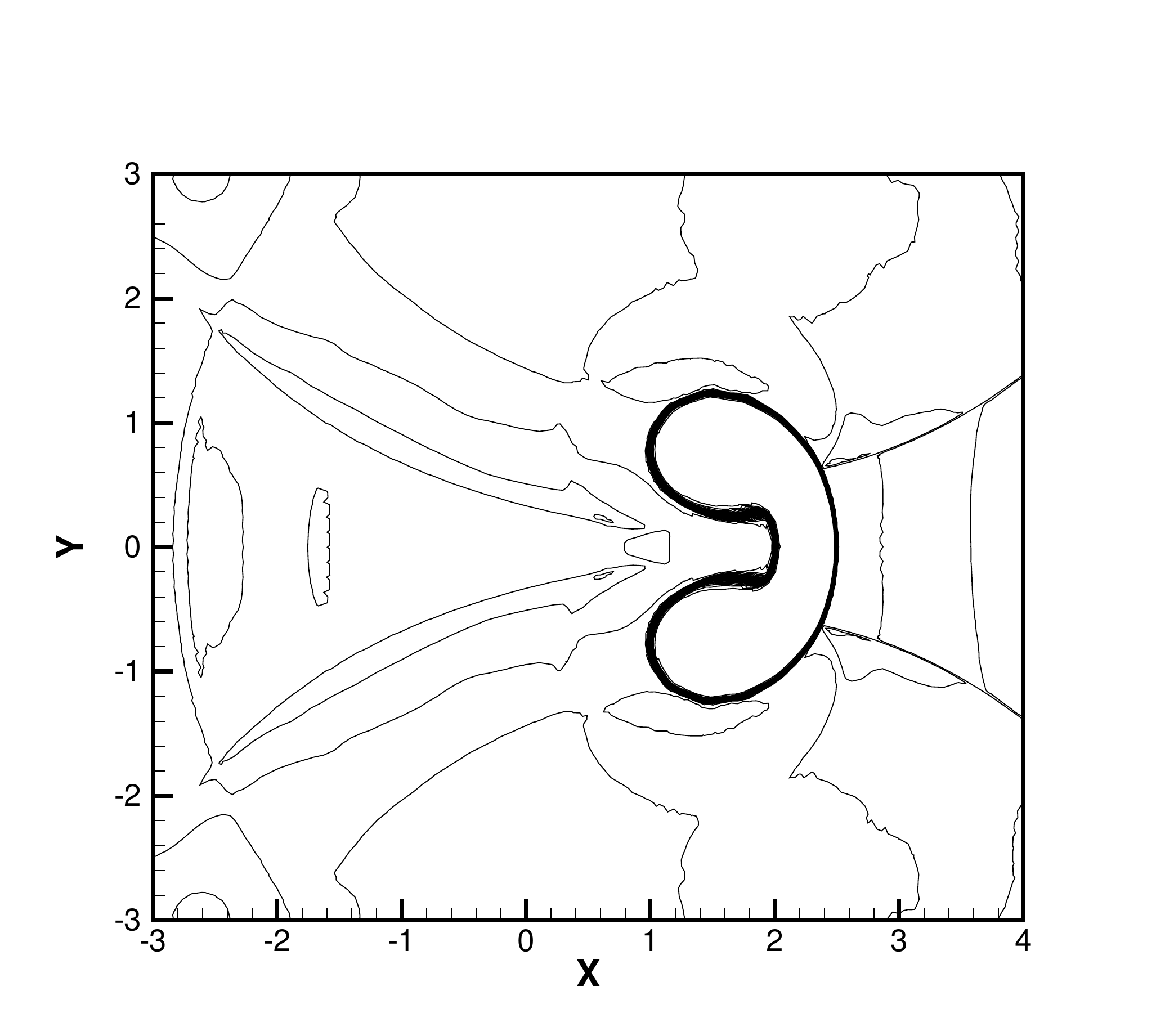}}
   }
   
   \caption{Example~\ref{examgasgas} Density contours at $t=0.5$, 1, 2, and 4 (from top to bottom). Left column: $P^1$-DG;
   Right column: $P^2$-DG. $N=70\times 60\times 4$.}
   \label{figgasgas}
   \end{center}
   \end{figure}
   
   \begin{figure}[hbtp]
 \begin{center}
 \mbox{%\subfigure[t=0]
 {\includegraphics[width=0.35\textwidth]{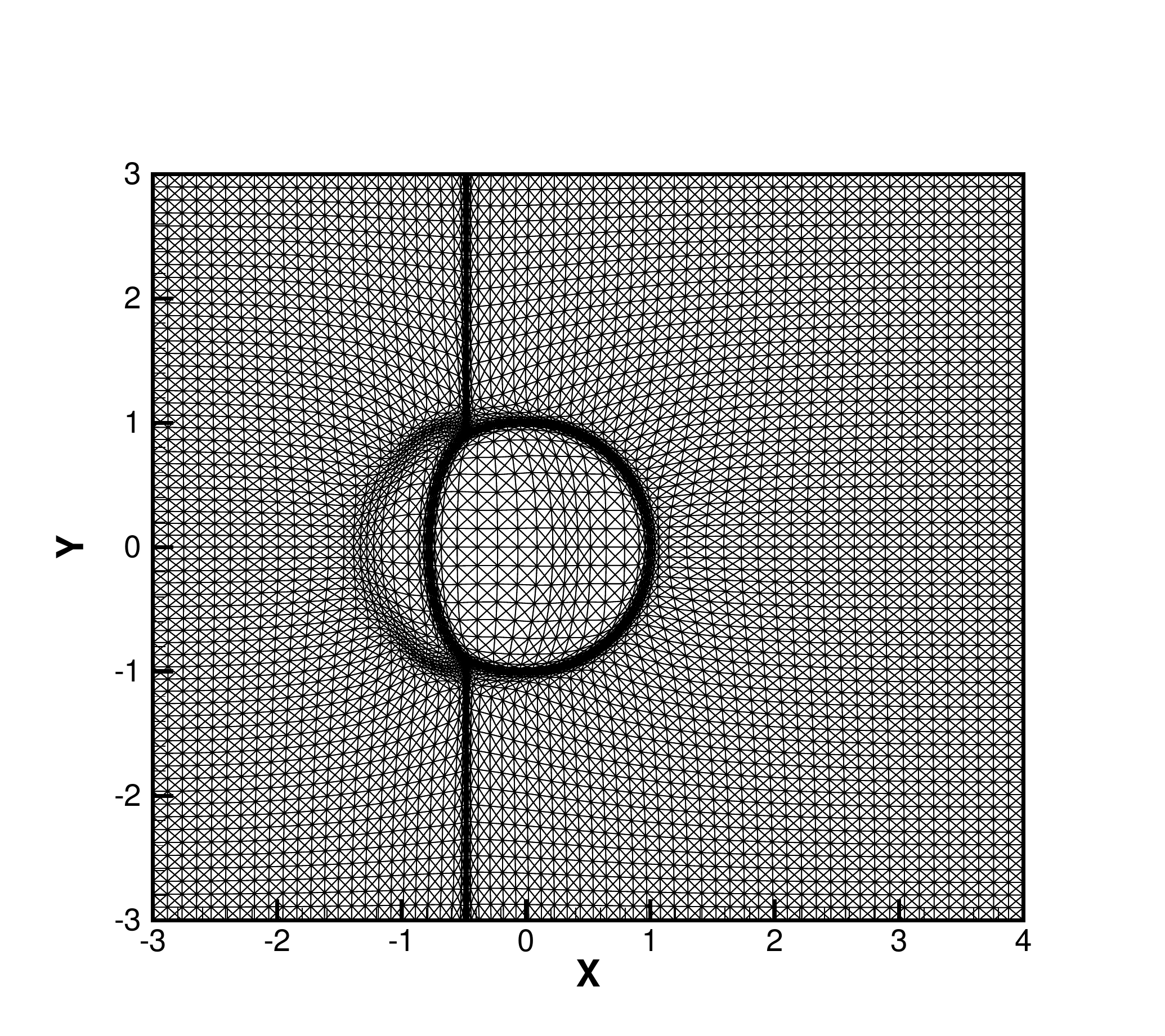}}\quad
  % \subfigure[close view of (a) near shock]
   {\includegraphics[width=0.35\textwidth]{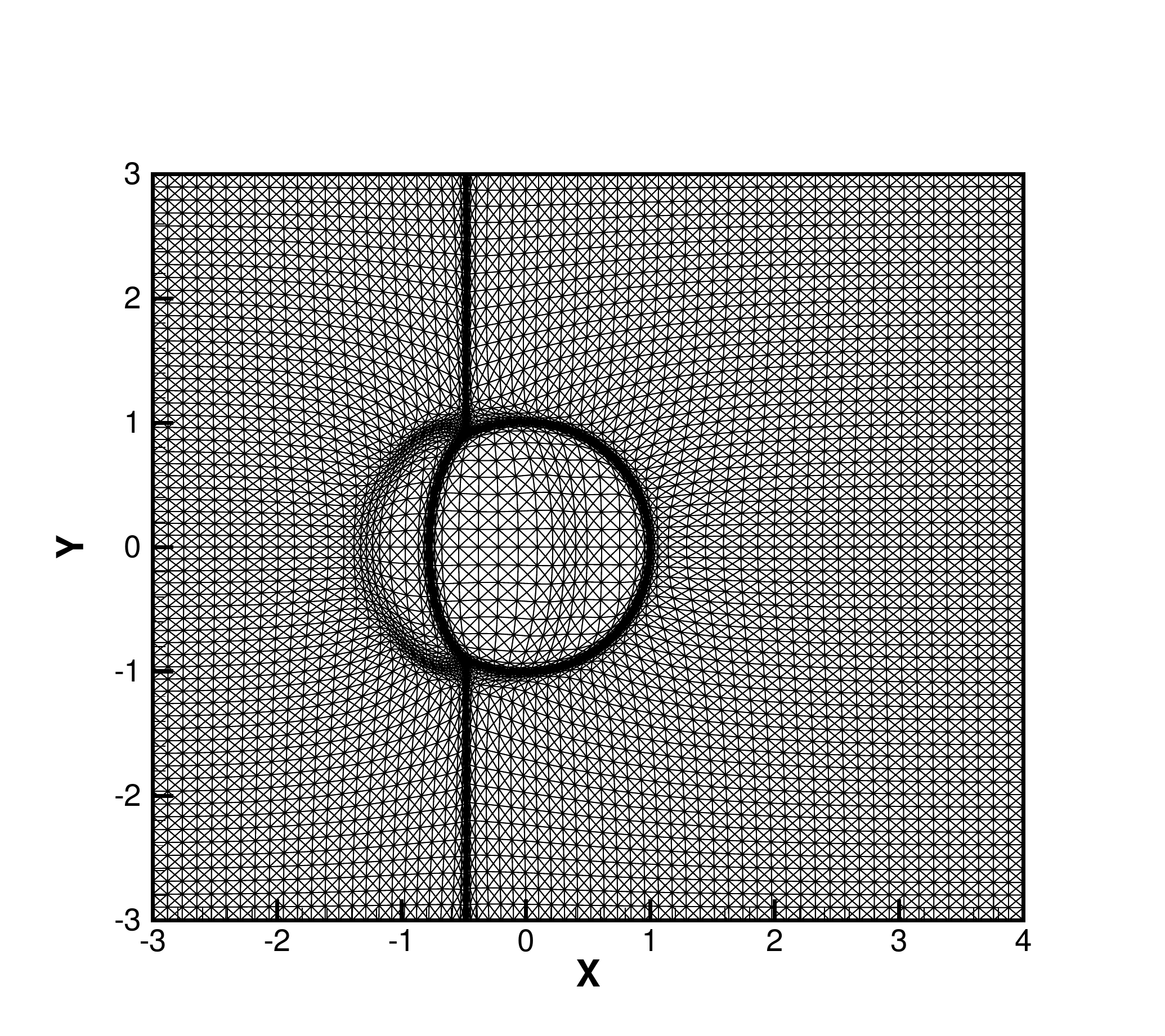}}
   }
   
 \mbox{%\subfigure[t=0]
 {\includegraphics[width=0.35\textwidth]{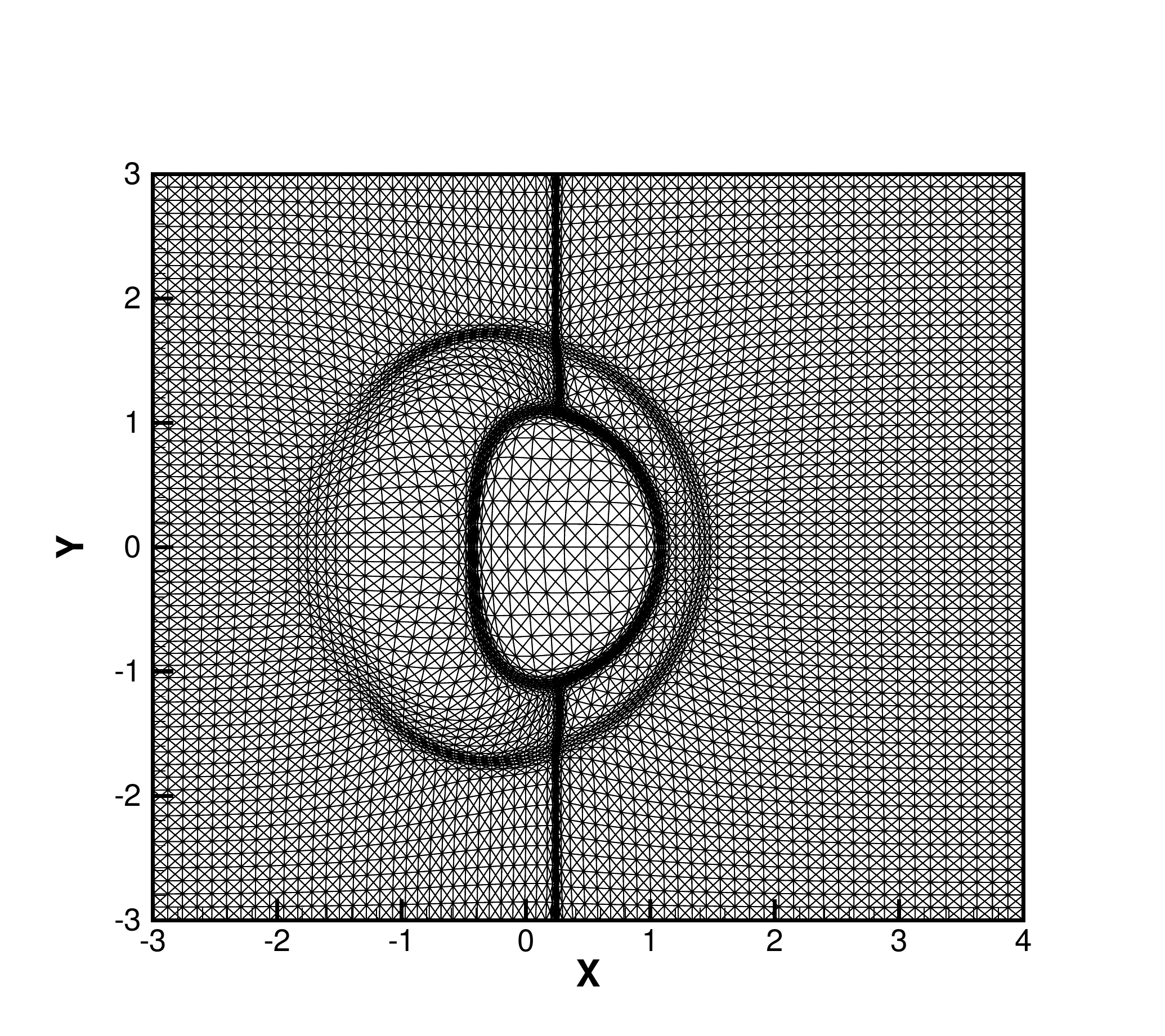}}\quad
  % \subfigure[close view of (a) near shock]
   {\includegraphics[width=0.35\textwidth]{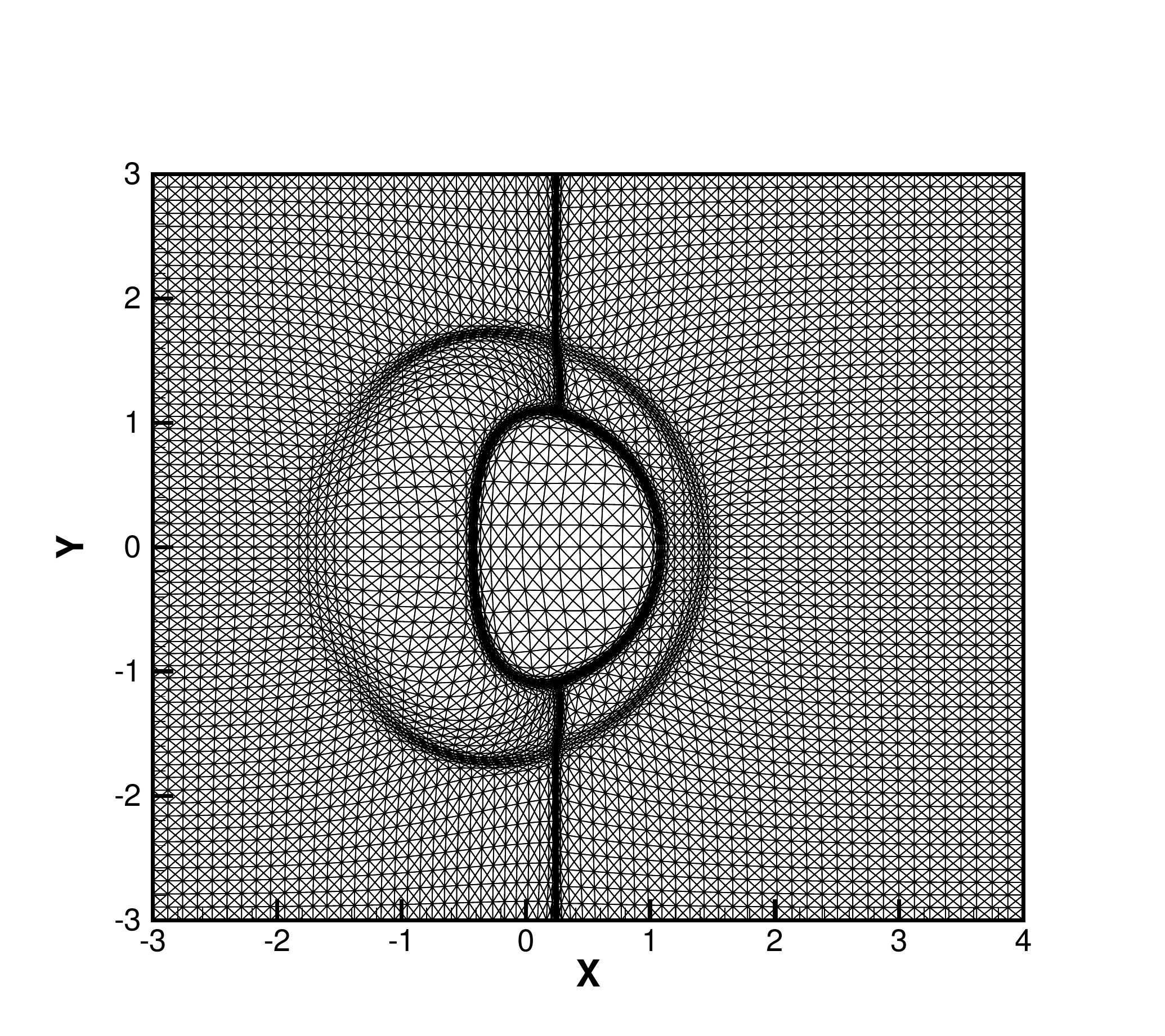}}
   }

 \mbox{%\subfigure[t=0]
 {\includegraphics[width=0.35\textwidth]{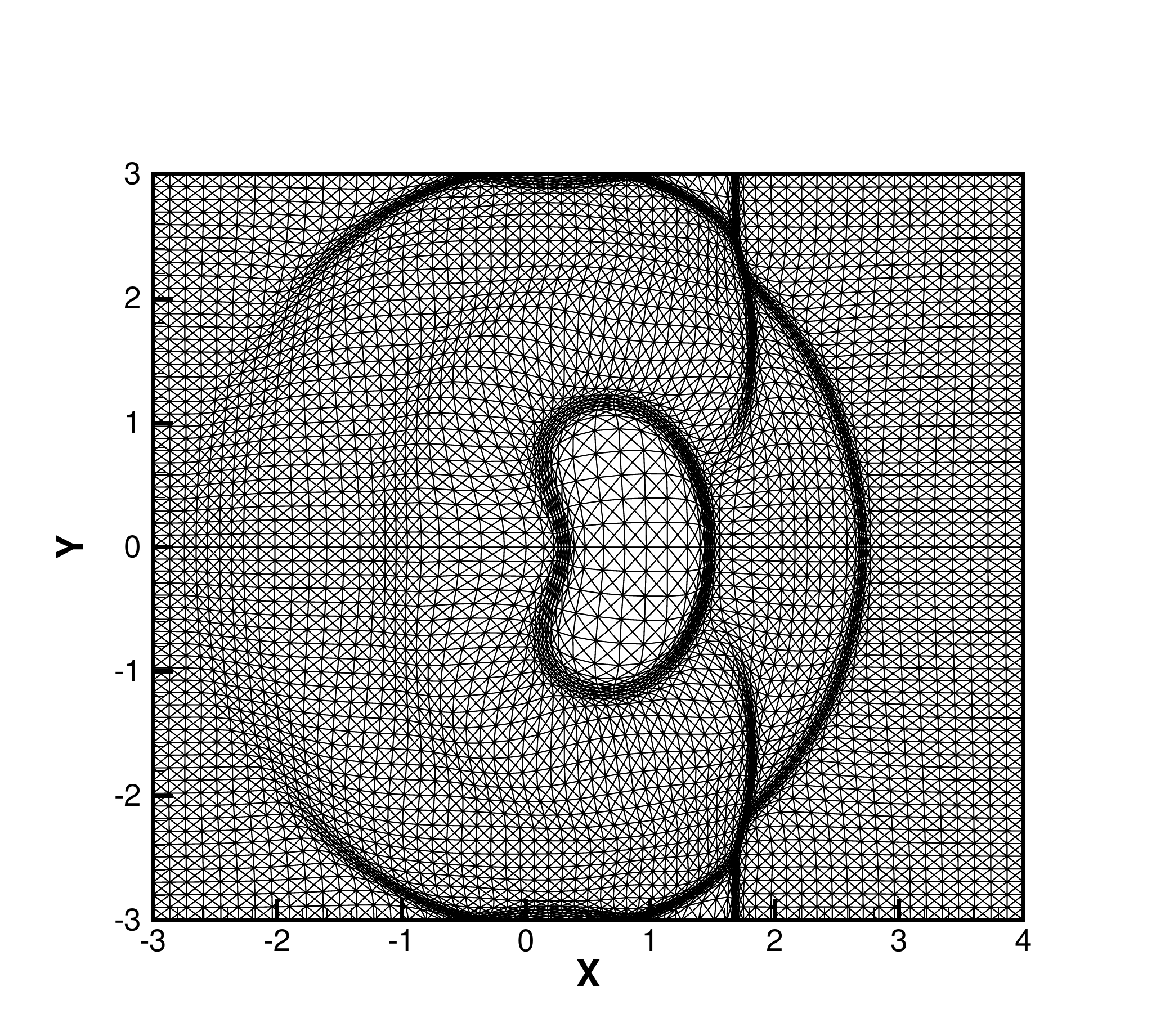}}\quad
  % \subfigure[close view of (a) near shock]
   {\includegraphics[width=0.35\textwidth]{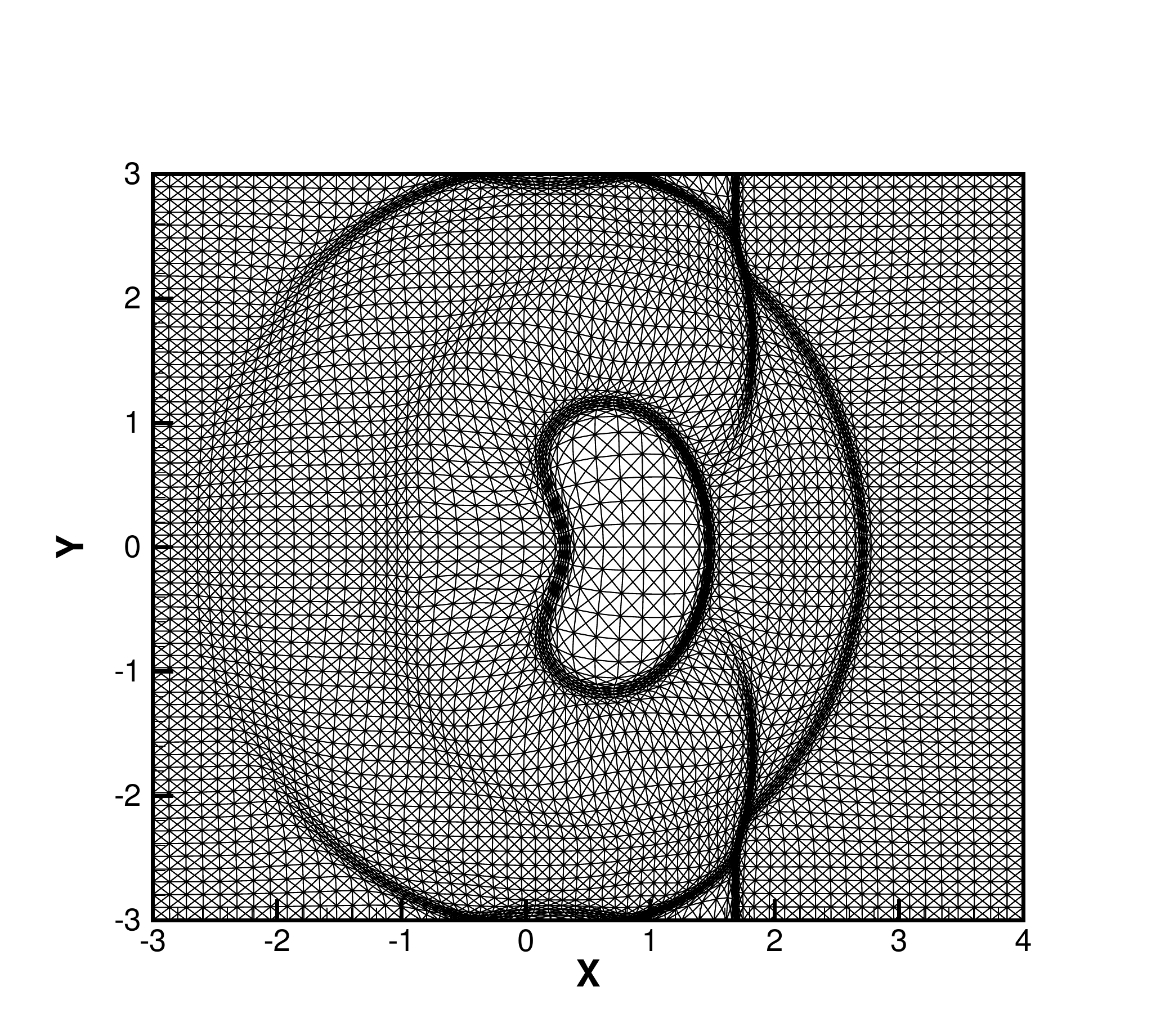}}
   }

 \mbox{%\subfigure[t=0]
 {\includegraphics[width=0.35\textwidth]{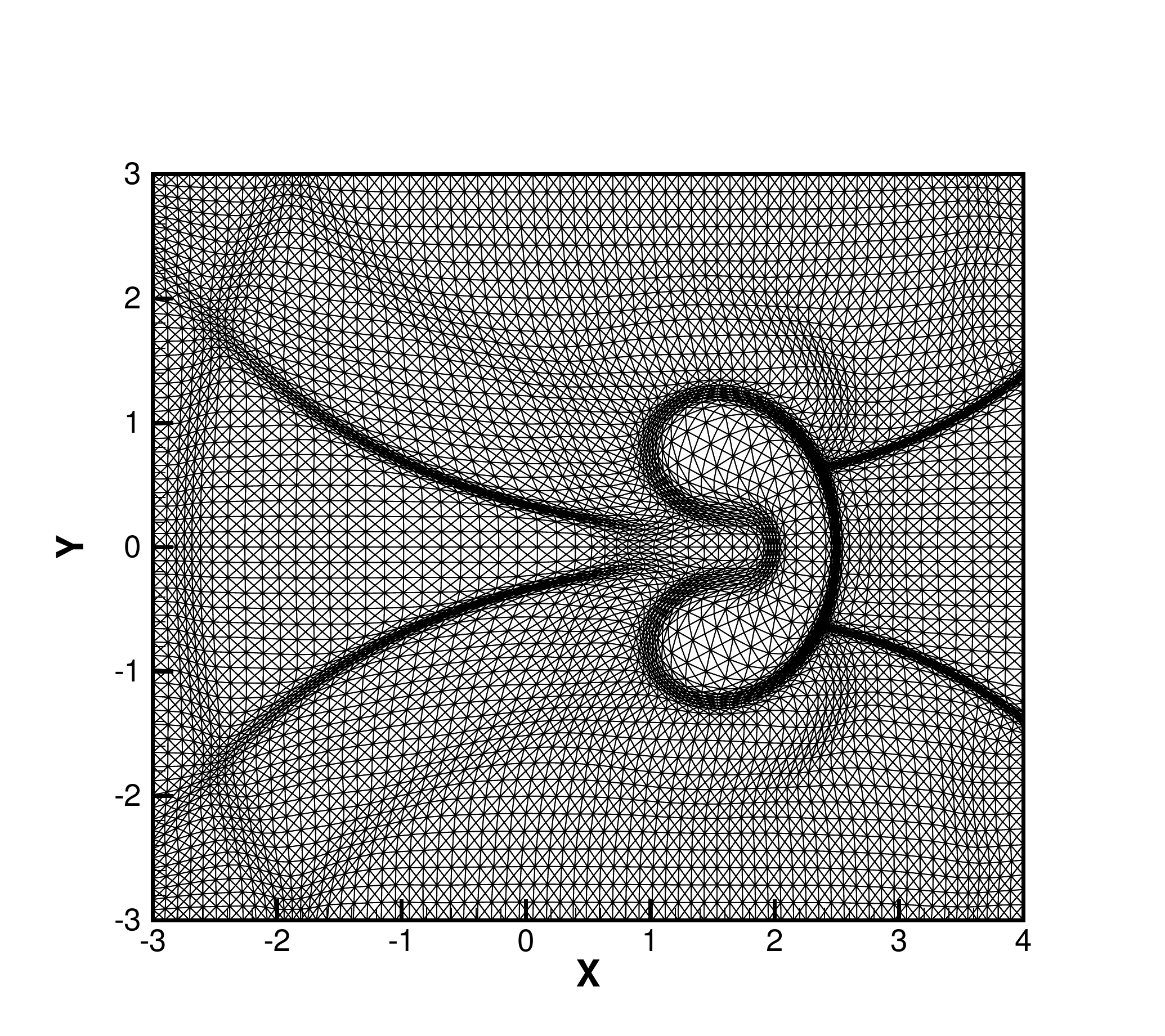}}\quad
  % \subfigure[close view of (a) near shock]
   {\includegraphics[width=0.35\textwidth]{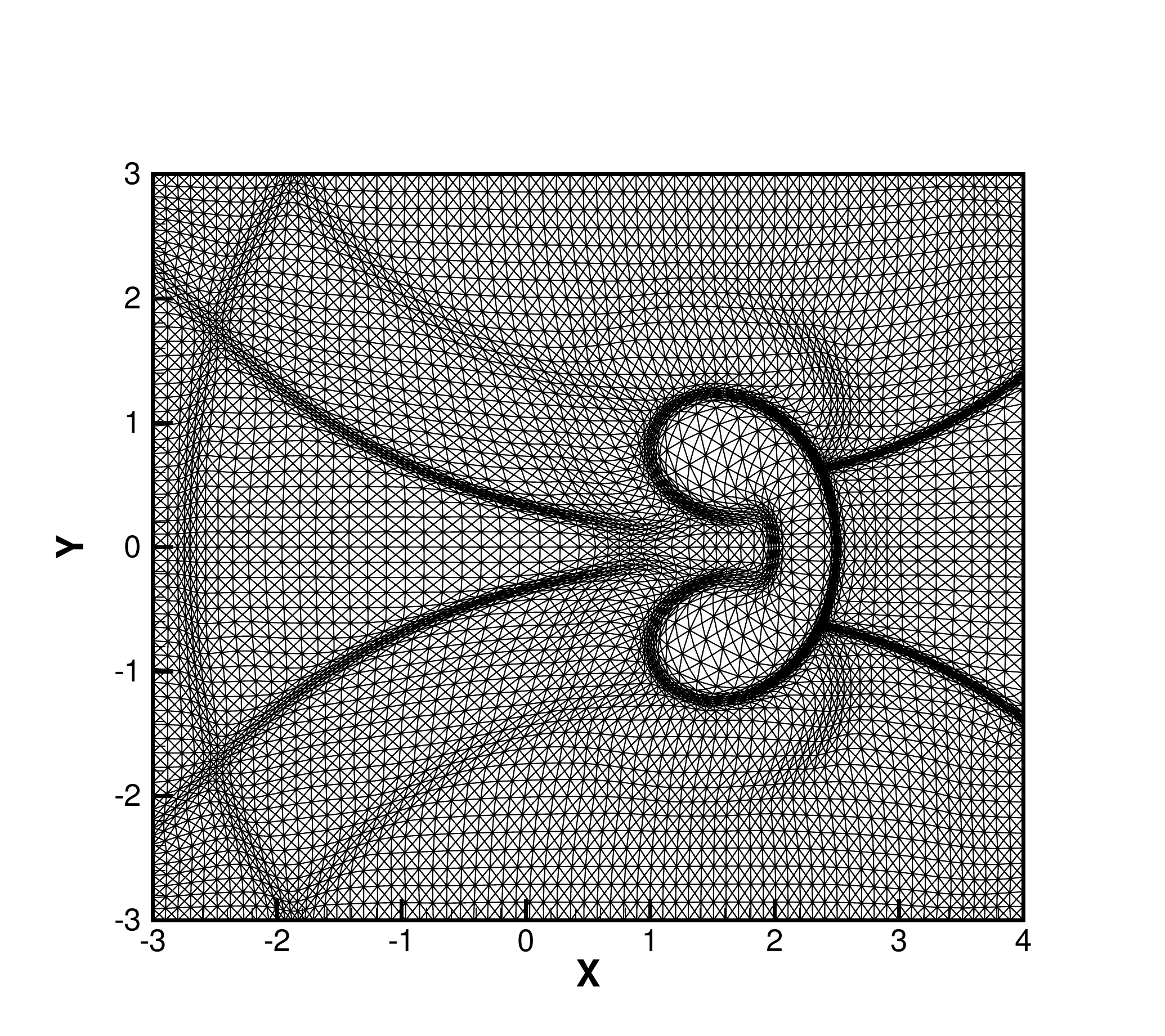}}
   }

   \caption{Example~\ref{examgasgas} The adaptive moving mesh of $N=70\times 60\times 4$ at time $t=0.5$, 1, 2, and 4.
   Left column: $P^1$-DG; Right column: $P^2$-DG.}
   \label{figgasgastraj}
   \end{center}
   \end{figure}

}
\end{exam}

\begin{exam}{\em
\label{exambw}
To demonstrate the performance of the current DG-ALE method with high pressure ratio in two dimensions,
we consider a model underwater explosion problem \cite{luo2020, shyue2006}. In this test, the computation domain is taken as $(-2,2)\times (-1.5,1)$.
Initially, the horizontal air-water interface is located at the $y=0$ and the center of a circular gas bubble with the radius 0.12 in the water is located
at $(0,-0.3)$. Above the air-water interface, the fluid is a perfect gas at the standard atmospheric condition and below the air-water interface
in region outside the gas bubble the fluid is water. Thus the initial condition is
\begin{equation*}
(\rho,U,V,P,\gamma,B)=
\begin{cases}
(1.225,0,0,101325,1.4,0),  \; & y > 0\\
(1250,0,0,10^9,1.4,0), \;  &x^2+(y+0.3)^2\leqslant 0.12^2  \text{ and } y \le 0\\
(1000,0,0,101325,4.4,6\times 10^8), \; &\text{otherwise}.
\end{cases}
\end{equation*}
A reflecting boundary condition is employed on the bottom of the domain and non-reflecting boundary conditions are used for the remaining boundary.
The parameters $\beta_i$'s in \eqref{ent} are all taken as 1 in this example.

At the beginning of the process, both the gas and water are in a stationary position. Due to the pressure difference between the fluids, the bubble begins to break and this results in a circularly outward-going shock wave in water, an inward-going rarefaction wave in gas, and an interface lying in between
those waves that separates the gas and the water. Soon after, the shock wave is diffracted through the nearby air-water surface, causing
the subsequent deformation of the interface topology from a circle to an oval-like shape.
 
The contours of the density and the adaptive mesh of  $N = 120\times 75\times 4$ are plotted at $t=0.2$, 0.4, 0.8, and 1.2
in Figs.~\ref{figbwden} and \ref{figbwmesh}. Once again, the mesh points are concentrated correctly near the material interface and shocks
and $P^2$-DG appears to give higher resolution than $P^1$-DG.
 
\begin{figure}[hbtp]
 \begin{center}
 \mbox{%\subfigure[t=0]
 {\includegraphics[width=0.45\textwidth]{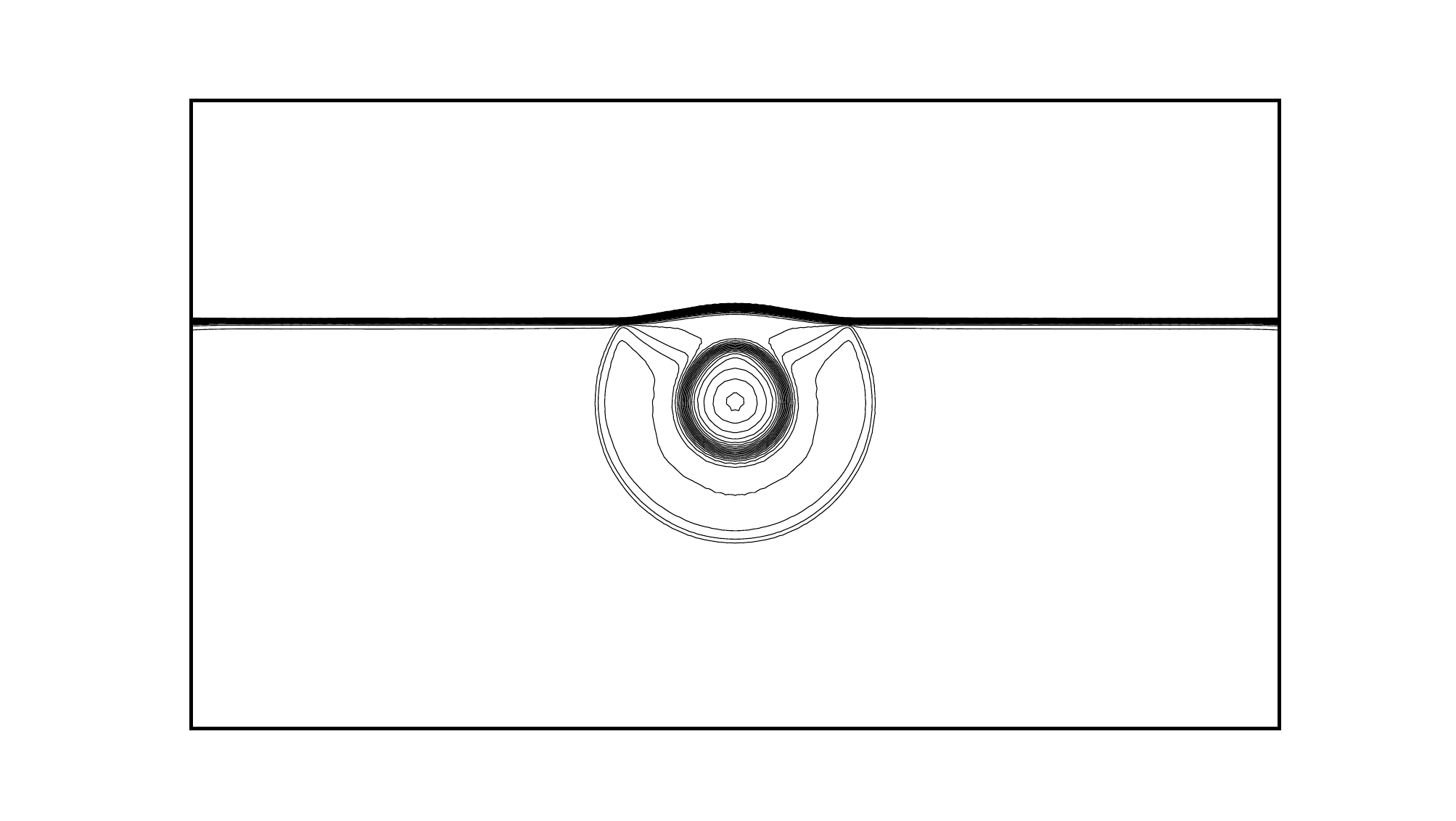}}\quad
  % \subfigure[close view of (a) near shock]
   {\includegraphics[width=0.45\textwidth]{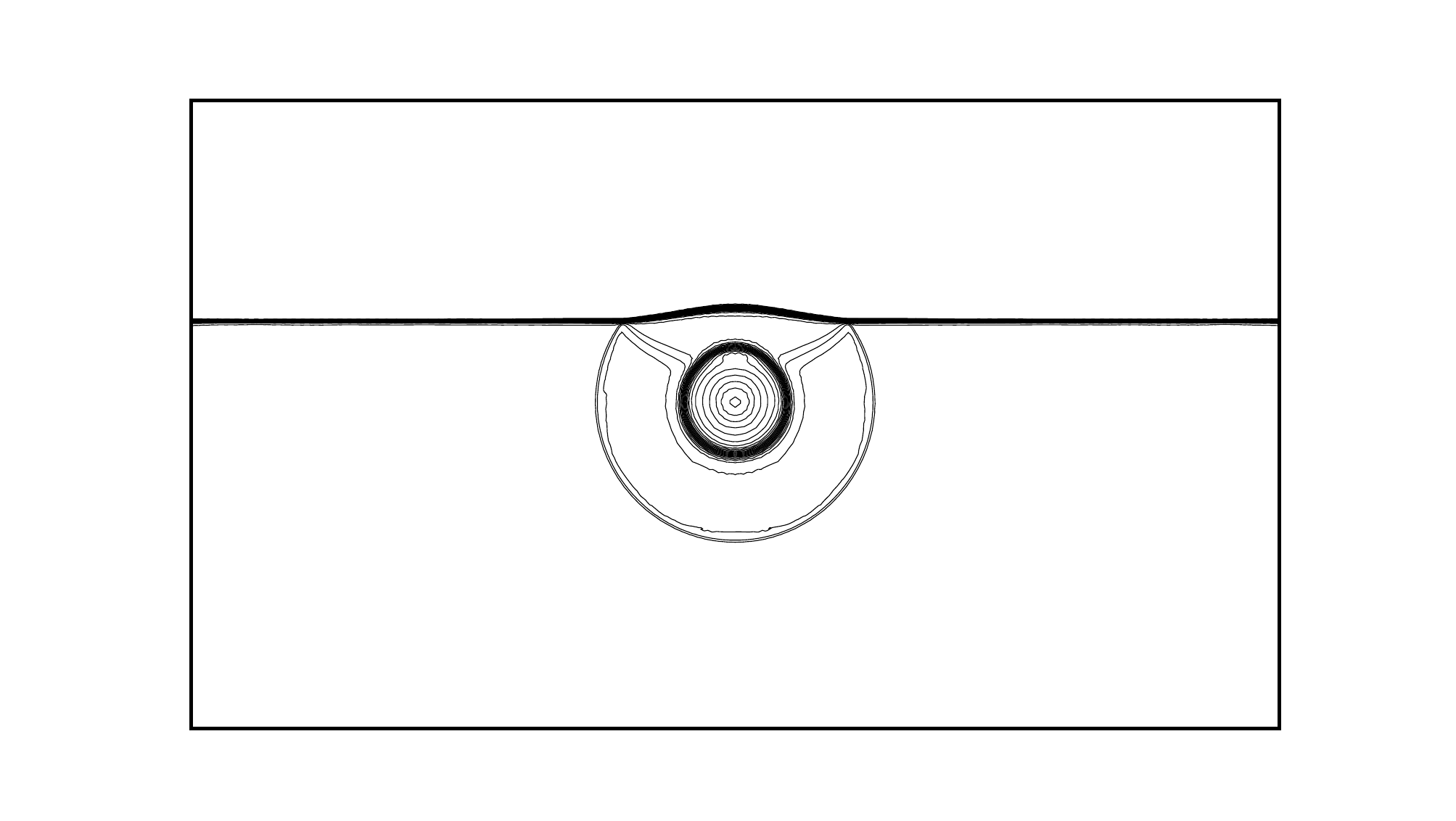}}
   }
   
 \mbox{%\subfigure[t=0]
 {\includegraphics[width=0.45\textwidth]{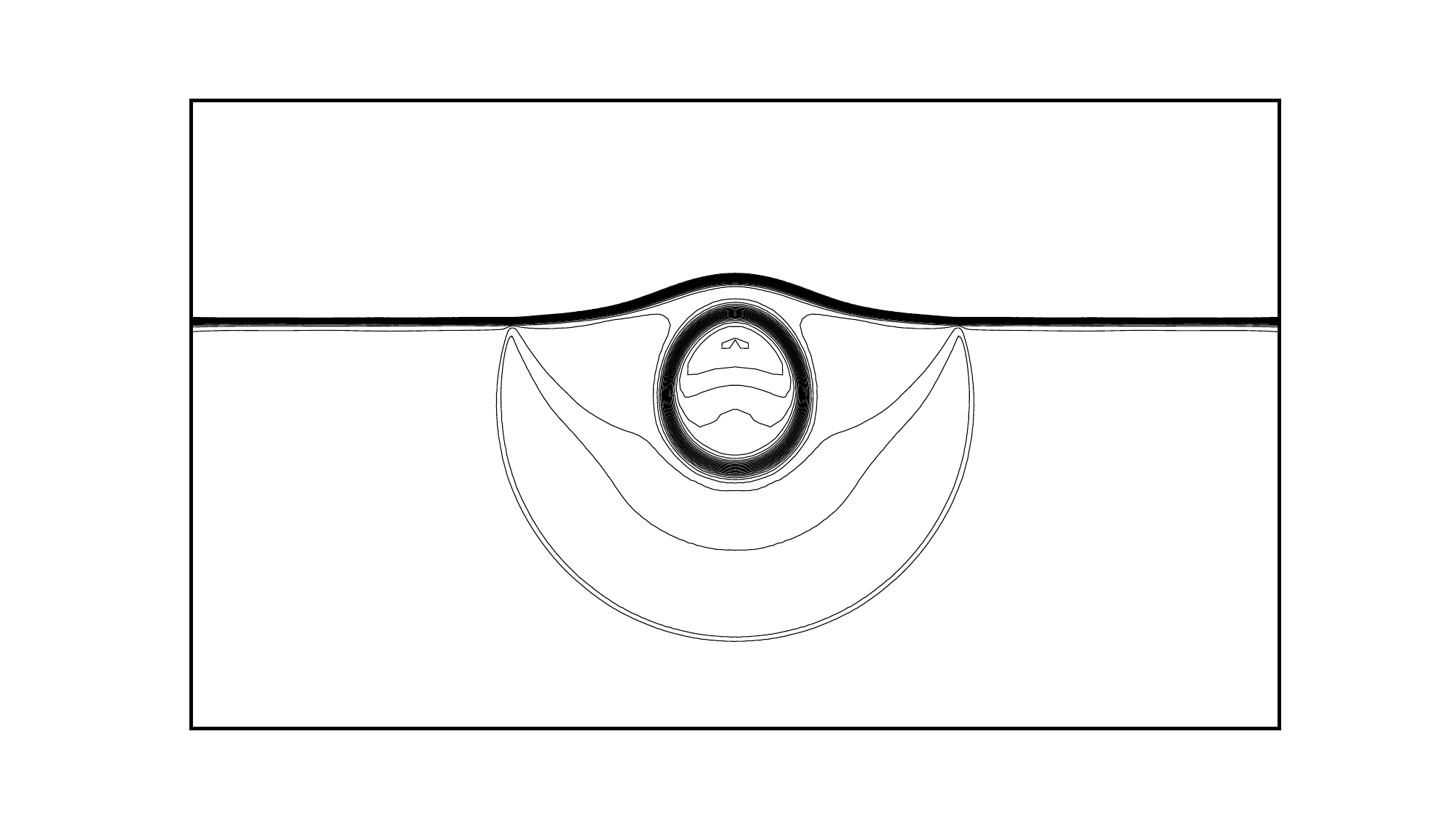}}\quad
  % \subfigure[close view of (a) near shock]
   {\includegraphics[width=0.45\textwidth]{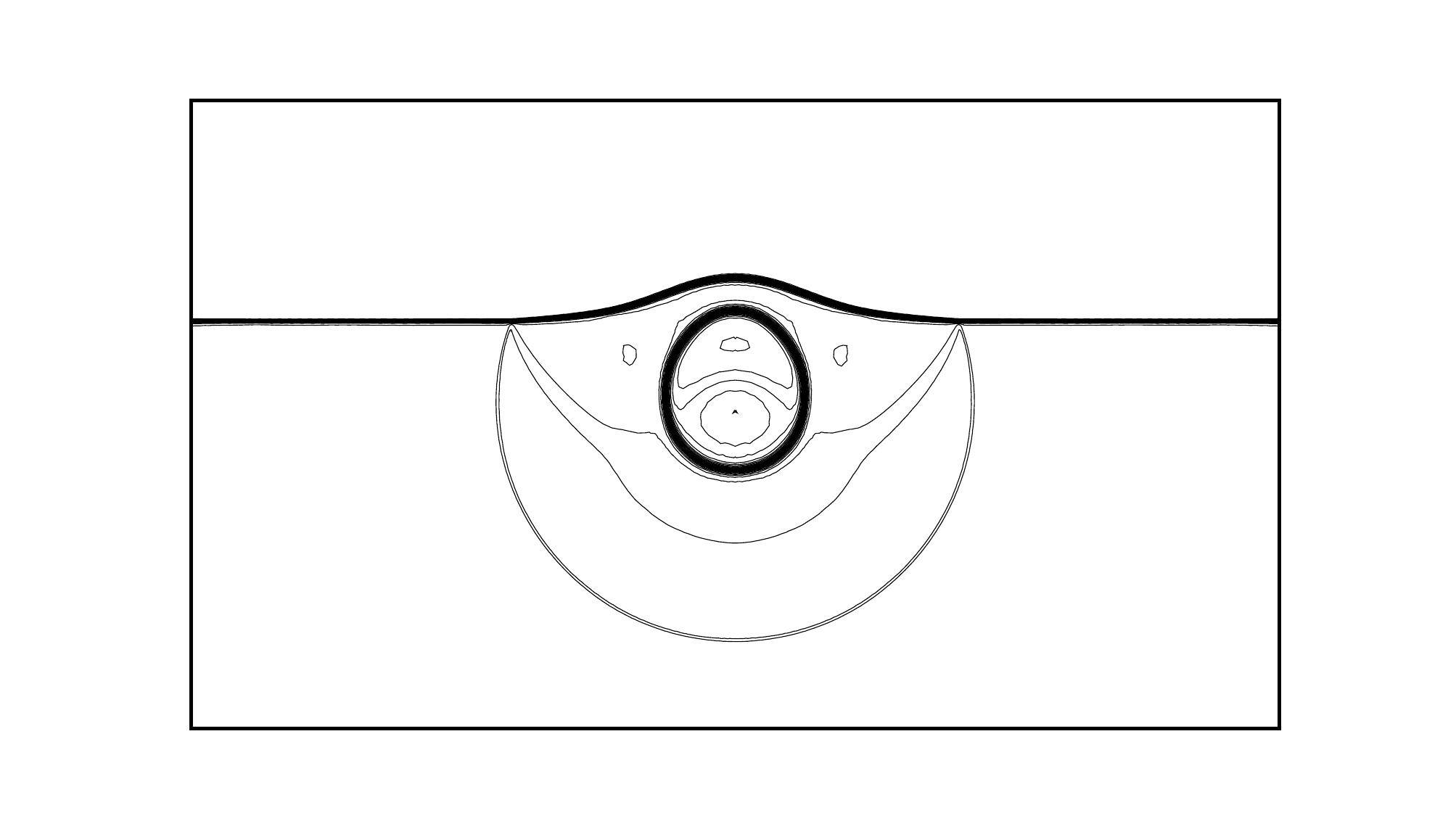}}
   }
 \mbox{%\subfigure[t=0]
 {\includegraphics[width=0.45\textwidth]{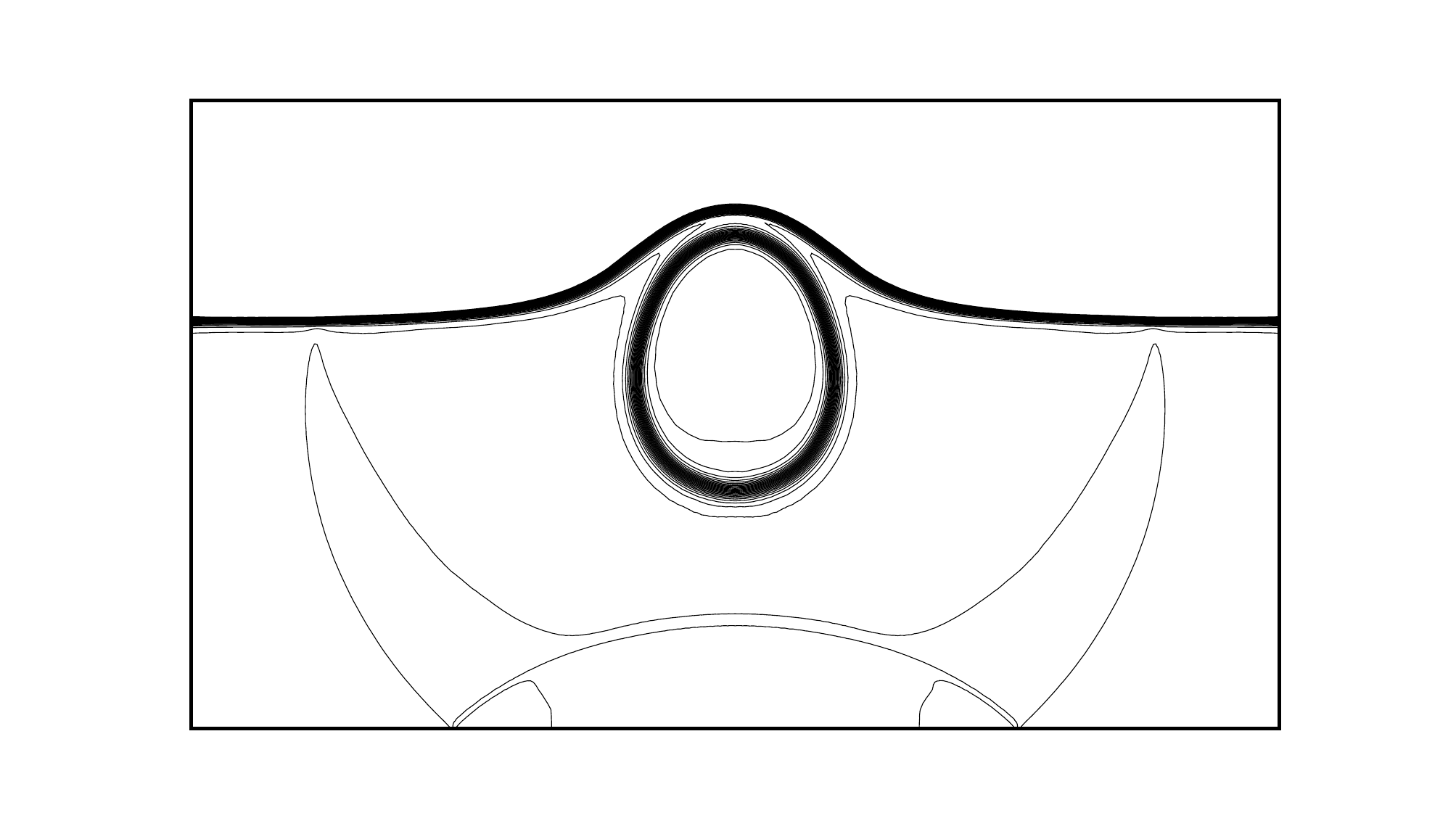}}\quad
  % \subfigure[close view of (a) near shock]
   {\includegraphics[width=0.45\textwidth]{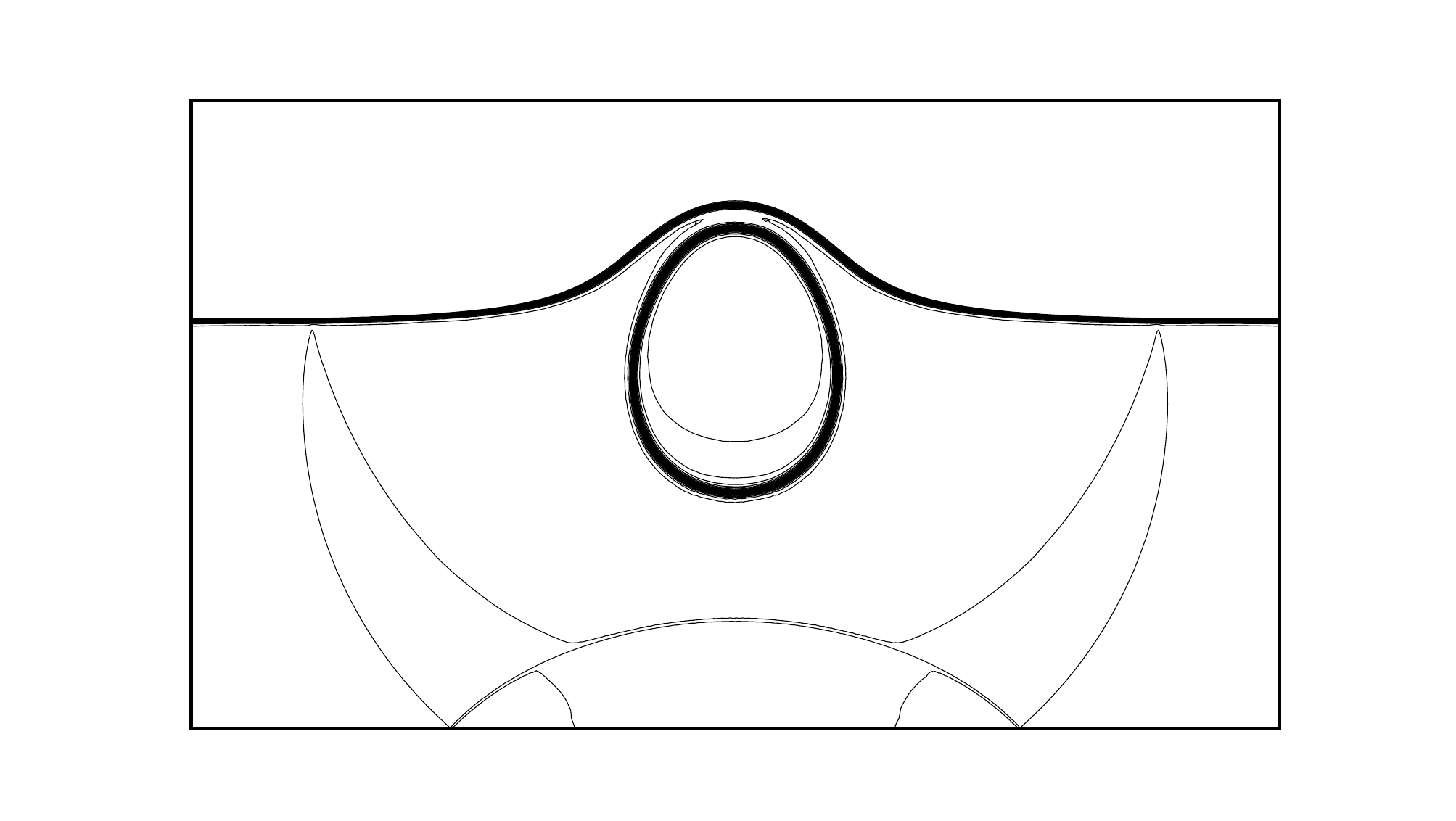}}
   }
   
 \mbox{%\subfigure[t=0]
 {\includegraphics[width=0.45\textwidth]{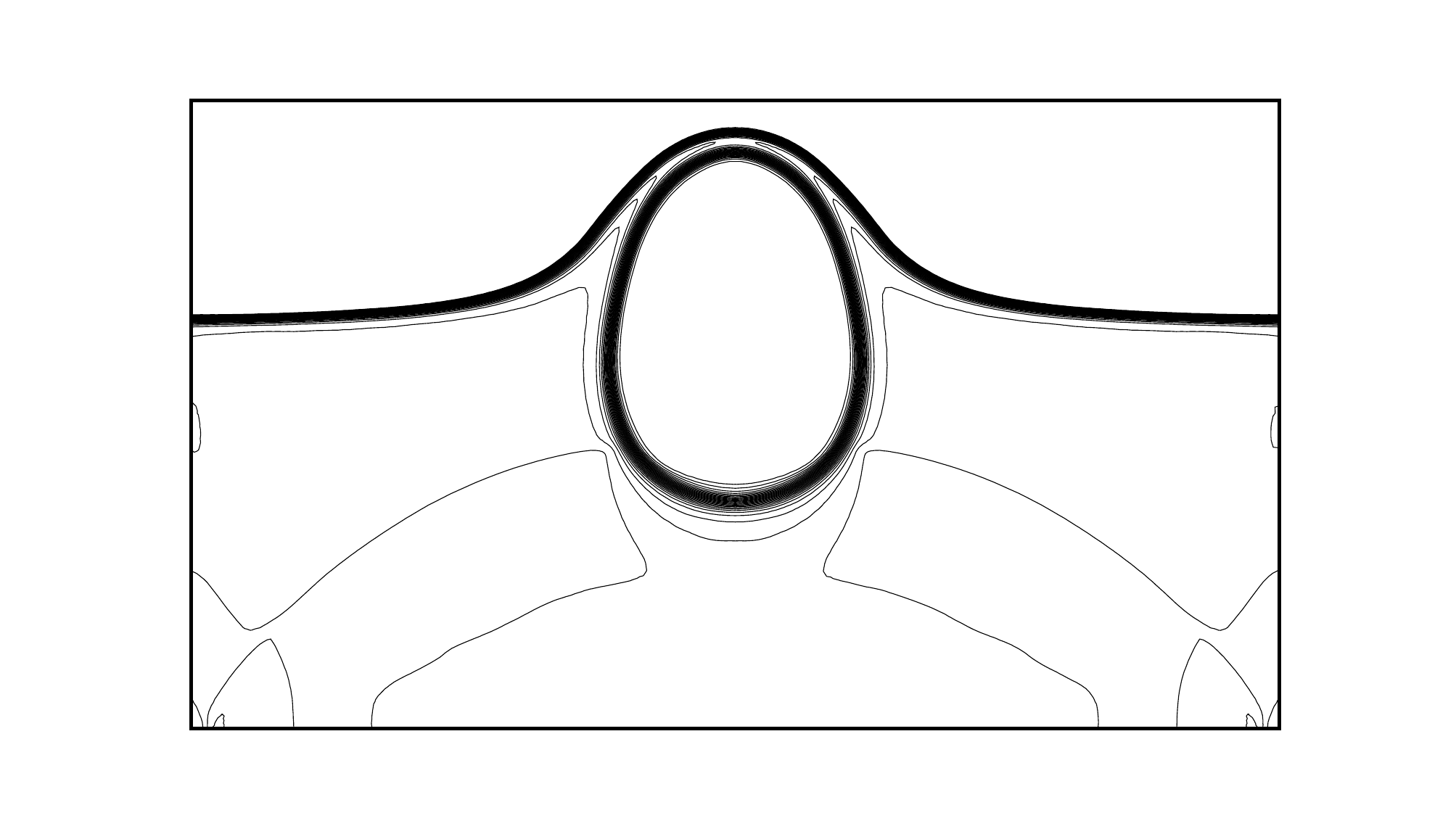}}\quad
  % \subfigure[close view of (a) near shock]
   {\includegraphics[width=0.45\textwidth]{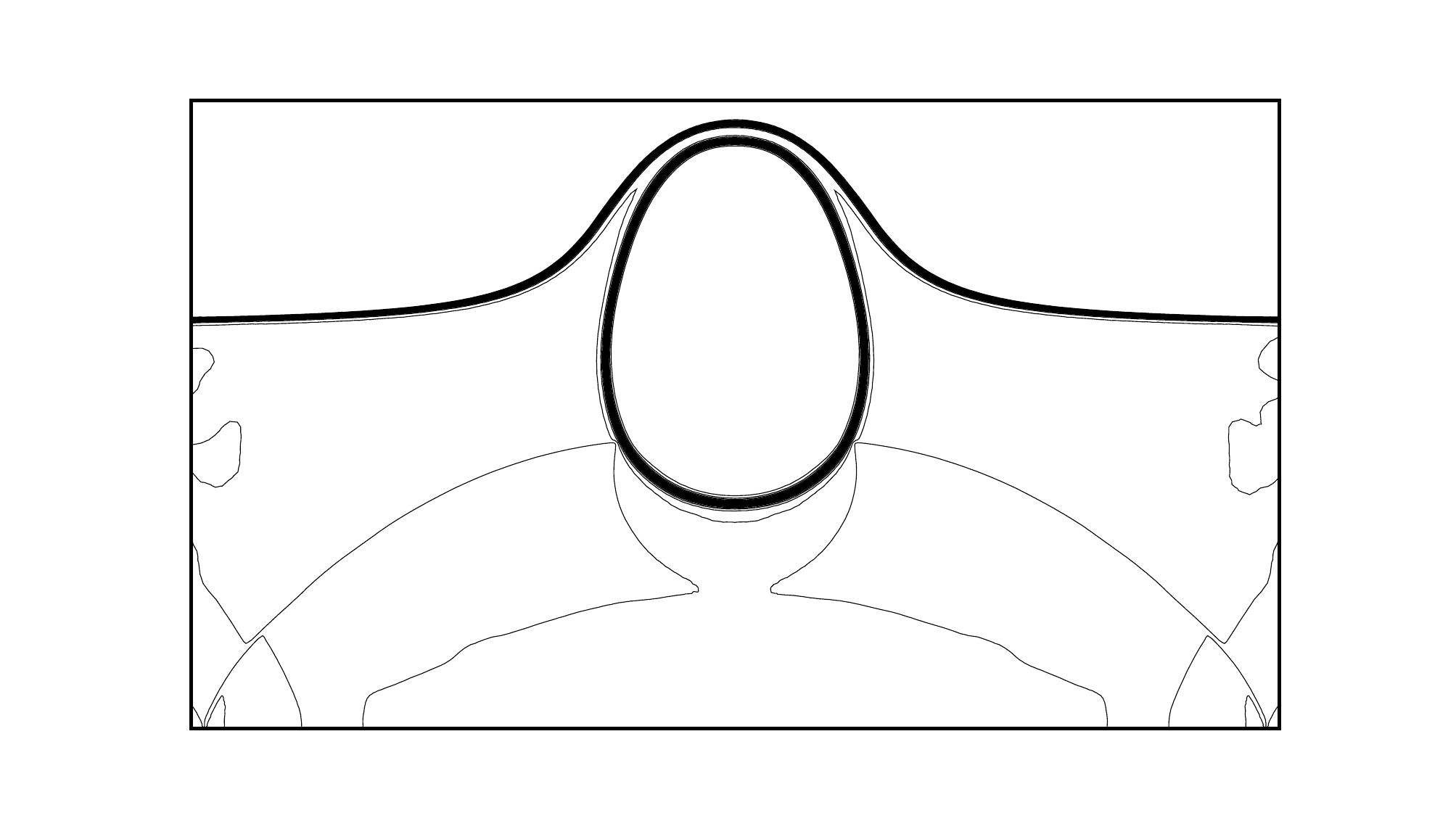}}
   }
   
   \caption{Example~\ref{exambw} Contours of the density at $t=0.2$, 0.4, 0.8, and 1.2 ms (from top to bottom).
   Left column: $P^1$-DG; Right column: $P^2$-DG. $N = 120\times 75\times 4$.}
   \label{figbwden}
   \end{center}
   \end{figure}

   %%P2
\begin{figure}[hbtp]
 \begin{center}
 \mbox{%\subfigure[t=0]
 {\includegraphics[width=0.45\textwidth]{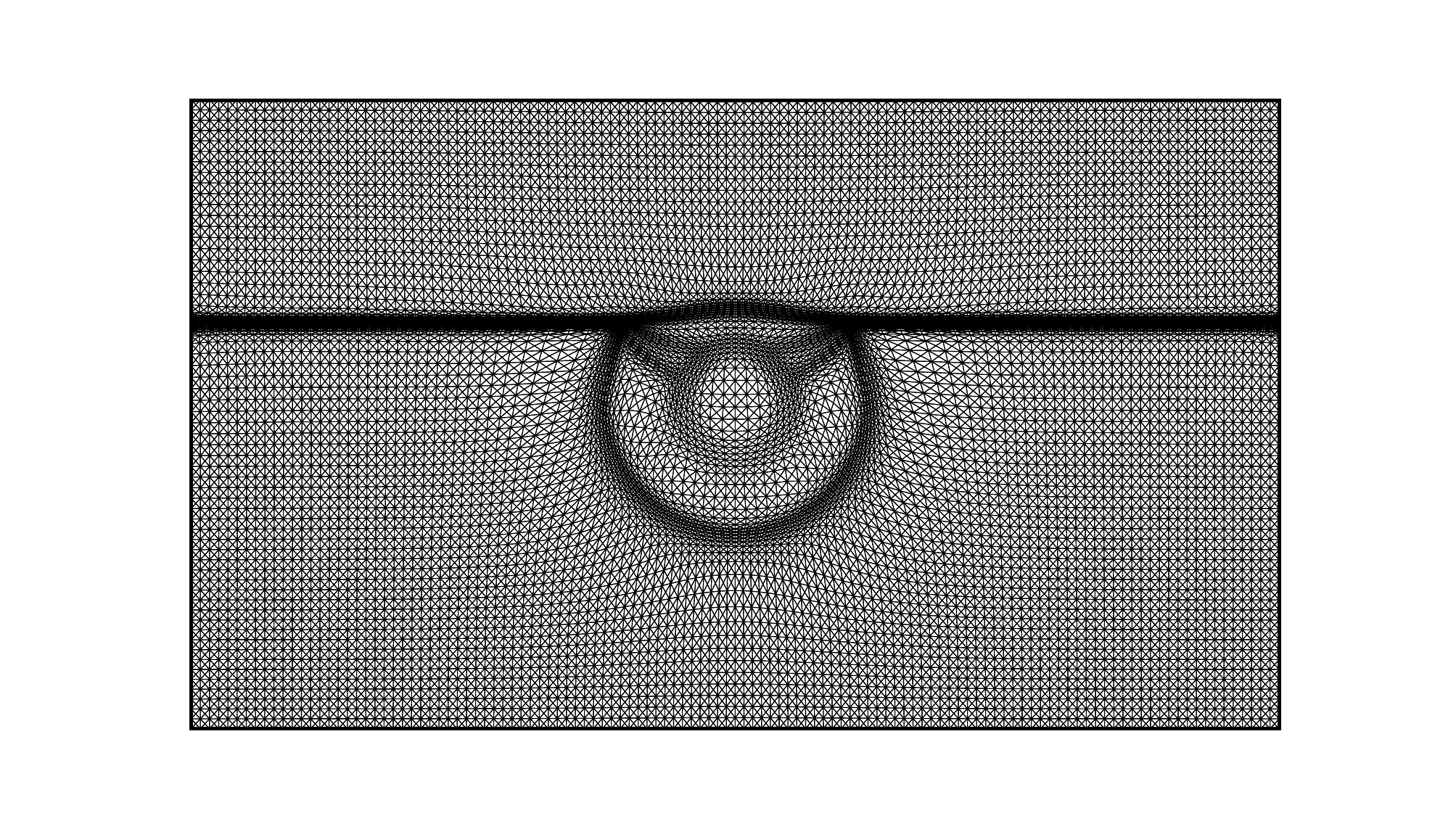}}\quad
  % \subfigure[close view of (a) near shock]
   {\includegraphics[width=0.45\textwidth]{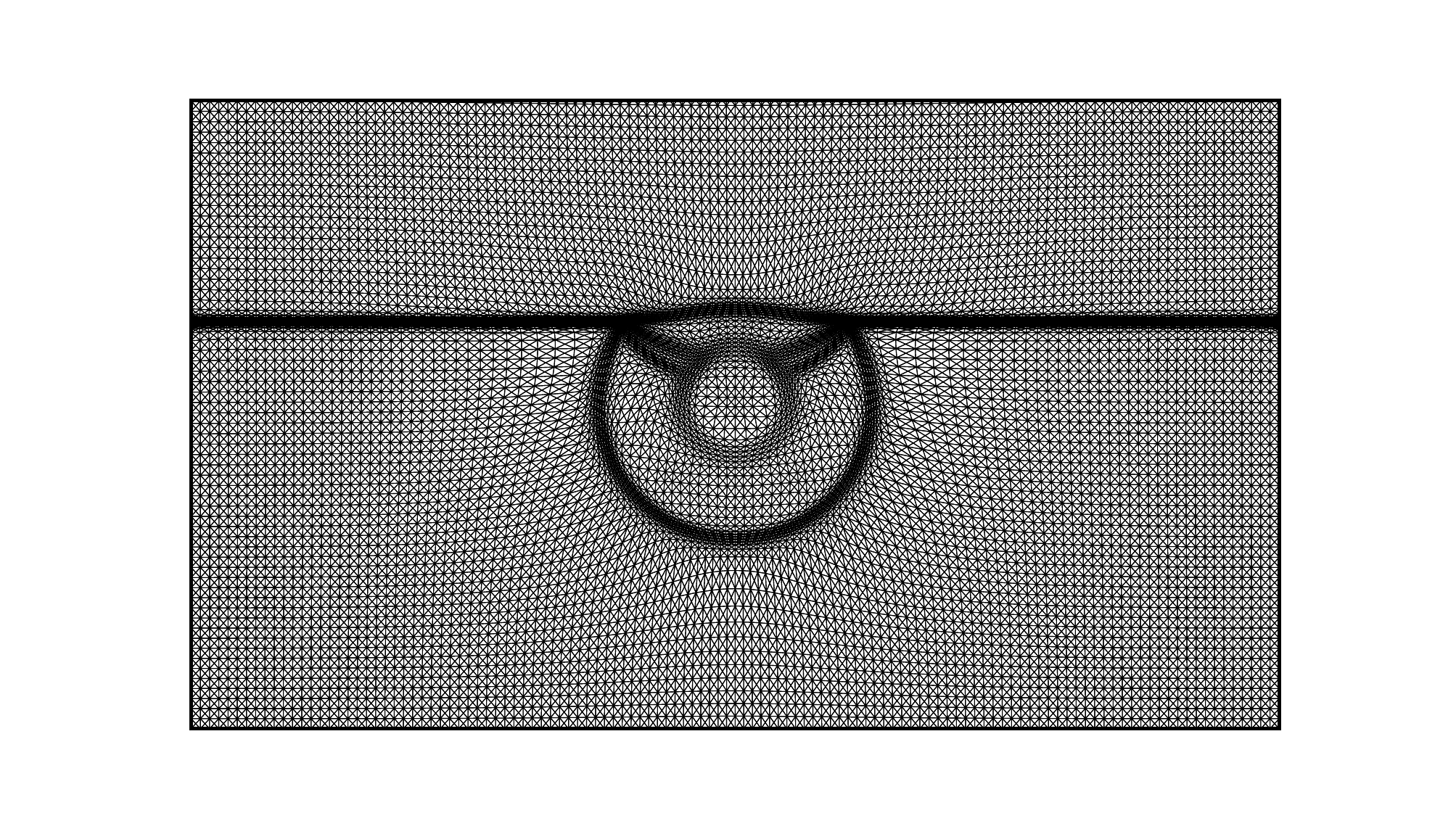}}
   }
   
 \mbox{%\subfigure[t=0]
 {\includegraphics[width=0.45\textwidth]{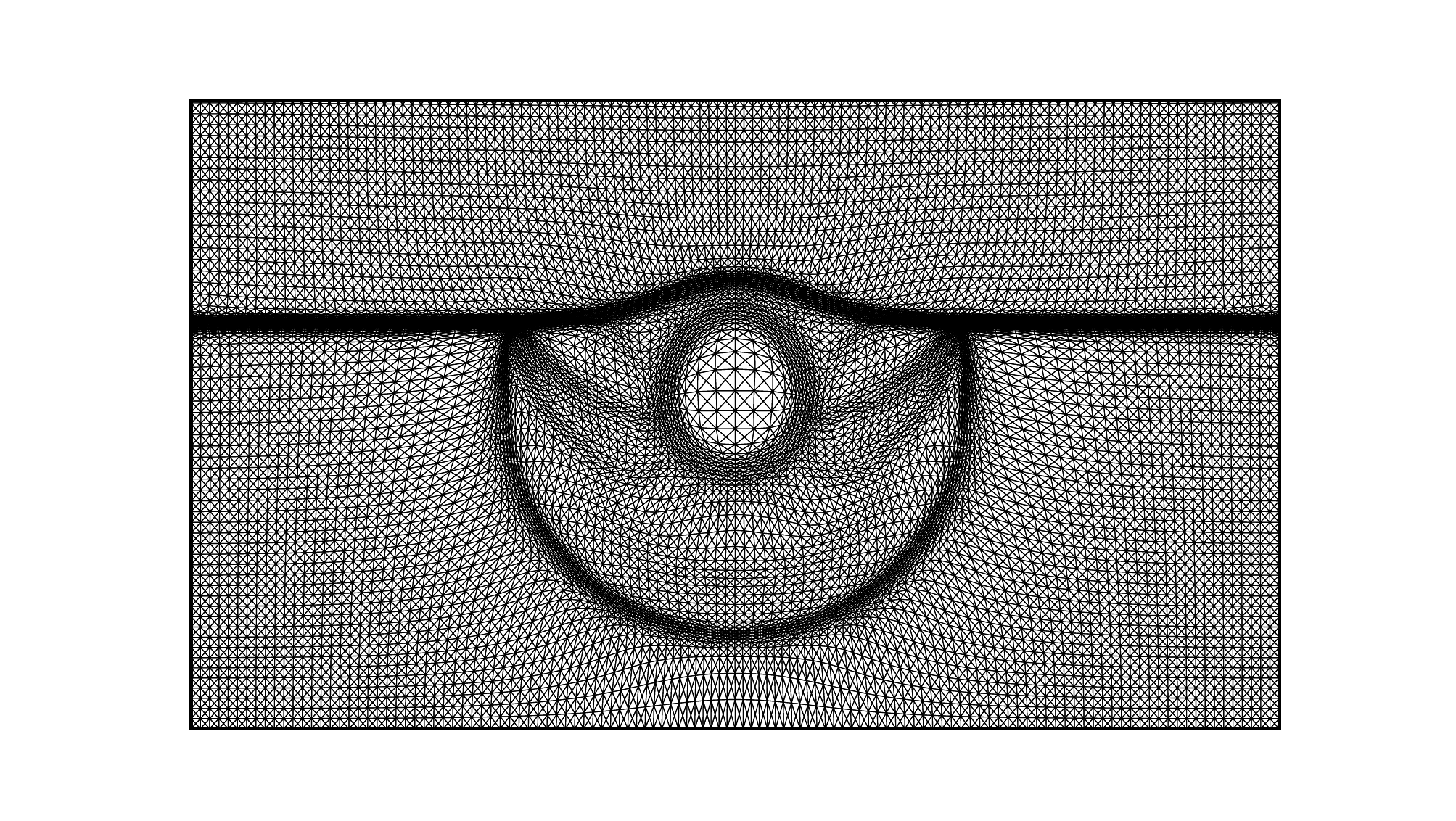}}\quad
  % \subfigure[close view of (a) near shock]
   {\includegraphics[width=0.45\textwidth]{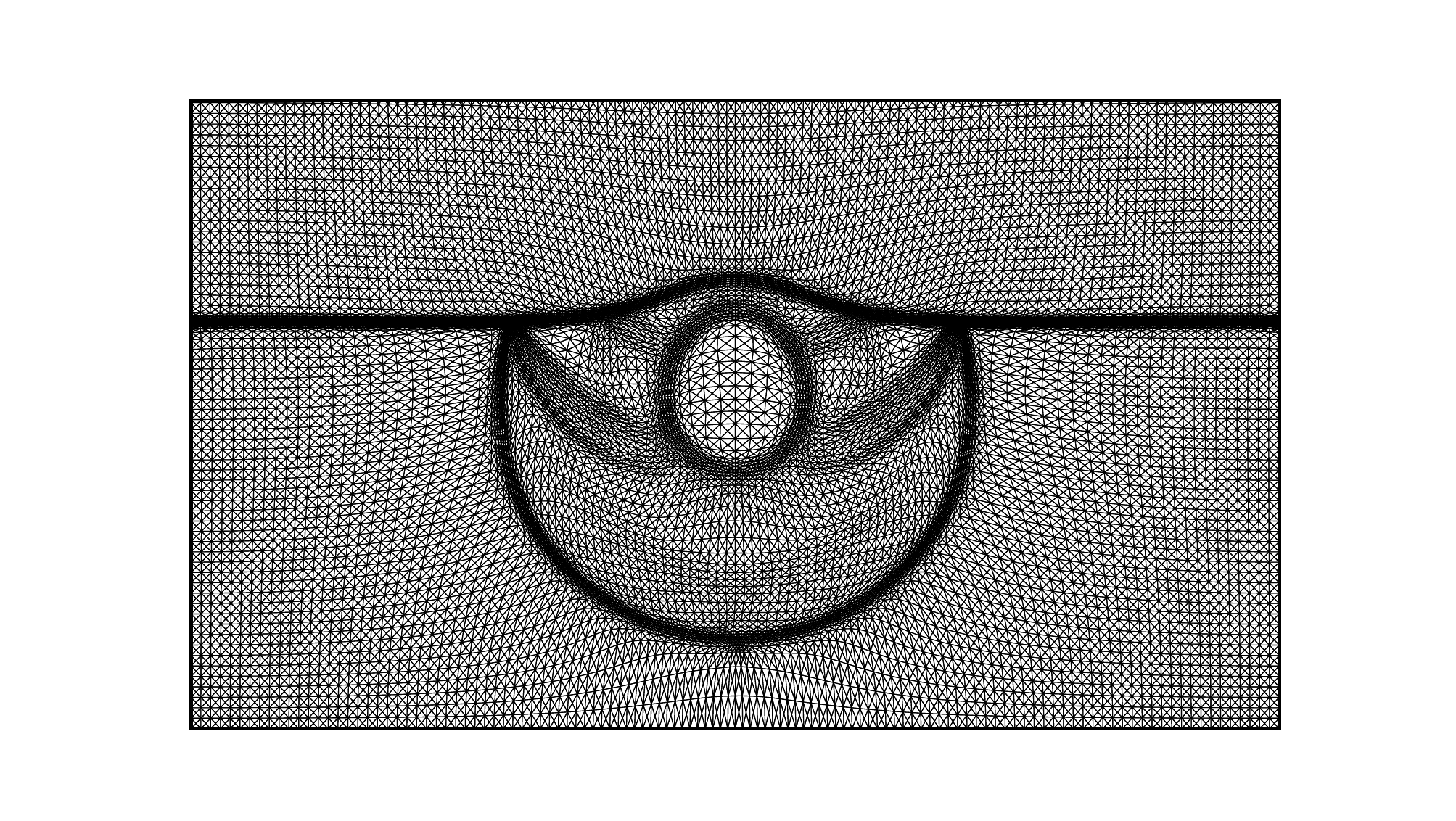}}
   }
 \mbox{%\subfigure[t=0]
 {\includegraphics[width=0.45\textwidth]{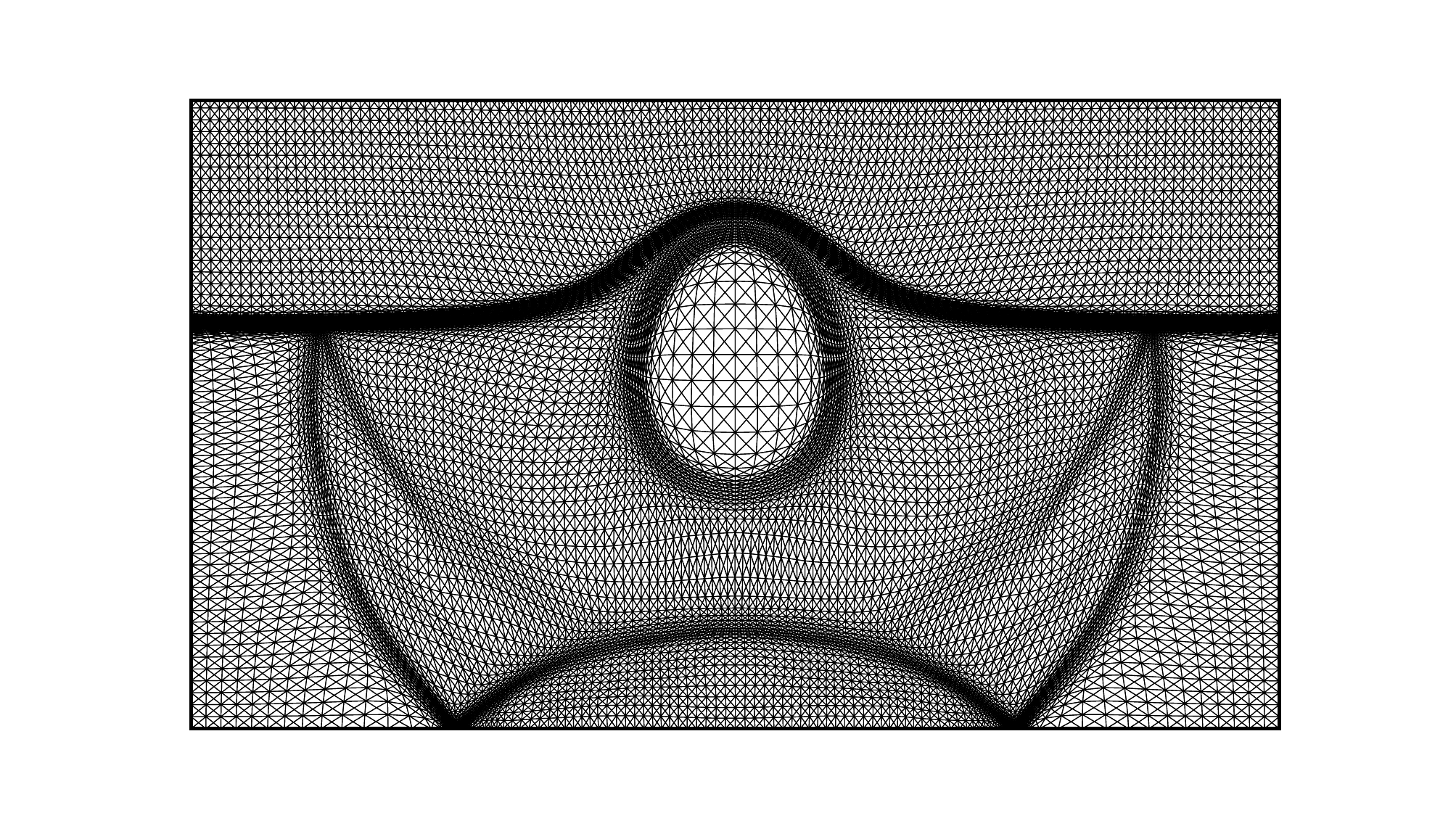}}\quad
  % \subfigure[close view of (a) near shock]
   {\includegraphics[width=0.45\textwidth]{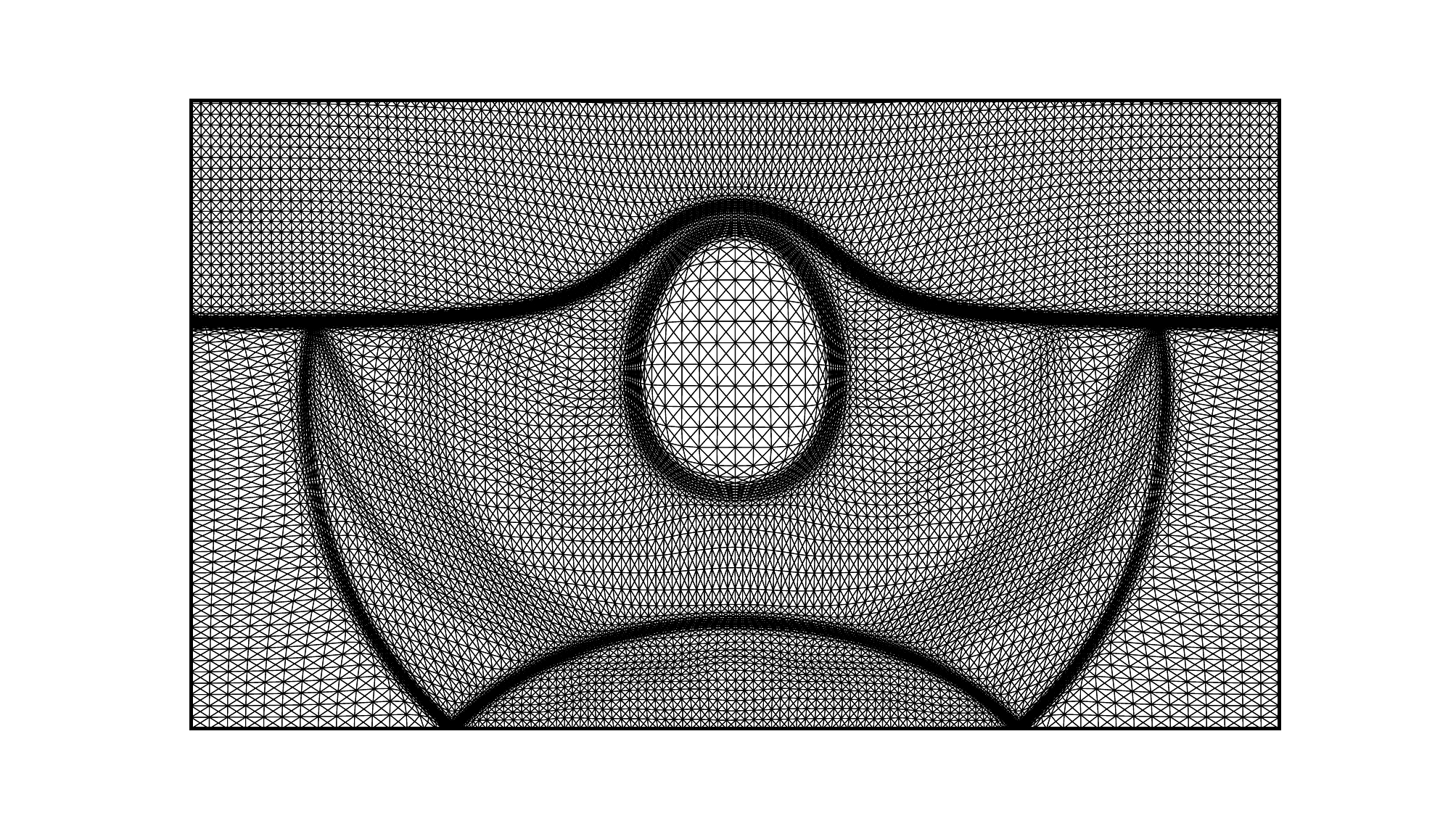}}
   }
   
 \mbox{%\subfigure[t=0]
 {\includegraphics[width=0.45\textwidth]{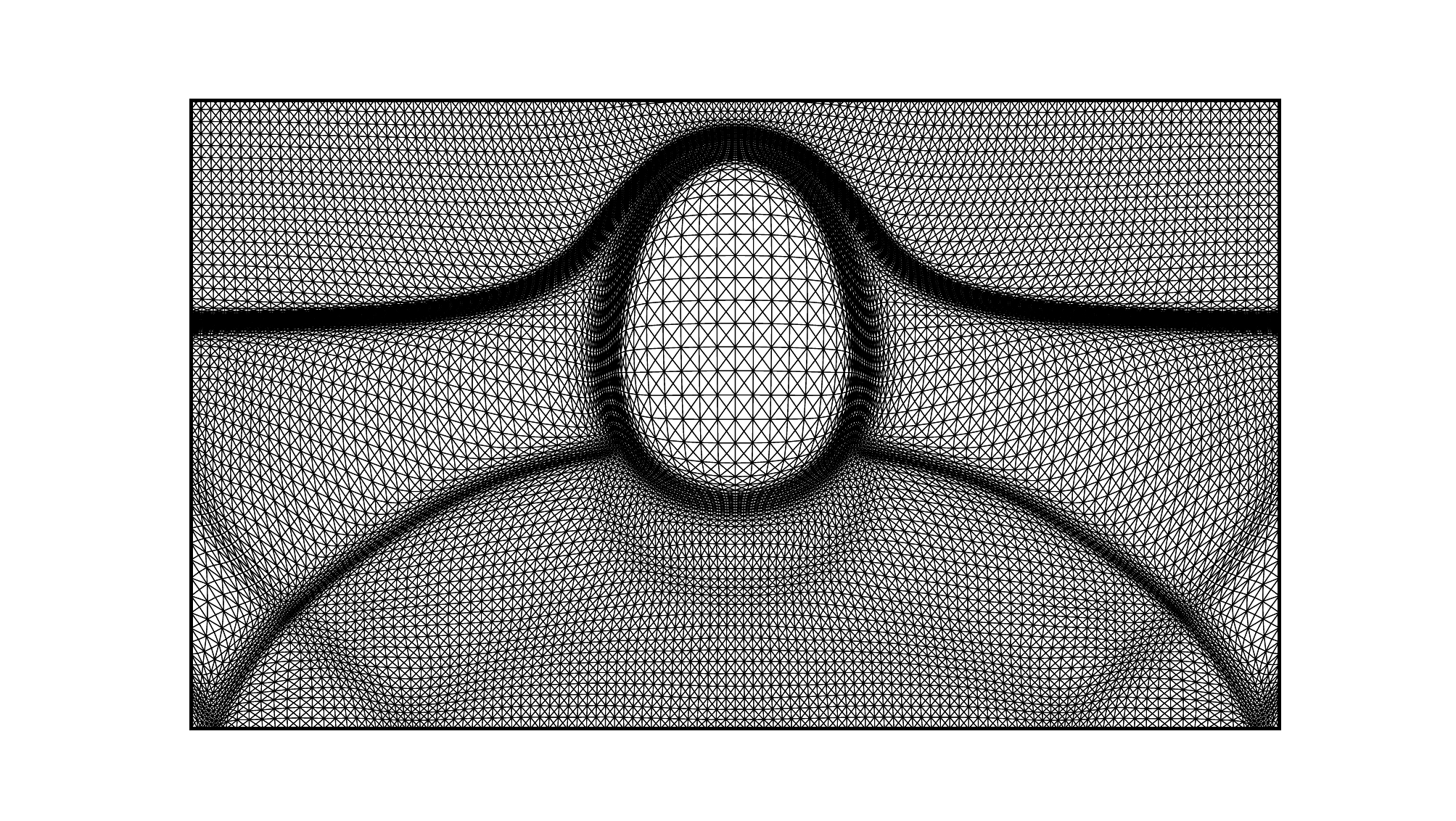}}\quad
  % \subfigure[close view of (a) near shock]
   {\includegraphics[width=0.45\textwidth]{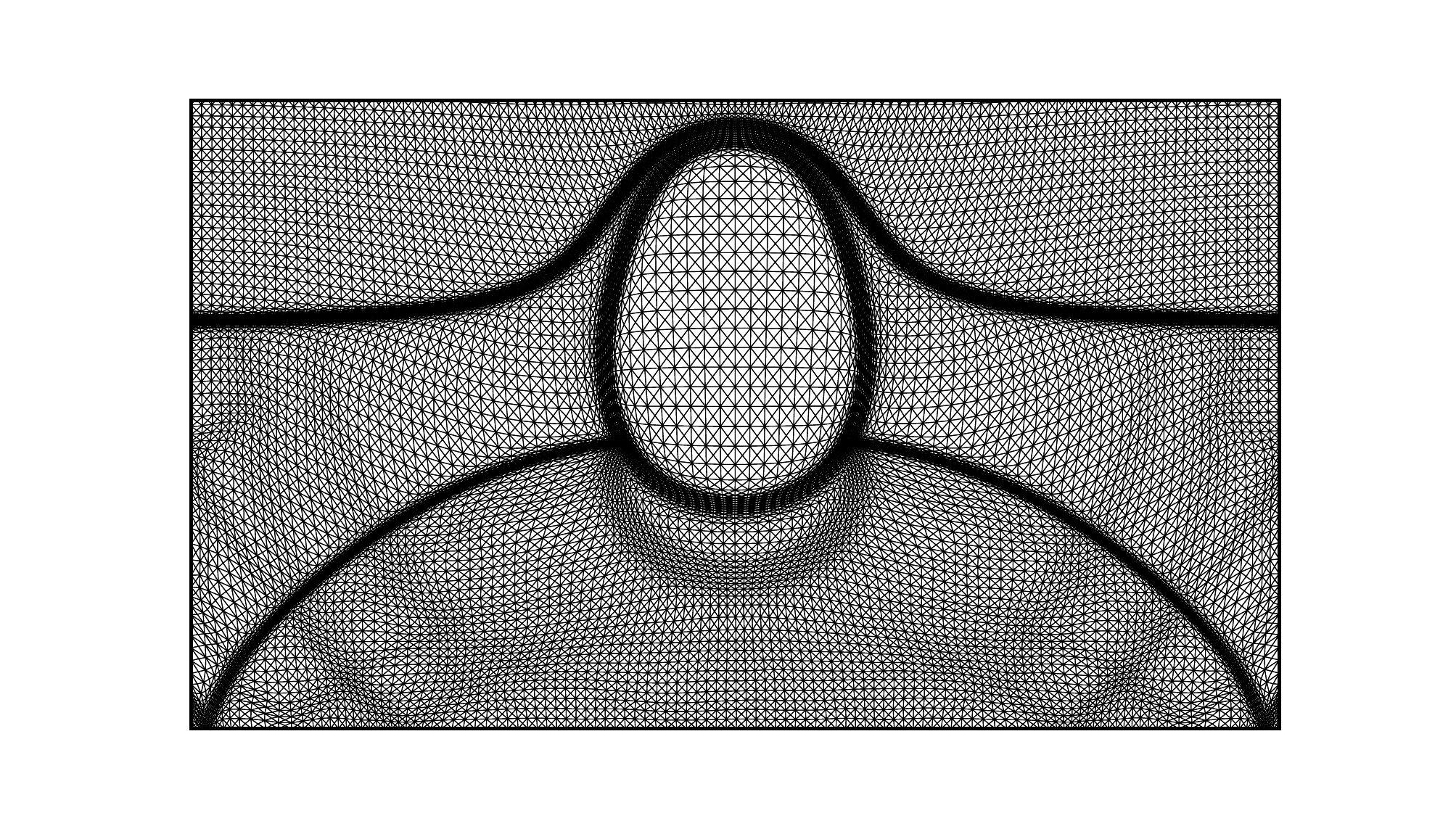}}
   }

   \caption{Example~\ref{exambw} The adaptive moving mesh of $N = 120\times 75\times 4$ at $t=0.2$, 0.4, 0.8, and 1.2
   (from  top to bottom). Left column: $P^1$-DG; Right column: $P^2$-DG.}
   \label{figbwmesh}
   \end{center}
   \end{figure}%
   
}
\end{exam}

 %final section
\section{Conclusions}
\label{seccon}
\setcounter{equation}{0}
\setcounter{figure}{0}
\setcounter{table}{0}

A quasi-conservative DG-ALE method for multi-component flows has been proposed in the previous sections. The non-oscillatory kinetic flux is adopted to compute the numerical flux, which avoids the need to construct a Riemann solver. The main difference between the current method and other DG-ALE methods
is that we use here a predictor-corrector strategy to define the grid velocity, i.e., a Lagrangian step is used to predict the new mesh
and then the MMPDE moving mesh method is employed to improve the quality of the Lagrangian mesh. In this way,
the method can take the advantages of both the Lagrangian meshing and the MMPDE moving mesh method and
keep good quality of the mesh while tracking discontinuities such as shocks and material interfaces in the flow field.
The quasi-conservative DG method \cite{luo2020} is adopted to solve the extended Euler equations of compressible multi-component flows
on moving meshes. Numerical results in one and two dimensions have demonstrated the high-order accuracy of the method
for smooth problems and its ability to capture shocks and material interfaces. Particularly the results demonstrate
that the incorporation of the Lagrangian meshing with the MMPDE moving mesh method works well to concentrate mesh points
in regions of shocks and material interfaces. In the future, we plan to extend the method to
compressible multi-component flows with more general equations of state such as the Mie-Gr\"uneisen equation of state.

\section*{Acknowledgements}
The research is partly supported by NSAF grant 12071392, Science Challenge Project, No. TZ2016002 and National Key Project (GJXM92579).

%\end{CJK*}
\end{document}